# A Stochastic Flow for the Stochastic Allen-Cahn Equation with Multiplicative Noise

Salvador Esquivel & Hendrik Weber

Institut für Analysis und Numerik
Universität Münster

*Email:* salvador.esquivel@uni-muenster.de
*Email:* hendrik.weber@uni-muenster.de

**Abstract:** We establish the existence of a stochastic flow on $L^{\infty}(\mathbb{T})$ for the stochastic Allen-Cahn equation with multiplicative noise

$$(\partial_t - \partial_x^2)\, u = u - u^3 + \sigma(u)\, \xi \quad \text{on} \quad \mathbb{R}_+ \times \mathbb{T},$$

where $\xi$ is space-time white noise and $\sigma\colon \mathbb{R} \to \mathbb{R}$ is sufficiently smooth, bounded, and has bounded derivatives. Our strategy is to obtain pathwise a priori estimates via regularity structures. In fact, we consider a general *singular* multiplicative equation with superlinear damping, driven by noises of parabolic regularity $\alpha - 2$, for $\alpha \in (0, 1)$, which can be lifted to a weakly admissible model. We show that the required estimates hold whenever

$$m > \frac{2-\alpha}{\alpha} \varepsilon_{\alpha}, \quad \text{where} \quad \varepsilon_{\alpha} = 1 - \alpha \left(1 - \frac{2}{3-\alpha}\right) \in (0, 1).$$

Thus the strength of the damping needs to be chosen only as a function of the regularity of the driving noise. Under an additional smoothness assumption on $\sigma$, we show that the stochastic flow is differentiable with respect to its initial condition.

## Table of contents

# 1. Introduction

The existence of stochastic flows for 1-dimensional reaction-diffusion equations

$$(\partial_t - \partial_x^2)\, u = b(u) + \sigma(u)\, \xi \quad \text{on} \quad \mathbb{R}_+ \times \mathbb{T} \tag{1.1}$$

driven by space-time white noise $\xi$ has remained, in most cases, an open problem for a long time. The main difficulty is that probabilistic existence results, even when global in time, come with a null set that depends on the initial condition (see e.g. [DS25, Theorem 4.2.1] and [DZ92, Theorem 7.2]). When the initial conditions form a finite-dimensional space, there are general results known as perfection theorems [AS95], which allow one to obtain a null set independent of the initial condition. These perfection theorems rely on some form of Kolmogorov continuity theorem, and consequently, no analogous results are available in the infinite-dimensional setting where the existence of stochastic flows can fail even for linear equations, as shown in [DZ92, Section 9.1.2]. The existence of a stochastic flow allows one to work with a *random dynamical system* [Arn98] and analyse concepts like ergodic theorems, random attractors, bifurcations, invariant manifolds, etc.

The stochastic Allen-Cahn equation

$$(\partial_t - \partial_x^2)\, u = u - u^3 + \sigma(u)\, \xi \quad \text{on} \quad \mathbb{R}_+ \times \mathbb{T} \tag{1.2}$$

is obtained when the reaction term $b$ is the negative gradient of the double-well quartic potential $F(x) = \frac{1}{4}(x^2 - 1)^2$. The deterministic equation ($\sigma = 0$) is a classical model [AC79] that arises from an $L^2$-gradient flow of the Ginzburg-Landau free energy

$$\mathcal{E}[u] = \int_{\mathbb{T}} \left( \frac{1}{2} |\partial_x u|^2 + \frac{1}{4}(u^2 - 1)^2 \right) \mathrm{d}x$$

to study phase separation and the evolution of its interfaces. The inclusion of the noise term allows us to take into account thermal fluctuations in the system. When the diffusion coefficient $\sigma$ is constant, known as the additive noise case, the existence of the stochastic flow follows immediately from the pathwise well-posedness of the equation and the *coming down from infinity* property (see [Ber22, Section 2.4]).

When the diffusion is not constant, the equation is not classically well-posed due to the low regularity of the space-time white noise, and the equation needs to be interpreted in the Walsh-Itô sense [Wal86]. However, as discussed before, this probabilistic notion of solution is not suitable to prove the existence of the stochastic flow. A pathwise approach to solutions can instead be obtained within the framework of regularity structures [Hai14]. This was done in [HP15] to obtain a version of the Wong-Zakai theorem for solutions to the SPDE (1.1) where, roughly speaking, it is shown that if we consider smooth $\varepsilon$-approximations $\xi_\varepsilon$ to the space-time white noise $\xi$ and classical solutions to the *renormalised* equation

$$(\partial_t - \partial_x^2)\, u_\varepsilon = \bar{b}(u_\varepsilon) + \sigma(u_\varepsilon)\, \xi_\varepsilon - c_\varepsilon\, (\sigma\, \sigma^{(1)})(u_\varepsilon),$$

with fixed continuous initial condition and for a modified drift $\bar{b}$ that depends on the particular approximation but not on the scale $\varepsilon$, see (4.11), then the solutions $u_\varepsilon$ converge to the Walsh-Itô solution of (1.1), with the same initial condition, in $C([0, T] \times \mathbb{T})$ in probability. The renormalisation constants $c_\varepsilon$ diverge at rate $\varepsilon^{-1}$ as $\varepsilon \to 0$, which is related to the fact that $(1+1)$-dimensional space-time white noise $\xi$ is not trace class and there is no Stratonovich notion of solution to (1.1).

As addressed in [HP15, Remark 1.3], although the notion of solutions they use is pathwise and thus yields a global null set independent of the initial conditions, their result does not imply the existence of a stochastic flow. This is due to their reformulation of (1.1) as a fixed-point problem (in a space of *modelled distributions*), which leads only to local-in-time existence of solutions. If one can rule out finite-time blow-up of solutions in this framework, then the existence of the stochastic flow will follow. This is precisely the approach we take in this work. We mention the work [GGH25] where a similar pathwise strategy is used to prove the existence of a random dynamical system associated with a coupled system of an SDE and a Fokker-Planck equation using the theory of rough paths.

The global-in-time existence of solutions in the framework of regularity structures (or alternative pathwise approaches like paracontrolled calculus [GIP15] or the flow equation approach [Duc25]) has seen a lot of work in recent years, [CMW23, CDW26, EW24, BC24, ZZZ22, HZ25, CG25] to mention some. It is important to note that the results of global-in-time existence of solutions in the Walsh–Itô theory do not imply the corresponding global-in-time existence for (pathwise) solutions obtained via the fixed-point problem in the space of modelled distributions.

**Setting and Main Result**

In what follows, we consider a more general multiplicative SPDE with superlinear damping, formally given by

$$(\partial_t - \partial_x^2)\, u = \beta\, u - u\, |u|^{m-1} + \sigma(u)\, \xi \quad \text{on} \quad \mathbb{R}_+ \times \mathbb{T}, \tag{1.3}$$

with $\beta \in \mathbb{R}$, $\sigma \in C_b^N(\mathbb{R})$, bounded and with derivatives up to order $N \in \mathbb{N}$ bounded, a random driving noise $\xi$ which takes values in a parabolic Hölder-Besov space $\mathcal{C}^{\alpha-2}$ of space-time distributions with negative regularity for some $\alpha \in (0,1)$ and which can be lifted into a model $(\Pi, \Gamma)$ in the sense of regularity structures, and $m > 1$ the strength of the damping effect which will depend only on the regularity of the driving noise. Setting $\beta = 1$, $m = 3$ and taking $\xi$ to be $(1+1)$-dimensional space-time white noise, corresponding to $\alpha = 1/2 - \kappa$ for any $0 < \kappa \ll 1$, we recover the stochastic Allen-Cahn equation (1.2).

For a noise with this regularity, by Schauder theory, one expects the solution $u$ to be $\alpha$-Hölder continuous. Consequently, $\sigma(u)$ is well-defined and $\alpha$-Hölder continuous (same for $u\,|u|^{m-1}$). However, since the sum of the regularities of $\sigma(u)$ and $\xi$ is negative, the product "$\sigma(u)\,\xi$" is not classically well-defined as soon as $\alpha < 1$, which makes equation (1.3) singular and in need of renormalisation. Moreover, this equation is subcritical if and only if $\alpha > 0$, and renormalised solutions to it can be described using the theory of regularity structures. We state our main result:

**THEOREM 1.1**. *Fix $\alpha \in (0,1)$ and $\beta \in \mathbb{R}$, set $\varepsilon_\alpha := 1 - \alpha\left(1 - \frac{2}{3-\alpha}\right)$ and let $m \in \left(\frac{2-\alpha}{\alpha}\varepsilon_\alpha, \frac{2-\alpha}{\alpha}\right]$. Choose $\gamma = \gamma(\alpha, m) \in (2-\alpha, 2)$ sufficiently close to $2-\alpha$ such that $\gamma/\alpha \notin \mathbb{N}$, Remark 2.2 holds and*

$$m > \frac{2-\alpha}{\alpha}\left(\gamma - 1 + \frac{2\,\alpha}{\gamma+1}\right).$$

*Assume that $\sigma \in C_b^k(\mathbb{R})$ for $k = \lceil \gamma/\alpha \rceil$. Let $(\Pi, \Gamma)$ be a 1-periodic in space weakly admissible model (Definition 2.1) on the regularity structure described in Section 2.1, $U \in \mathcal{D}^\gamma$ a 1-periodic in space coherent modelled distribution (see (2.16)) whose reconstruction $u := \mathcal{R}U$ solves weakly in $(0,1) \times \mathbb{T}$ the equation*

$$(\partial_t - \partial_x^2)\, u = \beta\, u - u\, |u|^{m-1} + \mathcal{R}(\hat{\sigma}(U)\, \Xi), \tag{1.4}$$

*which is the renormalised version of (1.3). Then*

$$\|u\|_{(t,1]\times\mathbb{T}} \lesssim_{\alpha, |\beta|, m, \gamma, \|\sigma\|_{C_b^k}} \max\left\{t^{-\frac{1}{m-1}}, (1 + |||(\Pi, \Gamma)|||)^\theta\right\} \qquad \forall\, t \in (0,1), \tag{1.5}$$

*where $|||(\Pi, \Gamma)|||$ is defined in (2.14) and $\theta > 0$ is defined as*

$$\theta := \frac{2}{\left(\alpha\, m - (2-\alpha)\left(\gamma - 1 + \frac{2\,\alpha}{\gamma+1}\right)\right)}.$$

*In particular, the implicit constant in (1.5) is independent of the model and of the initial condition. Moreover, suppose that $u$ has initial condition in $L^\infty(\mathbb{T})$, then there exists $T_0 \in (0,1]$ depending on $\alpha, |\beta|, m, \gamma, \|\sigma\|_{C_b^k}, \|u_0\|$ and $|||(\Pi, \Gamma)|||$ such that*

$$\|u\|_{[0,T_0]\times\mathbb{T}} \lesssim (1 + \|u_0\|_{L^\infty(\mathbb{T})}),$$

*and $T_0$ can be chosen uniformly over initial conditions and models such that $\|u_0\|_{L^\infty(\mathbb{T})}, |||(\Pi, \Gamma)|||$ remain bounded.*

Under the slightly stronger condition on the damping $m > (2-\alpha)\,\alpha^{-1}$, we prove an analogous result, see Theorem 3.13, without imposing periodic boundary conditions.

**Remark 1.2.** The statement of Theorem 1.1 is conditional on the existence of a weakly admissible model for the noise under consideration. General results on the construction of such models, such as [CH16, HS24], do not cover the full subcritical regime $\alpha \in (0,1)$ in spatial dimension $d=1$. This is due to the variance blow-up phenomenon, which in general imposes that additional restriction $\alpha > 1/4$. We refer to [Hai25] for recent progress in this direction. However, the restriction to $d=1$ in Theorem 1.1 is only there to simplify the notation, in particular to avoid some vector-valued expressions. The same estimates hold in arbitrary spatial dimension, with constants depending also on the dimension, provided that the equation remains subcritical.

Finally, in the main application considered below (Corollary 1.3), the lift of $(1+1)$-dimensional space-time white noise to a model was also explicitly carried out in [HP15].

As a consequence of Theorem 1.1 and the results from Section 4, in particular Theorem 4.4, we obtain the following result regarding the existence and differentiability of the stochastic flow induced by equation (1.2). Its proof can be found at the end of Section 4.2.

**Corollary 1.3.** *For $\sigma \in C_b^5(\mathbb{R})$ the stochastic Allen-Cahn equation with multiplicative noise (1.2) interpreted in the Itô-Walsh sense generates a stochastic flow on $L^\infty(\mathbb{T})$ in the sense of [DZ92, Section 9.1.2], i.e., there exists*

$$\varphi\colon \{(s,t) \in [0,\infty)\colon s \leqslant t\} \times L^\infty(\mathbb{T}) \times \Omega \to L^\infty(\mathbb{T}),$$

*and a set $\Omega_0 \subset \Omega$ of probability one such that, for every $\omega \in \Omega_0$, every $u_0 \in L^\infty(\mathbb{T})$ and every $0 \leqslant s \leqslant r \leqslant t$, the following hold:*

*i. $q \mapsto \varphi(s,q;u_0,\omega)$, $q \geqslant s$, is the Walsh-Itô solution of (1.2) started from $u_0$ at time $s$;*

*ii. $\varphi(s,s;u_0,\omega) = u_0$;*

*iii. $\varphi(r,t;\varphi(s,r;u_0,\omega),\omega) = \varphi(s,t;u_0,\omega)$;*

*iv. the map $u_0 \mapsto \varphi(s,t,u_0,\omega)$ is continuous.*

*Moreover, if $\sigma \in C_b^6(\mathbb{R})$ then the stochastic flow is almost surely Fréchet differentiable with respect to the initial condition, i.e., for any $h \in L^\infty(\mathbb{T})$, $0 \leqslant s \leqslant t$ the directional derivative*

$$D^h \varphi(s,t,u_0,\omega) = \lim_{\varepsilon \to 0} \frac{1}{\varepsilon}\left(\varphi(s,t,u_0+\varepsilon h,\omega) - \varphi(s,t,u_0,\omega)\right) \in L^\infty(\mathbb{T}),$$

*is a.s. well-defined, and the first variation process $D^h\varphi(s,\cdot,u_0,\cdot)$ is the Itô-Walsh solution to the linearised equation*

$$(\partial_t - \partial_x^2)\, w = (1 - 3\,\varphi(s,\cdot,u_0,\cdot)^2)\, w + \sigma^{(1)}(\varphi(s,\cdot,u_0,\cdot))\, w\,\xi \qquad \text{on} \quad \mathbb{R}_+ \times \mathbb{T},$$

*with initial condition $w(s) = h$.*

**Remark 1.4.** We briefly explain the different regularity assumptions on $\sigma$ in Corollary 1.3. For space-time white noise, the parameters can be chosen so that $k = \lceil \gamma/\alpha \rceil = 4$. This is the number of bounded derivatives of $\sigma$ used for the a priori estimates of Theorem 1.1. The construction of the flow also relies on the local fixed-point theory in spaces of weighted modelled distributions. The local Lipschitz continuity of the map $U \mapsto \hat{\sigma}(U)$ in this argument is obtained by assuming one additional bounded derivative of $\sigma$; see [FH20, Proposition 14.8]. Finally, the differentiability argument applies the same local Lipschitz estimate to $\sigma^{(1)}$ and therefore uses one further bounded derivative. This explains the assumptions $\sigma \in C_b^5$ and $\sigma \in C_b^6$, respectively, in Corollary 1.3.

**Outline of the proof**

We follow the strategy introduced in [MW20b] to show a priori estimates for the Da Prato-Debussche remainder associated with the dynamic $\varphi_3^4$ model. Section 2 is devoted to introducing the regularity structure in which the coherent modelled distribution takes values, and the main results there, Lemma 2.5, provides controls on all the high regularity modelled seminorms of $\hat{\sigma}(U)$ in terms of one high regularity modelled seminorm of $U$ (the one at level one), its generalised derivative (denoted by $u'$, which comes from the coefficient of the polynomial symbol of $U$ in (2.16)), and some seminorms coming from the model $(\Pi, \Gamma)$, which we assume as given data. This analytic control relies on the algebraic identities from Lemma 2.4 for the action of the structure group on $\hat{\sigma}(U)$, which we obtain via a duality argument using the $\star$-product introduced in [BM23] and a Taylor remainder formula. This result generalises [BCMW22, Corollary 2.13] from the rough path setting.

In Section 3.1 the main result is Theorem 3.2, which yields a control on the high regularity modelled seminorms of $U$, and by the results of Section 2 on those of $\hat{\sigma}(U)$, in terms of the $L^\infty$-norm of the solution when restricted to small scales. We expect these estimates to be optimal. The proof of this result follows the standard loop of Schauder estimates, reconstruction, chain rule and a three-point argument used, e.g., in [BOS25, CMW23, EW24]. The precise scale of smallness needed, analogous to the local-in-time existence for a mild formulation up to the parabolic scaling, is given by Definition 3.1. After a very careful algebraic analysis of the local expansions, see Lemmas 2.4 and 2.5, this scale is roughly given by

$$\lambda \sim \|u\|^{-\frac{(m-1)}{2}} \wedge \|u\|^{-\frac{(\gamma-1)}{\alpha}}, \tag{1.6}$$

where $\|u\|$ is the $L^\infty$-norm. The first term comes from the damping term of the equation, while the second comes from the singular product. We can see two different behaviours depending on the strength of the damping:

- If $m > \frac{2-\alpha}{\alpha}$ then, in general, $\|u\|^{-\frac{m-1}{2}} \leqslant \|u\|^{-\frac{\gamma-1}{\alpha}}$. This follows because $\frac{2-\alpha}{\alpha} = 1 + \frac{2((2-\alpha)-1)}{\alpha}$ and the map $x \mapsto 1 + \frac{2(x-1)}{\alpha}$ is increasing, which means that there must exist $\gamma \in (2-\alpha, 2)$ close enough to $2-\alpha$ such that

$$m > \frac{2(\gamma-1)}{\alpha} + 1 \Longleftrightarrow \frac{m-1}{2} > \frac{\gamma-1}{\alpha}. \tag{1.7}$$

  In particular, for this choice of $\gamma \in (2-\alpha, 2)$ the smallness scale is given by $\lambda \sim \|u\|^{-\frac{m-1}{2}}$. In this case, one can apply the exact same argument as [CMW23, BCMW22] to close the global estimate. We only sketch the argument here, since this is not the main result of this work. Consider $(u)_\lambda$ a regularisation of the solution at some scale $\lambda$, which will be chosen depending on a local $L^\infty$-norm of the solution, to get the equation

$$(\partial_t - \partial_x^2)(u)_\lambda = \beta(u)_\lambda - (u)_\lambda |(u)_\lambda|^{m-1} + ((u)_\lambda |(u)_\lambda|^{m-1} - (u|u|^{m-1})_\lambda) + (\mathcal{R}(\hat{\sigma}(U)\,\Xi))_\lambda.$$

  The coercivity of the term $-(u)_\lambda |(u)_\lambda|^{m-1}$ is exploited using [MW20a, Theorem 4.4] which reduces the proof to obtaining an estimate on the other terms on the right hand side in terms of $\|u\|^m$. Control of the commutator term and the removal of the regularisation is achieved by exploiting the small but positive $\alpha$-Hölder regularity of the solution; see Lemma 3.19. For the reconstruction term, one uses the bound

$$\lambda^2 \|\mathcal{R}(\hat{\sigma}(U)\,\Xi)_\lambda - \Pi_x(\hat{\sigma}(U)(x)\,\Xi)_\lambda\| \lesssim \|u\|,$$

  consequence of (3.28), (3.33), and (3.36) and Theorem 3.4, which implies that

$$\|\mathcal{R}(\hat{\sigma}(U)\,\Xi)_\lambda - \Pi_x(\hat{\sigma}(U)(x)\,\Xi)_\lambda\| \lesssim \lambda^{-2}\|u\| \sim \|u\|^{(m-1)+1} = \|u\|^m.$$

  The bounds for the local approximation $\Pi_x(\hat{\sigma}(U)(x)\,\Xi)_\lambda$ follow similarly, and the argument can be closed.

- If $m \leqslant \frac{2-\alpha}{\alpha}$ then for any $\gamma \in (2-\alpha, 2)$ one has that

$$\frac{m-1}{2} \leqslant \frac{(2-\alpha)-1}{\alpha} < \frac{\gamma-1}{\alpha} \tag{1.8}$$

  and therefore the scale in (1.6) is determined by $\lambda \sim \|u\|^{-\frac{(\gamma-1)}{\alpha}}$. If one tries to replicate the argument as in the previous case, one ends up with a bound for the reconstruction term of the form

$$\|\mathcal{R}(\hat{\sigma}(U)\,\Xi)_\lambda - \Pi_x(\hat{\sigma}(U(x))\,\Xi)_\lambda\| \lesssim \lambda^{-2}\,\|u\| \sim \|u\|^{\frac{2(\gamma-1)}{\alpha}+1}, \tag{1.9}$$

  which by our assumptions will result in a power strictly bigger than $m$, which prevents the argument from closing.

The Stochastic Allen-Cahn equation (1.2) falls marginally into the second regime. Indeed, for any $0 < \kappa \ll 1$ and $\alpha = 1/2 - \kappa$ we have that

$$\frac{2-\alpha}{\alpha} = \frac{3/2+\kappa}{1/2-\kappa} > 3 = m.$$

This situation is similar to the one encountered when trying to establish global estimates for the *generalised Parabolic Anderson Model* (gPAM) without damping

$$(\partial_t - \Delta)u = \sigma(u)\,\xi, \tag{1.10}$$

with $\xi \in \mathcal{C}^{-1-\kappa}$ in the first subcritical regime $0 < \kappa < 1/3$. In our notation this corresponds to $\alpha = 1 - \kappa$, while $m = 1$, since linear terms do not create an obstruction for the argument. This problem was recently solved in [CDW26, SZZ26a, SZZ26b] for $0 < \kappa < \sqrt{5} - 2$. We build on the main idea of those works, originally introduced in [CDW26], consisting in rewriting the worst superlinear terms of the equation into transport terms; see (3.17). By Lemma 3.17, these terms do not contribute to the $L^\infty$-estimate.

A key additional difficulty in the present work is that this transport decomposition has to be carried out in the full subcritical regime. This requires a detailed algebraic analysis of the action of the structure group; see Lemma 2.4. It is precisely this analysis that makes possible the refined decomposition of the local expansions in Lemma 2.9. While the full subcritical regime allows arbitrarily many terms in the local expansions, we show that the transport decomposition needs to be applied only to three terms; see (3.43) and (3.62). These are precisely the terms already present in [CDW26]. Although sharper decompositions are possible, as suggested by the refinements in [SZZ26b], the present one is sufficient to cover the cubic damping in the stochastic Allen-Cahn equation.

One of these worst superlinear contributions comes from a quadratic term in the generalised derivative $u'$, which after an application of Taylor's remainder theorem, see Lemma 2.8, is of the form

$$H(x,y) := (u'(x)\,(y-x)_1)^2 \int_0^1 (1-\lambda)\,\sigma^{(2)}(u(x) + \lambda\,u'(x)\,(y-x)_1)\,\mathrm{d}\lambda,$$

see (3.55). This *germ* appears in the change-of-basepoint formula $\langle \mathbf{1}, \hat{\sigma}(U)(y) - \Gamma_{yx}\hat{\sigma}(U)(x)\rangle$, and contributes to one of the modelled seminorms needed to control the singular product in (1.4) through the Reconstruction Theorem. If this quadratic term in $u'$ is estimated directly, one obtains precisely the superlinear bounds that prevent the argument from closing. The core idea is to keep track of one of these generalised derivatives through the reconstruction operator via the explicit multi-scale decomposition (3.50). We then use the identity (3.41) to relate this generalised derivative $u'$ to the genuine derivative $\partial_{e_1}(u_\lambda)$ of the solution regularised at scale $\lambda > 0$. This converts the quadratic term in $u'$ into a transport term, up to a lower-order error. The details of this part of the argument are presented in Lemmas 3.14 and 3.15. This is what allows us to prove that, for every $\alpha \in (0,1)$, there exists $\varepsilon_\alpha \in (0,1)$, continuous in $\alpha$ and typically close to 1, such that every damping strength $m > \varepsilon_\alpha\,(2-\alpha)\,\alpha^{-1}$ is sufficient to obtain estimates that are sublinear in the $L^\infty$-norm of the solution. This final part of the proof of Theorem 1.1 is conducted in Section 3.5.

In Section 4, we turn to the flow generated by the equation. Under stronger smoothness assumptions on the diffusion coefficient $\sigma$, we establish differentiability of the flow with respect to the initial condition. This derivative is then identified with the solution to the linearised equation, a singular linear equation which requires its own renormalisation; see Theorem 4.4. Finally, in Section 4.3, we specialise to the case of $(1+1)$-dimensional space-time white noise and explain how it fits into the regularity structures framework used in this work, relying largely on the construction of the model and the renormalisation results already carried out in [HP15].

**Comparison with [CG25]**

In the recent preprint [CG25], the authors develop an alternative strategy to obtain a priori estimates, using scaling arguments, for a slightly more general class of equations than (1.3), still within the framework of regularity structures. However, for the equations we consider in this work their results are significantly less sharp than ours since the strength of the damping they require for a given regularity of the noise is much higher than that required in our Theorem 1.1. According to [CG25, Assumption 6.1] we denote by $\alpha_{\mathrm{CG}} = 2\,(m-1)^{-1}$ the parameter that makes the deterministic part of the equation scale invariant, and by $\beta_{\mathrm{CG}} = \alpha - 2$ the regularity of the driving noise. Following [CG25, Section 6.3.2] there must exist $\eta \geqslant 0$, a parameter involved in some weighted norms, such that

$$\alpha_{\mathrm{CG}}\,(\eta - 1) < \beta_{\mathrm{CG}} + 2 \Longleftrightarrow \eta < 1 + \frac{\beta_{\mathrm{CG}} + 2}{\alpha_{\mathrm{CG}}} = 1 + \frac{\alpha}{2}\,(m-1), \tag{1.11}$$

(see eq. (6.29) in [CG25]) and, if the diffusion coefficient $\sigma$ is not a polynomial, such that

$$\eta \geqslant \left\lfloor \frac{2}{2+\beta_{\mathrm{CG}}} \right\rfloor = \left\lfloor \frac{2}{\alpha} \right\rfloor. \tag{1.12}$$

The conditions (1.11) and (1.12) require

$$m > 1 + \frac{2}{\alpha}\left( \left\lfloor \frac{2}{\alpha} \right\rfloor - 1 \right) = \mathcal{O}(\alpha^{-2}), \tag{1.13}$$

which is a lower bound on the strength of the damping in terms of the regularity of the driving noise. The condition required in this work is $m > \frac{2-\alpha}{\alpha}\,\varepsilon_\alpha = \mathcal{O}(\alpha^{-1})$ which represents an improvement of a full order of magnitude. In particular, for $(1+1)$ space-time white noise, condition (1.13) translates to $m > 13$, while with our strategy we are able to treat $m = 3$ on the torus (but only $m > 3$ in the full space). We believe that the method in [CG25] is not substantially simpler than ours and the difference in length of our approaches comes mainly from the sharp analysis we perform in Lemma 2.5 to obtain the "level **1** controls all" argument, which is crucial for obtaining the correct exponents.

**Notation:** We consider $\mathbb{N} = \{0, 1\ldots\}$ the set of natural numbers. Given $k = (k_0, k_1) \in \mathbb{N}^2$ we denote by $k! := k_0!\,k_1!$, and for $x = (x_0, x_1) \in \mathbb{R}^2$ we set $x^k := x_0^{k_0}\,x_1^{k_1}$. We work with the parabolic distance $d(x, y) := \max\left\{\sqrt{|y_0 - x_0|}, |y_1 - x_1|\right\}$, where $|\cdot|$ denotes the absolute value on $\mathbb{R}$. We call $x_0$ its time coordinate and $x_1$ its spatial coordinate. Given $x \in \mathbb{R}^2$ and $\lambda > 0$ we denote by $B_\lambda(x) \subset \mathbb{R}^2$ the half-parabolic ball, i.e., $B_\lambda(x) := \{y \in \mathbb{R}^2 \colon d(x, y) \leqslant \lambda, y_0 \leqslant x_0\}$. Given a test function $\psi \in C_c^\infty(\mathbb{R}^2)$ we define $\psi_x^\lambda(y) := \lambda^{-3}\,\psi(\lambda^{-2}\,(y_0 - x_0), \lambda^{-1}\,(y_1 - x_1))$.

## 2. Coherent Modelled Distributions

In this section, we introduce the regularity structure associated with the equation (1.3). In this framework, renormalised solutions to (1.3) are represented by *coherent modelled distributions*. Since solutions to (1.3) are expected to be $\alpha$-Hölder functions, the damping term is not singular, and only the product $\sigma(u)\,\xi$ needs to be renormalised. In practical terms, for us, this means that the regularity structure we need to model the solution is the same as the one for the *gPAM* equation (1.10) in the full subcritical regime.

The main result of this section is Lemma 2.5, which gives estimates for the modelled seminorms associated with these coherent modelled distributions, and which are needed to control the singular product in (1.3) via the Reconstruction Theorem. The proof of this result is based on the algebraic identities in Lemma 2.4.

## 2.1. The regularity structure

**Decorated Trees:** Following [BB26, Section 2] we work with decorated trees $(\tau, \mathfrak{n}, \mathfrak{e})$ where $\tau$ is a non-planar rooted tree with set of nodes $N_\tau$, root $\rho_\tau \in N_\tau$, $\mathfrak{n}: N_\tau \to \mathbb{N}^2$ is the node-decoration map encoding multiplication by polynomials, $E_\tau$ is its set of edges and $\mathfrak{e} = (\mathfrak{t}(\cdot), \mathfrak{p}(\cdot)): E_\tau \to \{\mathcal{I}, \Xi\} \times \mathbb{N}^2$ is the edge-decoration map. The first component $\mathfrak{t}$ denotes the type of the edge, $\mathcal{I}$ denotes an abstract integration map and $\Xi$ an abstract noise. The map $\mathfrak{p}$ encodes derivatives acting on the operator $\mathcal{I}$, and therefore we impose the restriction $\mathfrak{p}(e) \neq 0 \Longrightarrow \mathfrak{t}(e) = \mathcal{I}$ for any edge $e \in E_\tau$.

For decorated trees consisting of only one node, its root $\rho$, we denote them by $\mathbf{1}$ if its node decoration $\mathfrak{n}(\rho) = 0 \in \mathbb{N}^2$, and by $\boldsymbol{X}^k$ if $\mathfrak{n}(\rho) = k$. With some abuse of notation, we denote by $\Xi$ the tree consisting of two nodes, with both node decorations equal to $0 \in \mathbb{N}^2$, and its unique edge given by $\mathfrak{e} = (\mathfrak{t}, \mathfrak{p}) = (\Xi, 0)$. Moreover, for $k \in \mathbb{N}^2$ we consider $\mathcal{I}_k$ as an operator that grafts a tree $\tau$ by its root onto a new root with node decoration 0, and with edge decoration (between the old root and the new root) equal to $(\mathcal{I}, k)$. If $k = 0$ we set $\mathcal{I} := \mathcal{I}_0$. Given two decorated trees $\tau_1, \tau_2$ we define their tree product as the tree obtained by identifying their roots into a single root, setting the node decoration at this new root as the sum of the decorations of the previous roots, i.e., $\mathfrak{n}(\rho_{\tau_1 \tau_2}) = \mathfrak{n}(\rho_{\tau_1}) + \mathfrak{n}(\rho_{\tau_2})$ and preserving all the other decorations. With a slight abuse of notation, we define the total polynomial decoration of a tree as

$$\mathfrak{n}(\tau) := \sum_{v \in N_\tau} \mathfrak{n}(v) \in \mathbb{N}^2. \tag{2.1}$$

We denote by $\mathscr{T}$ the set of trees such that at any node there is at most one outgoing edge (we orient edges away from the root) with noise decoration, i.e., type decoration $\mathfrak{t}$ equal to $\Xi$, and such that every such noise edge terminates at a leaf. Every decorated tree in $\mathscr{T}$ can be represented uniquely, up to the order of its factors, as

$$\tau = \boldsymbol{X}^k \, \Xi^\ell \prod_{i=1}^n \mathcal{I}_{a_i}(\tau_i)^{\beta_i} \tag{2.2}$$

for some $k \in \mathbb{N}^2$, $\ell \in \{0, 1\}$, $n \in \mathbb{N}$, $\{a_i\}_{i=1}^n \subset \mathbb{N}^2$, $\{\beta_i\}_{i=1}^n \subset \mathbb{N} \setminus \{0\}$ and $\{\tau_i\}_{i=1}^n \subset \mathscr{T}$ such that $(\tau_i, a_i) \neq (\tau_j, a_j)$ for $i \neq j$. Here and below, an empty product (and $\Xi^0$) is understood to be $\mathbf{1}$. The symmetry factor is defined recursively via $\Xi! := 1$ and for a tree of the form (2.2) as

$$\tau! := k! \prod_{i=1}^n (\tau_i!)^{\beta_i} \, \beta_i!. \tag{2.3}$$

If the tree is given in a different factorisation

$$\tau = \boldsymbol{X}^k \, \Xi^\ell \prod_{i=1}^n \mathcal{I}_{a_i}(\tau_i)$$

where repetition between trees $\mathcal{I}_{a_i}(\tau_i)$ is allowed, its symmetry factor can be written as

$$\tau! = k! \frac{n!}{\delta((\mathcal{I}_{a_i}(\tau_i))_{i=1}^n)} \prod_{i=1}^n \tau_i! \tag{2.4}$$

where $\delta((\mathcal{I}_{a_i}(\tau_i))_{i=1}^n)$ counts the number of different ordered $n$-tuples $(\mathcal{I}_{a_i}(\tau_i))_{i=1}^n$ that correspond to the same unordered collection, or multi-set, $\{\mathcal{I}_{a_i}(\tau_i)\}_{i=1}^n$. The factor $n!/\delta((\mathcal{I}_{a_i}(\tau_i))_{i=1}^n)$ is the order of the subgroup of permutations of $n$-elements that preserve the ordered collection $(\mathcal{I}_{a_i}(\tau_i))_{i=1}^n$.

Given $\alpha \in (0,1)$, which we consider fixed, we define the homogeneity of a tree by setting

$$|\Xi| := \alpha - 2, \qquad |\boldsymbol{X}^k| := |k| := 2\,k_0 + k_1, \quad k = (k_0, k_1) \in \mathbb{N}^2,$$

and recursively

$$|\tau| := |\boldsymbol{X}^k| + \ell\,|\Xi| + \sum_{i=1}^n \beta_i\,(|\tau_i| + 2 - |a_i|). \tag{2.5}$$

The motivation for this definition is that $\Xi$ is an abstract representation of a noise $\xi$ with parabolic regularity $\alpha - 2$, the operator $\mathcal{I}$ is an abstract representation of inverting the operator $(\partial_t - \partial_x^2)$, which, by Schauder theory, improves parabolic regularity by 2, and each derivative decoration on edges decreases it according to the parabolic scaling.

**Coherence map:** One of the main results in [BCCH20] is that in the framework of regularity structures, solutions to (1.3) are encoded via *coherent* functions $U\colon \mathbb{R}^2 \to \operatorname{span}\{\mathscr{T}\}$ of the form

$$U(x) = \sum_{k \in \mathbb{N}^2} \frac{u_k(x)}{k!}\,\boldsymbol{X}^k + \sum_{\tau \in \mathscr{T}} \frac{\Upsilon^\sigma[\tau](\{u_k(x)\}_{k\in\mathbb{N}^2})}{\tau!}\,\mathcal{I}(\tau), \tag{2.6}$$

where $\Upsilon^\sigma$ is the coherence map associated with the nonlinearity $\sigma$. In (2.6) we implicitly use the convention $\mathcal{I}(\boldsymbol{X}^k) = 0$, see [FH20, Remark 14.26]. More precisely, for each $\tau \in \mathscr{T} \setminus \mathcal{P}$, where we denote by $\mathcal{P} = \{\boldsymbol{X}^k\}_{k\in\mathbb{N}^2}$ the set of purely polynomial trees, the function $\Upsilon^\sigma[\tau]\colon \mathbb{R}^{\mathbb{N}^2} \to \mathbb{R}$ should be thought of as a function of the solution to (1.3) and its derivatives, encoded in some abstract variables $\mathcal{X} := \{\mathcal{X}_k\}_{k\in\mathbb{N}^2}$, and which is defined recursively as $\Upsilon^\sigma[\Xi^\ell](\mathcal{X}) = \delta_{\ell,1}\,\sigma(\mathcal{X}_0)$ and for a tree $\tau \in \mathscr{T} \setminus \mathcal{P}$ of the form (2.2) as

$$\Upsilon^\sigma[\tau] := (\partial^k D_{a_1} \cdots D_{a_n} \Upsilon^\sigma[\Xi^\ell]) \prod_{i=1}^n \Upsilon^\sigma[\tau_i], \tag{2.7}$$

where $D_a$ denotes the derivative with respect to the variable $\mathcal{X}_a$, $\partial^{e_i} := \sum_{a\in\mathbb{N}^2} \mathcal{X}_{a+e_i} D_a$ for the canonical vectors $e_i \in \mathbb{N}^2$ and which is then extended to $\partial^k$ for arbitrary $k \in \mathbb{N}^2$ by composition. For simplicity we remove the $\sigma$-dependence on the coherence and denote it by $\Upsilon = \Upsilon^\sigma$. In the representation (2.6) of the solution, only polynomial trees and trees of the form $\mathcal{I}(\tau)$ such that $\Upsilon[\tau] \neq 0$ are needed. We define

$$\mathcal{T}_\Xi := \{\tau \in \mathscr{T} \setminus \mathcal{P} \colon \Upsilon[\tau] \neq 0\}.$$

As a consequence of the noise appearing multiplicatively in (1.3) and the nonlinearity depending only on the solution itself and not on its derivatives, one can easily show that $\tau \in \mathcal{T}_\Xi$ implies that every vertex that is not the terminal vertex of a noise edge has exactly one outgoing noise edge. Moreover, all non-noise edges are decorated by $(\mathcal{I}, 0)$, i.e., there are no instances of $\mathcal{I}_a$ for $a \neq 0$.

**The regularity structure:** For $n \in \mathbb{N}$ denote by $\mathcal{I}(\mathcal{T}_\Xi)^n := \{\prod_{i=1}^n \mathcal{I}(\tau_i) \colon \{\tau_i\}_{i=1}^n \subset \mathcal{T}_\Xi\}$ and $\mathcal{P}\mathcal{I}(\mathcal{T}_\Xi)^n := \{\boldsymbol{X}^k \rho \colon k \in \mathbb{N}^2, \rho \in \mathcal{I}(\mathcal{T}_\Xi)^n\}$. The basis for our regularity structure is given by the set

$$\mathcal{T} := \mathcal{T}_\Xi \cup \mathcal{T}_{\mathcal{U}}, \qquad \mathcal{T}_{\mathcal{U}} := \bigcup_{n\in\mathbb{N}} \mathcal{P}\mathcal{I}(\mathcal{T}_\Xi)^n.$$

For every $\tau \in \mathcal{T}$ there exists $k \in \mathbb{N}^2$, $\ell \in \{0,1\}, n \in \mathbb{N}$ and $\{\tau_i\}_{i=1}^n \subset \mathcal{T}$ such that

$$\tau = \boldsymbol{X}^k\,\Xi^\ell \prod_{i=1}^n \mathcal{I}(\tau_i). \tag{2.8}$$

Given $\tau \in \mathcal{T}$ we define its number of noises by $\mathfrak{C}(\tau) := \#\{e \in E_\tau : \mathfrak{t}(e) = \Xi\}$. Since all integration edges have derivative decoration zero it can be easily seen using (2.8) that

$$\alpha\, \mathfrak{C}(\tau) - 2\, \delta_{\ell,1} = |\tau| - |\mathfrak{n}(\tau)|, \tag{2.9}$$

where $\ell \in \{0,1\}$ is the indicator of a noise edge et the root in the representation (2.8) for $\tau$.

We define the underlying real vector space of our regularity structure as $T := \mathrm{span}\{\mathcal{T}\}$ and endow it with the grading induced by the homogeneity map (2.5). Since $|\mathcal{I}(\Xi)| = \alpha > 0$, one can see by induction that $|\mathcal{I}(\tau)| > 0$ for all $\tau \in \mathcal{T}$. Therefore for any $\tau \in \mathcal{T}_\Xi$ using the representation (2.8) with $\ell = 1$, we conclude that

$$|\tau| - \min_{i=1,\ldots,n} |\tau_i| = |k| + |\Xi| + \sum_{i=1}^{n} |\mathcal{I}(\tau_i)| - \min_{i=1,\ldots,n} |\tau_i| \geqslant \alpha > 0.$$

Since integration and multiplication with polynomials increase homogeneity, this shows that the index set $A = |\mathcal{T}|$ is bounded from below and locally finite, and the underlying rule subcritical in the sense of [BHZ19]. For any $\beta, \gamma \in \mathbb{R}$ we denote by $\mathcal{T}_\gamma := \{\tau \in \mathcal{T} : |\tau| < \gamma\}$, $\mathcal{T}_{\beta,\gamma} := \{\tau \in \mathcal{T} : \beta \leqslant |\tau| < \gamma\}$ and $T_\gamma := \mathrm{span}\{\mathcal{T}_\gamma\}$, $T_{\beta,\gamma} := \mathrm{span}\{\mathcal{T}_{\beta,\gamma}\}$.

We define the collection of trees $\mathcal{T}^+$ as those trees in $\mathscr{T}$ of the form (2.2) such that $\ell = 0$, $\{\tau_i\}_{i=1}^n \subset \mathcal{T}$ and $|\mathcal{I}_{a_i}(\tau_i)| > 0$ for all $i \in \{1, \ldots, n\}$, and denote by $T^+ = \mathrm{span}\{\mathcal{T}^+\}$ and $\mathcal{T}_\beta^+ := \{\tau \in \mathcal{T}^+ : |\tau| < \beta\}$ for $\beta > 0$. We equip $T^+$ with its standard unital commutative product, with unit $\mathbf{1}$. Its elements may have non-zero derivative decorations, but only on integration edges connected to the root. The construction of the structure group $G$ is standard (see e.g. [FH20, Section 15.3]) and can be encoded as a coaction $\Delta : T \to T \otimes T^+$ via the identity

$$\Gamma \tau = (\mathrm{Id}_T \otimes \gamma) \Delta \tau \tag{2.10}$$

for $\Gamma \in G$ and $\gamma \in \mathrm{Char}(T^+)$, i.e., $\gamma : T^+ \to \mathbb{R}$ which is unital, linear, and multiplicative.

Given a multi-index $\beta \in \mathbb{N}^{\mathcal{T}}$, i.e., $\beta : \mathcal{T} \to \mathbb{N}$ with finitely many non-zero entries, we define

$$\mathfrak{a}(\beta) := \prod_{\rho \in \mathcal{T}} \mathcal{I}(\rho)^{\beta(\rho)}. \tag{2.11}$$

Then (2.8) states that for every tree $\tau \in \mathcal{T}$ there exists a unique $(k, \ell, \beta) \in \mathbb{N}^2 \times \{0,1\} \times \mathbb{N}^{\mathcal{T}}$ such that $\tau = \boldsymbol{X}^k\, \Xi^\ell\, \mathfrak{a}(\beta)$. By subcriticality $\mathcal{T}_\gamma$ is finite for any $\gamma \in \mathbb{R}$ and any $\beta \in \mathbb{N}^{\mathcal{T}_\gamma}$ is a multi-index

**Duality:** The $\star$-product as defined in [BB26, eq. (6)] is dual to the coaction when considering the inner product in $T$, and $T^+$, defined by

$$\langle \tau_1, \tau_2 \rangle := \tau_1!\, \delta_{\tau_1, \tau_2},$$

where $\delta_{\cdot,\cdot}$ is Kronecker's delta. For $\tau_1 \in T^+$ and $\tau_2, \mu \in T$, the duality (see also [BM23, Theorem 4.2]) is given by the identity

$$\langle \tau_1 \star \tau_2, \mu \rangle = \langle \tau_2 \otimes \tau_1, \Delta \mu \rangle, \tag{2.12}$$

where on the right hand side we used the usual extension of the inner products to $T \otimes T^+$. The relevant part of identity (2.12) is that it identifies the correct combinatorial factors from the coproduct $\Delta$. From [BB26, Proposition 2.2] one has the following morphism property between the coherence map and the $\star$-product for $\mu, \tau \in \mathscr{T}$ with $\mu = \boldsymbol{X}^k \prod_{i=1}^n \mathcal{I}_{a_i}(\mu_i)$:

$$\Upsilon\left(\left\{\boldsymbol{X}^k \prod_{i=1}^n \mathcal{I}_{a_i}(\mu_i)\right\} \star \tau\right) = (\partial^k D_{a_1} \cdots D_{a_n} \Upsilon[\tau]) \prod_{i=1}^n \Upsilon[\mu_i], \tag{2.13}$$

where $D_a$ and $\partial^k$ are defined as in (2.7). Although $\star$ can be described in terms of grafting operators, for our purposes (2.12) defines the action of $T^+$ on $T$, with the morphism (2.13) as its relevant property.

**Models and modelled distributions:** We assume that the reader is familiar with the notions of a model and the space of modelled distributions (see e.g. [FH20, Chapter 13]). We will only consider models that satisfy the following weak form of the usual admissibility condition. Every admissible model for the heat kernel is weakly admissible in this sense. In particular, the smooth models and the Itô model used later satisfy the assumptions below.

**DEFINITION 2.1.** *A model $(\Pi, \Gamma)$ is called weakly admissible for the heat operator $(\partial_t - \partial_x^2)$ if the following hold:*

1. *On polynomial symbols it acts as the polynomial model, i.e., for all $k \in \mathbb{N}^2$*

$$(\Pi_x \boldsymbol{X}^k)(y) = (y - x)^k, \quad \Gamma_{yx} \boldsymbol{X}^k = (\boldsymbol{X} + (y - x)\, \mathbf{1})^k.$$

2. *For every $\tau \in \mathcal{T}_0 \subset \mathcal{T}_\Xi$, $x \in \mathbb{R}^2$ and $\psi \in C_c^\infty(\mathbb{R}^2)$ the following holds:*

$$\langle (\partial_t - \partial_x^2)(\Pi_x \mathcal{I}(\tau)), \psi \rangle = \langle \Pi_x \tau, \psi \rangle.$$

It will be useful to work with the following modification of the usual model seminorms which makes it homogeneous on the noise:

$$|||(\Pi; \Gamma)||| := \sup_{\substack{\tau \in \mathcal{T} \\ |\tau| < 0}} \left( [\Pi; \tau] \vee \max_{k \in \mathbb{N}^2} [\Gamma\, \mathcal{I}(\tau); \boldsymbol{X}^k] \right)^{\frac{1}{\mathfrak{C}(\tau)}}, \tag{2.14}$$

where $\mathfrak{C}(\tau)$ is the number of noises in the tree, and

$$[\Pi; \tau] := \sup_{B_\lambda(x) \subset P} |\langle \Pi_x \tau, \psi_x^\lambda \rangle|\, \lambda^{-|\tau|}, \quad [\Gamma\, \mathcal{I}(\tau); \boldsymbol{X}^k] := \sup_{x, y \in P} \frac{|\langle \boldsymbol{X}^k, \Gamma_{yx} \mathcal{I}(\tau) \rangle|}{d(x, y)^{|\mathcal{I}_k(\tau)|}},$$

where $P = (0, 1) \times (-1, 1)$ (or $P = (0, 1) \times \mathbb{T}$ in the periodic case), and $\psi \in C_c^\infty$ is a fixed test function with $\mathrm{supp}(\psi) \subset B_1(0)$, the half-parabolic unit ball, and $\int \psi = 1$. We recall that fixing one test function in the definition of $[\Pi; \tau]$ produces an equivalent seminorm to not fixing the test function; see, e.g., [CZ21, Theorem 12.4]. See Remark 2.2 for the justification of considering only trees with negative homogeneity in (2.14).

**Remark 2.2.** The reason (2.14) only involves trees of negative homogeneity is that we can choose $\gamma$ sufficiently close to $2 - \alpha$ so that $\frac{\gamma}{\alpha} < \lfloor \frac{2}{\alpha} \rfloor$. This is possible since $\frac{2 - \alpha}{\alpha} = \frac{2}{\alpha} - 1 < \lfloor \frac{2}{\alpha} \rfloor$. Let $\tau$ be a tree without a noise decoration at the root and without polynomial decorations. By (2.9), $|\tau\, \Xi| = \alpha\, (\mathfrak{C}(\tau) + 1) - 2$, and therefore

$$0 < |\tau\, \Xi| < \gamma + \alpha - 2 \iff \frac{2}{\alpha} - 1 < \mathfrak{C}(\tau) < \frac{\gamma}{\alpha}.$$

However $\mathfrak{C}(\tau)$ is an integer, whereas the interval on the right-hand side contains no integer by the choice of $\gamma$. Thus no such tree can have homogeneity in $(0, \gamma + \alpha - 2)$. Hence, in every occurrence below where the model norm is used, the corresponding tree has negative homogeneity, and the definition (2.14) is sufficient.

Given a modelled distribution $U \in \mathcal{D}^\gamma$ we consider the usual modelled seminorms given for $\tau \in \mathcal{T}$ by

$$[U; \tau]_D := \sup_{x, y \in D} \frac{|\langle \tau, U(y) - \Gamma_{yx} U(x) \rangle|}{d(x, y)^{\gamma - |\tau|}}. \tag{2.15}$$

## 2.2. Control on the modelled seminorms

Assume we are given a model $(\Pi, \Gamma)$. Fix $\gamma \in (1, 2)$ and a coherent modelled distribution $U_\gamma \in \mathcal{D}^\gamma$, then $U_\gamma$ has the form

$$U_\gamma(x) := u(x)\,\mathbf{1} + u'(x)\,\boldsymbol{X} + \sum_{\tau \in \mathcal{T}_{\gamma-2}} \frac{\Upsilon_x[\tau]}{\tau!}\,\mathcal{I}(\tau), \tag{2.16}$$

where $\boldsymbol{X} := \boldsymbol{X}^{(0,1)}$ represents a spatial polynomial, and $\Upsilon_x[\tau] = \Upsilon[\tau](u(x), u'(x))$; see Lemma 2.7.

It is easy to see that, since $|\mathcal{I}(\tau)| > 0$ for all $\tau \in \mathcal{T}$, then $\mathcal{P} \cup \mathcal{I}(\mathcal{T})$ spans a function-like sector in which the coherent modelled distributions (2.6) take values. Therefore $\hat{\sigma}(U)$ is well-defined as

$$\hat{\sigma}(U_\gamma) := \mathcal{Q}_{<\gamma} \sum_{k \geqslant 0} \frac{\sigma^{(k)}(u(x))}{k!}\,(U_\gamma(x) - u(x)\,\mathbf{1})^k \in \mathcal{D}^\gamma,$$

where $\mathcal{Q}_{<\gamma}$ denotes the projection into $\operatorname{span}\{\mathcal{T}_\gamma\}$; see [FH20, Section 14.2]. By Lemma 2.6, with $\Gamma = \mathrm{Id}$, the coherence condition on $U_\gamma \in \mathcal{D}^\gamma$ constrains the coefficients of $\hat{\sigma}(U_\gamma)$ in the following way:

$$\hat{\sigma}(U_\gamma) = \sum_{\tau \in \mathcal{T}_{0,\gamma}} \frac{\langle \tau, \hat{\sigma}(U_\gamma) \rangle}{\tau!}\,\tau = \sum_{\tau \in \mathcal{T}_{0,\gamma}} \frac{\Upsilon_x[\tau\,\Xi]}{\tau!}\,\tau,$$

and in particular

$$\hat{\sigma}(U_\gamma)\,\Xi = \sum_{\tau \in \mathcal{T}_{0,\gamma}} \frac{\Upsilon_x[\tau\,\Xi]}{\tau!}\,\tau\,\Xi. \tag{2.17}$$

**Remark 2.3.** Although it is possible that $\tau \in \mathcal{T}_{0,\gamma}$ might contain a noise decoration at the root, e.g. $\tau = \Xi\,\mathcal{I}(\Xi)^n$ for large enough $n \in \mathbb{N}$, we are implicitly using that those trees do not contribute to the sum in (2.17) since the coherence map would vanish on the tree $\tau\,\Xi$ in that case (trees with two noise edges decorations outgoing from the same node are not allowed).

The following lemmas give us an explicit description of the action of the structure group on $\hat{\sigma}(U_\gamma)$.

**LEMMA 2.4.** *Let $U_\gamma \in \mathcal{D}^\gamma$ a coherent modelled distribution for $\gamma \in (1, 2)$ and $\tau \in \mathcal{T}_\gamma$ a tree without a noise edge at the root. If $\tau$ contains no polynomial decorations, then we have*

$$\begin{aligned} &\langle \tau, \Gamma_{yx}\,\hat{\sigma}(U_\gamma)(x) \rangle \\ &= \sum_{n \in \mathbb{N}} \frac{1}{n!} \sum_{\substack{\tau_1, \ldots, \tau_n \in \mathcal{T} \\ \sum_{i=1}^n |\mathcal{I}(\tau_i)| < \gamma - |\tau|}} \prod_{i=1}^n \left( \frac{\Upsilon_x[\tau_i]}{\tau_i!}\,\gamma_{yx}(\mathcal{I}(\tau_i)) \right) \\ &\quad \times \{ (\Upsilon[\tau\,\Xi])^{(n)}(u(x)) + \mathbb{1}_{\sum_{i=1}^n |\mathcal{I}(\tau_i)| < \gamma - 1 - |\tau|}\,(\Upsilon[\tau\,\Xi])^{(n+1)}(u(x))\,u'(x)(y-x)_1 \}, \end{aligned}$$

*where $\gamma_{yx}$ denotes the character associated to $\Gamma_{yx}$ via (2.10) and we used that, since $\tau$ contains no polynomial decorations, then $\Upsilon[\tau\,\Xi]$ is a function of only one variable, see Lemma 2.7, and $(\Upsilon[\tau\,\Xi])^{(n)}$ denotes the $n$-th derivative of the function $u \mapsto \Upsilon[\tau\,\Xi](u)$. On the other hand, if $\tau$ contains polynomial decorations, then*

$$\begin{aligned} &\langle \tau, \Gamma_{yx}\,\hat{\sigma}(U_\gamma)(x) \rangle \\ &= \sum_{n \in \mathbb{N}} \frac{1}{n!}\,(\Upsilon_x'[\tau\,\Xi])^{(n)}(u(x)) \sum_{\substack{\tau_0, \ldots, \tau_n \in \mathcal{T} \\ |\mathcal{I}'(\tau_0)| + \sum_{i=1}^n |\mathcal{I}(\tau_i)| \in (0, \gamma - |\tau|)}} \frac{\Upsilon[\tau_0]}{\tau_0!}\,\gamma_{yx}(\mathcal{I}'(\tau_0)) \prod_{i=1}^n \left( \frac{\Upsilon[\tau_i]}{\tau_i!}\,\gamma_{yx}(\mathcal{I}(\tau_i)) \right) \\ &\quad + u'(x) \sum_{n \in \mathbb{N}} \frac{1}{n!}\,(\Upsilon_x'[\tau\,\Xi])^{(n)}(u(x)) \sum_{\substack{\tau_1, \ldots, \tau_n \in \mathcal{T} \\ \sum_{i=1}^n |\mathcal{I}(\tau_i)| < \gamma - |\tau|}} \prod_{i=1}^n \left( \frac{\Upsilon_x[\tau_i]}{\tau_i!}\,\gamma_{yx}(\mathcal{I}(\tau_i)) \right), \end{aligned}$$

*where* $\mathcal{I}' = \mathcal{I}_{(0,1)}$ *denotes an edge decorated by a spatial derivative,* $\Upsilon'$ *is a function of one variable such that* $\Upsilon[\tau\,\Xi](u,u') = u'\,\Upsilon'[\tau\,\Xi](u)$*, see Lemma 2.7, and in particular the derivatives* $(\Upsilon'[\tau\,\Xi])^{(n)}$ *are well-defined.*

With Lemma 2.4 we can conclude the following control on the modelled seminorms of $\hat{\sigma}(U_\gamma)$.

**LEMMA 2.5**. *Let* $\sigma \in C_b^k$ *for* $k = \lceil \gamma/\alpha \rceil$ *and* $U \in \mathcal{D}^\gamma$ *a coherent modelled distribution for* $\gamma \in (1,2)$ *and* $\tau \in \mathcal{T}_\gamma$ *be a tree without a noise edge at the root. Consider the set of multi-indexes* $A_\tau := \{\beta \in \mathbb{N}^{\mathcal{T}_{\gamma-|\tau|-2}} : |\mathfrak{a}(\beta)| < \gamma - |\tau|\}$*, where* $\mathfrak{a}$ *is the map from* (2.11)*, and let* $\partial A_\tau$ *be the set of multi-index as defined in Lemma 2.8. Then, if* $\tau$ *contains no polynomial decorations we have the bound*

$$\begin{aligned}\lambda^{\gamma-|\tau|}\,[\hat{\sigma}(U_\gamma);\tau]_{B_\lambda(z)} \;\lesssim\;& \lambda^{\gamma-|\tau|}\,[U_{\gamma-|\tau|};\mathbf{1}]_{B_\lambda(z)} + \max_{\rho\in\mathfrak{a}(\partial A_\tau)} \{(\lambda\,\|u'\|_{B_\lambda(z)})^{|\mathfrak{n}(\rho)|}\,\lambda^{|\rho|-|\mathfrak{n}(\rho)|}\,[\Gamma\,\rho;\mathbf{1}]_{B_\lambda(z)}\}\\ &+\mathbb{1}_{\gamma-|\tau|>1}\max_{\rho\in\mathfrak{a}(A_\tau)}(\lambda\,\|u'\|)^{\gamma-|\tau|-(|\rho|-|\mathfrak{n}(\rho)|)}\,\lambda^{|\rho|-|\mathfrak{n}(\rho)|}\,[\Gamma\,\rho;\mathbf{1}]_{B_\lambda(z)},\end{aligned}$$

*where* $|\mathfrak{n}(\rho)|$ *denotes the homogeneity of the total polynomial decoration of the tree* $\rho$*; see* (2.1)*. On the other hand, if* $\tau$ *contains polynomial decorations we have the bound*

$$\begin{aligned}\lambda^{\gamma-|\tau|+1}[\hat{\sigma}(U_\gamma);\tau]_{B_\lambda(z)} \;\lesssim\;& \lambda\,\|u'\|_{B_\lambda(z)}(\lambda^{\gamma-|\tau|+1}[U_{\gamma-|\tau|+1};\mathbf{1}]_{B_\lambda(z)})^{\frac{\gamma-|\tau|}{\gamma-|\tau|+1}}\\ &+(\lambda\,\|u'\|_{B_\lambda(z)})^{\gamma-|\tau|+1}\left(1+\max_{\substack{\rho\in\mathcal{T},\mathfrak{n}(\rho)\neq 0\\ |\mathcal{I}(\rho)|+|\tau|\in[\gamma,\gamma+1)}}(\lambda^{(|\mathcal{I}(\rho)|-1)}[\Gamma\mathcal{I}(\rho);\mathbf{1}])^{\gamma-|\tau|}\right)\\ &+\lambda\,\|u'\|_{B_\lambda(z)}\max_{\substack{\rho\in\mathcal{T},\mathfrak{n}(\rho)=0\\ |\mathcal{I}(\rho)|+|\tau|\in[\gamma,\gamma+1)}}\lambda^{|\mathcal{I}(\rho)|}[\Gamma\mathcal{I}(\rho);\mathbf{1}]\\ &+\lambda\,\|u'\|_{B_\lambda(z)}\max_{\rho\in\mathfrak{a}(\partial A_\tau)}\lambda^{|\rho|}[\Gamma\,\rho;\mathbf{1}]_{B_\lambda(z)}\\ &+\max_{\rho\in\mathfrak{a}(A_\tau)}\lambda^{\gamma-|\tau|-|\rho|+1}[U_{\gamma-|\tau|-|\rho|+1};\boldsymbol{X}]_{B_\lambda(z)}\lambda^{|\rho|}[\Gamma\,\rho;\mathbf{1}]_{B_\lambda(z)}.\end{aligned}$$

## 2.3. Proof of Lemma 2.4

The following lemma provides us with a description of the action of the structure group on $\hat{\sigma}(U_\gamma)$ via the duality (2.12) of the $\star$-product.

**LEMMA 2.6**. *Let* $U_\gamma \in \mathcal{D}^\gamma$ *be a coherent modelled distribution for* $\gamma > 0$*. Let* $\tau \in \mathcal{T}_\gamma$ *be a tree without a noise edge at the root and* $\Gamma \in G$ *an element in the structure group, then we have the identity*

$$\langle \tau, \Gamma\,\hat{\sigma}(U_\gamma)\rangle = \sum_{\mu\in\mathcal{T}^+_{\gamma-|\tau|}} \frac{\Upsilon[\mu\star(\tau\,\Xi)]}{\mu!}\,\gamma(\mu),$$

*where* $\gamma$ *is the character associated to* $\Gamma \in G$ *via* (2.10)*.*

**Proof.** For $\tau \in \mathcal{T}$ and $\Gamma \in G$, using that $\Gamma\Xi = \Xi$ and the multiplicativity of $\Gamma$ we obtain the identity

$$\langle \tau, \Gamma\,\hat{\sigma}(U_\gamma)\rangle = \langle \tau\,\Xi, (\Gamma\,\hat{\sigma}(U_\gamma))\,\Xi\rangle = \langle \tau\,\Xi, \Gamma\,\hat{\sigma}(U_\gamma)\,\Gamma\Xi\rangle = \langle \tau\,\Xi, \Gamma(\Xi\,\hat{\sigma}(U_\gamma))\rangle = \langle \mathcal{I}(\tau\,\Xi), \mathcal{I}(\Gamma(\hat{\sigma}(U_\gamma)\,\Xi))\rangle.$$

Even though $\mathcal{I}$ and $\Gamma$ do not commute, we have that $\Gamma\,\mathcal{I}\tau - \mathcal{I}\,\Gamma\,\tau \in \mathrm{span}\{\mathcal{P}\}$, the polynomial sector, (see [FH20, Assumption 14.21]), and since $\mathcal{I}(\Xi\,\tau)$ is not a polynomial we have that

$$\langle \tau, \Gamma\,\hat{\sigma}(U_\gamma)\rangle = \langle \mathcal{I}(\tau\,\Xi), \mathcal{I}\,\Gamma(\hat{\sigma}(U_\gamma)\,\Xi)\rangle = \langle \mathcal{I}(\tau\,\Xi), \Gamma\,\mathcal{I}(\hat{\sigma}(U_\gamma)\,\Xi)\rangle.$$

By coherence $U_{\gamma+\alpha} - \mathcal{I}(\hat{\sigma}(U_\gamma)\,\Xi) \in \mathrm{span}\{\mathcal{P}\}$ and using the general identity from [EW24, Lemma 2.19, eq. (2.23)] we conclude

$$\langle \mathcal{I}(\tau\,\Xi), \Gamma\,\mathcal{I}(\hat{\sigma}(U_\gamma)\,\Xi)\rangle = \langle \mathcal{I}(\tau\,\Xi), \Gamma\,U_{\gamma+\alpha}\rangle = \sum_{\substack{\mu\in\mathcal{T}^+ \\ |\mathcal{I}(\mu\star(\Xi\tau))|<\gamma+\alpha}} \frac{\Upsilon[\mu\star(\tau\,\Xi)]}{\mu!}\,\gamma(\mu).$$

At last, observe that $|\mathcal{I}(\mu\star(\Xi\,\tau))| < \gamma+\alpha \Longleftrightarrow |\mu| < \gamma - |\tau|$. □

**LEMMA 2.7.** *Let $\gamma\in(1,2)$, $\sigma\in C_b^k$ for $k=\lceil \gamma/\alpha \rceil$ and $\tau\in\mathcal{T}_\gamma$ a tree without a noise edge at the root:*

1. *If $\tau$ contains no polynomial decorations, then $\Upsilon[\tau\,\Xi]$ is a function of only $\mathcal{X}_0$ which we denote by $\Upsilon[\tau\,\Xi](u)$. More specifically, $\Upsilon[\tau\,\Xi](u)$ is a function of $\{\sigma^{(i)}(u)\}_{i=0}^{\mathfrak{C}(\tau)}$, where $\mathfrak{C}(\tau)$ is the number of noises in $\tau$.*

2. *If $\tau$ has polynomial decorations, then $\mathfrak{n}(\tau) = e_{(0,1)}$ and $\Upsilon[\tau\,\Xi](u,u')$ is a linear function of $u'$. In this case, we denote by $\Upsilon'[\tau\,\Xi] := D_{(0,1)}\,\Upsilon[\tau\,\Xi]$, where $D_{(0,1)}$ is the derivative with respect to the abstract variable $\mathcal{X}_{(0,1)}$ as in (2.7), which leads to the identity*

$$\Upsilon[\tau\,\Xi](u,u') = u'\,\Upsilon'[\tau\,\Xi](u).$$

*Moreover, $\Upsilon'[\tau\,\Xi]$ is a function of $\{\sigma^{(i)}(u)\}_{i=0}^{1+\mathfrak{C}(\tau)}$.*

**Proof.** The first claim follows since $\Upsilon[\Xi](\{\mathcal{X}_k\}_{k\in\mathbb{N}^2}) = \sigma(\mathcal{X}_0)$ only depends on $\mathcal{X}_0$ and by a simple induction using (2.7) since dependence on higher order terms $\mathcal{X}_k$ can only be introduced by the derivative $\partial^k$ in (2.7) which acts when there is a polynomial decoration. That $\Upsilon[\tau\,\Xi]$ is a function of $\sigma$ and its derivatives follows immediately by (2.7) and the highest order derivative $\sigma^{(\mathfrak{C}(\tau))}$ is obtained when $\tau$ is of the form $\mathcal{I}(\Xi)^{\mathfrak{C}(\tau)}$. Observe that by (2.9) with $\ell=0$ since $\tau$ does not contain a noise at the root, we have that $\mathfrak{C}(\tau) = |\tau|/\alpha < \gamma/\alpha \leqslant k$ and since $\mathfrak{C}(\tau)\in\mathbb{N}$ then $\mathfrak{C}(\tau)\leqslant k-1$ which makes all the coefficients well-defined and bounded.

For the second claim, first let $\tilde{\tau}$ denotes the same trees $\tau$ but with the polynomial decorations removed, then $\tilde{\tau}\in\mathcal{T}$. It can be easily seen via an induction using (2.2) with $k=0$ and $\ell=0$ that $0<\alpha\leqslant|\tilde{\tau}|$, and since $|\tau| = |\tilde{\tau}| + |\mathfrak{n}(\tau)| < \gamma < 2$ by definition of $\mathcal{T}_\gamma$ we conclude that $|\mathfrak{n}(\tau)|<2$ and in particular $\mathfrak{n}(\tau)\in\{0, e_{(0,1)}\}$, i.e., $\tau$ contains at most one spatial polynomial decoration $\boldsymbol{X}$. We proceed inductively with the base case following from the identity $\Upsilon[\Xi\,\boldsymbol{X}] = \sigma^{(1)}(u)\,u'$. By (2.8) it is enough to consider $\tau$ of the form $\tau = \boldsymbol{X}^k \prod_{i\in I}\mathcal{I}(\tau_i)$, then we have the following cases:

i. If $k=0$ then $\exists!\ i_0\in I$ such that $\tau_i$ contains a polynomial decoration, and by induction $\Upsilon[\tau_i]$ is a linear function of $\mathcal{X}_{(0,1)}$, while for all $i\in I\setminus\{i_0\}$ the map $\Upsilon[\tau_i]$ does not depend on $u'$ by the previous part and the result now follows from (2.7) since

$$\Upsilon[\tau\,\Xi](\{\mathcal{X}_k\}_{k\in\mathbb{N}^2}) = (D^{|I|}\sigma)\prod_{i\in I}\Upsilon[\tau_i]((\{\mathcal{X}_k\}_{k\in\mathbb{N}^2})) = \mathcal{X}_{(0,1)}\left(\Upsilon'[\tau_{i_0}]\prod_{i\in I\setminus\{i_0\}}\Upsilon[\tau_i]\,\sigma^{(|I|)}\right)(\mathcal{X}_0).$$

ii. If $k=(0,1)$ then $\Upsilon[\tau_i]$ is not a function of $u'$ for all $i\in I$ and therefore

$$\Upsilon[\tau\,\Xi](\{\mathcal{X}_k\}_{k\in\mathbb{N}^2}) = (\partial^{(0,1)}D^{|I|}\sigma)\prod_{i\in I}\Upsilon[\tau_i](\{\mathcal{X}_k\}_{k\in\mathbb{N}^2}) = \mathcal{X}_{(0,1)}\left(\prod_{i\in I}\Upsilon[\tau_i]\,\sigma^{(|I|+1)}\right)(\mathcal{X}_0).$$

Again, the claim of $\Upsilon'[\tau\,\Xi]$ being a function of $\sigma$ and its derivatives follows analogously with the observation that the highest order derivative $\sigma^{(1+\mathfrak{C}(\tau))}$ is in this case obtained when $\tau$ is of the form $\boldsymbol{X}\,\mathcal{I}(\Xi)^n$, in which case $n=\mathfrak{C}(\tau)$. At last, by (2.9) we have $\alpha\,\mathfrak{C}(\tau) = |\tau|-1 < \gamma-1$ and, since $\alpha\in(0,1)$, then $1+\mathfrak{C}(\tau) < \gamma/\alpha \leqslant k$ and $1+\mathfrak{C}(\tau)\leqslant k-1$, which makes all the coefficients $\Upsilon'$ well-defined and bounded. □

Now we can prove Lemma 2.4.

**Proof of Lemma 2.4.** First consider $\tau$ that has no polynomial decorations. We use the representation of Lemma 2.6. Given $\mu = \boldsymbol{X}^k \prod_{i\in I} \mathcal{I}_{a_i}(\tau_i) \in \mathcal{T}^+$ by the morphism property of $\Upsilon$ and the $\star$-product (2.13) we have

$$\Upsilon[\mu \star (\Xi\,\tau)] = \prod_{i\in I} \Upsilon[\tau_i] \left(\partial^k \prod_{i\in I} D_{a_i} \Upsilon[\Xi\,\tau]\right). \tag{2.18}$$

By Lemma 2.7 we have that $\Upsilon[\Xi\,\tau](\mathcal{X}_0)$ and therefore $D_a \Upsilon[\Xi\,\tau] = 0$ for all $a \neq 0$, which allows us to restrict ourselves to $\mu \in \mathcal{T}^+$ of the form $\mu = \boldsymbol{X}^k \prod_{i\in I} \mathcal{I}(\tau_i)$. We then obtain

$$\begin{aligned}
&\langle \tau, \Gamma\,\hat{\sigma}(U_\gamma)\rangle \\
= &\sum_{\mu \in \mathcal{T}^+_{\gamma - |\tau|}} \frac{\Upsilon[\mu \star (\Xi\,\tau)]}{\mu!} \gamma(\mu) \\
= &\sum_{k\in\mathbb{N}^2} \sum_{\mu = \prod_{i\in I}\mathcal{I}(\tau_i) \in \mathcal{T}^+} \frac{\Upsilon[(\boldsymbol{X}^k \prod_{i\in I} \mathcal{I}(\tau_i))]}{(\boldsymbol{X}^k \prod_{i\in I}\mathcal{I}(\tau_i))!} \gamma\left(\boldsymbol{X}^k \prod_{i\in I} \mathcal{I}(\tau_i)\right) \mathbb{1}_{|k| + \sum_{i=1}^n |\mathcal{I}(\tau_i)| < \gamma - |\tau|} \\
= &\sum_{k\in\mathbb{N}^2} \sum_{\substack{\mu = \prod_{i\in I}\mathcal{I}(\tau_i) \in \mathcal{T}^+ \\ |\mu| < \gamma - |\tau| - |k|}} \delta((\mathcal{I}(\tau_i))_{i\in I}) \frac{\prod_{i\in I} \Upsilon[\tau_i]\,(\partial^k D^{|I|} \Upsilon[\Xi\,\tau])}{k! \prod_{i\in I} \tau_i!\,|I|!} \gamma(\boldsymbol{X}^k) \prod_{i\in I} \gamma(\mathcal{I}(\tau_i)) \\
= &\sum_{k\in\mathbb{N}^2} \sum_{n\in\mathbb{N}} \sum_{\tau_1,\ldots,\tau_n \in \mathcal{T}^+} \frac{\gamma(\boldsymbol{X}^k)}{k!} \frac{\partial^k D^n \Upsilon[\Xi\,\tau]}{n!} \prod_{i=1}^n \left(\frac{\Upsilon[\tau_i]}{\tau_i!} \gamma(\mathcal{I}(\tau_i))\right) \mathbb{1}_{|k| + \sum_{i=1}^n |\mathcal{I}(\tau_i)| < \gamma - |\tau|} \\
= &\sum_{n\in\mathbb{N}} \frac{1}{n!} \sum_{\tau_1,\ldots,\tau_n\in\mathcal{T}^+} \left(\sum_{k\in\mathbb{N}^2} \frac{\partial^k D^n \Upsilon[\Xi\,\tau]}{k!} \gamma(\boldsymbol{X})^k \mathbb{1}_{|k| + \sum_{i=1}^n |\mathcal{I}(\tau_i)| < \gamma - |\tau|}\right) \prod_{i=1}^n \left(\frac{\Upsilon[\tau_i]}{\tau_i!} \gamma(\mathcal{I}(\tau_i))\right)
\end{aligned}$$

where $\delta((\mathcal{I}(\tau_i))_{i\in I})$ is as in (2.4). Since the condition $|k| + \sum_{i=1}^n |\mathcal{I}(\tau_i)| < \gamma - |\tau|$ restricts the sum to

$$\begin{aligned}
&\sum_{k\in\mathbb{N}^2} \frac{\partial^k D^n \Upsilon_x[\Xi\,\tau]}{k!} \gamma_{yx}(\boldsymbol{X})^k \mathbb{1}_{|k| + \sum_{i=1}^n |\mathcal{I}(\tau_i)| < \gamma - |\tau|} \\
= &\ ((\Upsilon[\Xi\,\tau])^{(n)}(u(x)) + \mathbb{1}_{\sum_{i=1}^n (|\mathcal{I}(\tau_i)|) < \gamma - 1 - |\tau|} (\Upsilon[\Xi\,\tau])^{(n+1)}(u(x))\, u'(x)(y-x)_1)
\end{aligned}$$

we conclude the result for the case without polynomial decorations.

For the second part consider $\tau$ with polynomial decorations. By subcriticality, and since $\tau$ has no noise edge at the root, one can easily see that $|\tau| \geqslant |\boldsymbol{X}| = 1$ and therefore the restriction $|\mu| < \gamma - |\tau| < 1$ in the sum from Lemma 2.6 implies that $\mu \in \mathcal{T}^+$ with polynomial decoration at the root cannot appear, and therefore

$$\begin{aligned}
\langle \tau, \Gamma\,\hat{\sigma}(U)\rangle &= \sum_{\substack{\mu\in\mathcal{T}^+ \\ |\mu| < \gamma - |\tau|}} \frac{\Upsilon[\mu \star (\tau\,\Xi)]}{\mu!} \gamma(\mu) \\
&= \sum_{\substack{\mu = \prod_{i\in I} \mathcal{I}_{\rho_i}(\tau_i)^{\beta_i} \in \mathcal{T}^+ \\ |\mu| < \gamma - |\tau|}} \frac{\Upsilon[(\prod_{i\in I} \mathcal{I}_{\rho_i}(\tau_i)^{\beta_i}) \star (\tau\,\Xi)]}{\mu!} \gamma(\mu) \\
&= \sum_{\substack{\mu = \prod_{i\in I} \mathcal{I}_{\rho_i}(\tau_i)^{\beta_i} \in \mathcal{T}^+ \\ |\mu| < \gamma - |\tau|}} \frac{\prod_{i\in I} \Upsilon[\tau_i]\,(\prod_{i\in I} (D^{\rho_i})^{\beta_i} \Upsilon[\tau\,\Xi])}{(\prod_{i\in I} \mathcal{I}_{\rho_i}(\tau_i)^{\beta_i})!} \prod_{i\in I} \gamma(\mathcal{I}_{\rho_i}(\tau_i)).
\end{aligned}$$

By Lemma 2.7 $\Upsilon[\tau\,\Xi] = \Upsilon'[\tau\,\Xi]\,\mathcal{X}_{e_1}$ and therefore $\prod_{i\in I} (D^{\rho_i})^{\beta_i} \Upsilon[\tau\,\Xi] \neq 0$ iff $\sum_{i\in I} \beta_i\,\rho_i \in \{0, e_1\}$, in which case we have that

$$\prod_{i\in I} (D^{\rho_i})^{\beta_i} \Upsilon[\tau\,\Xi] = \begin{cases} (\Upsilon'[\tau\,\Xi])^{(\sum_{i\in I}\beta_i)} \mathcal{X}_{e_1} & \text{if } \sum_i \beta_i\,\rho_i = 0 \\ (\Upsilon'[\tau\,\Xi])^{\sum_{i\in I}\beta_i - 1} & \text{if } \sum_i \beta_i\,\rho_i = e_1 \end{cases}.$$

This leads us to rearrange the sum as

$$
\begin{aligned}
&\langle \tau, \Gamma\, \hat{\sigma}(U)\rangle \\
=& \sum_{\substack{\mu=\prod_{i\in I}\mathcal{I}(\tau_i)^{\beta_i}\in\mathcal{T}^+\\|\mu|<\gamma-|\tau|}} \frac{(\Upsilon'[\tau\,\Xi])^{(\sum_{i\in I}\beta_i)}}{(\prod_{i\in I}\mathcal{I}(\tau_i)^{\beta_i})!}\mathcal{X}_{e_1}\prod_{i\in I}\left(\frac{\Upsilon[\tau_i]}{\tau_i!}\gamma(\mathcal{I}(\tau_i))\right)\\
&+\sum_{\substack{\mu=\mathcal{I}'(\tau_{i_0})\prod_{i\in I\setminus\{i_0\}}\mathcal{I}(\tau_i)^{\beta_i}\\|\mu|<\gamma-|\tau|}} \frac{(\Upsilon'[\tau\,\Xi])^{\sum_{i\in I}\beta_i-1}}{(\mathcal{I}'(\tau_i)\prod_{i\in I\setminus\{i_0\}}\mathcal{I}(\tau_i)^{\beta_i})!}\Upsilon[\tau_{i_0}]\,\gamma(\mathcal{I}'(\tau_{i_0}))\prod_{i\in I\setminus\{i_0\}}(\Upsilon[\tau_i]\,\gamma(\mathcal{I}(\tau_i)))\\
=& \sum_{\substack{\mu=\prod_{i\in I}\mathcal{I}(\tau_i)^{\beta_i}\in\mathcal{T}^+\\|\mu|<\gamma-|\tau|}} \frac{(\Upsilon'[\tau\,\Xi])^{(\sum_{i\in I}\beta_i)}}{\prod_{i\in I}\tau_i!\,\beta_i!}\mathcal{X}_{e_1}\prod_{i\in I}\left(\frac{\Upsilon[\tau_i]}{\tau_i!}\gamma(\mathcal{I}(\tau_i))\right)\\
&+\sum_{\substack{\mu=\mathcal{I}'(\tau_{i_0})\prod_{i\in I\setminus\{i_0\}}\mathcal{I}(\tau_i)^{\beta_i}\\|\mu|<\gamma-|\tau|}} \frac{(\Upsilon'[\tau\,\Xi])^{\sum_{i\in I}\beta_i-1}}{\tau_{i_0}!\prod_{i\in I\setminus\{i_0\}}\tau_i!\,\beta_i!}\frac{\Upsilon[\tau_{i_0}]}{\tau_{i_0}!}\gamma(\mathcal{I}'(\tau_{i_0}))\prod_{i\in I\setminus\{i_0\}}\left(\frac{\Upsilon[\tau_i]}{\tau_i!}\gamma(\mathcal{I}(\tau_i))\right)\\
=& \mathcal{X}_{e_1}\sum_{n\in\mathbb{N}}\frac{1}{n!}(\Upsilon'[\tau\,\Xi])^{(n)}\left(\sum_{\substack{\tau_1,\ldots,\tau_n\in\mathcal{T}\\\sum_{i=1}^n|\mathcal{I}(\tau_i)|<\gamma-|\tau|}}\prod_{i=1}^n\left(\frac{\Upsilon[\tau_i]}{\tau_i!}\gamma(\mathcal{I}(\tau_i))\right)\right)\\
&+\sum_{n\in\mathbb{N}}\frac{1}{(n-1)!}(\Upsilon'[\tau\Xi])^{(n-1)}(u(x))\sum_{\substack{\tau_1,\ldots,\tau_n\in\mathcal{T}\\|\mathcal{I}'(\tau_1)|+\sum_{i=2}^n|\mathcal{I}(\tau_i)|\in(0,\gamma-|\tau|)}}\frac{\Upsilon[\tau_1]}{\tau_1!}\gamma(\mathcal{I}'(\tau_1))\prod_{i=2}^n\left(\frac{\Upsilon[\tau_i]}{\tau_i!}\gamma(\mathcal{I}(\tau_i))\right)\\
=& \mathcal{X}_{e_1}\sum_{n\in\mathbb{N}}\frac{1}{n!}(\Upsilon'[\tau\,\Xi])^{(n)}\sum_{\substack{\tau_1,\ldots,\tau_n\in\mathcal{T}\\\sum_{i=1}^n|\mathcal{I}(\tau_i)|<\gamma-|\tau|}}\prod_{i=1}^n\left(\frac{\Upsilon[\tau_i]}{\tau_i!}\gamma(\mathcal{I}(\tau_i))\right)\\
&+\sum_{n\in\mathbb{N}}\frac{1}{n!}(\Upsilon'[\tau\Xi])^{(n)}\sum_{\substack{\tau_0,\ldots,\tau_n\in\mathcal{T}\\|\mathcal{I}'(\tau_0)|+\sum_{i=1}^n|\mathcal{I}(\tau_i)|\in(0,\gamma-|\tau|)}}\frac{\Upsilon[\tau_0]}{\tau_0!}\gamma(\mathcal{I}'(\tau_0))\prod_{i=1}^n\left(\frac{\Upsilon[\tau_i]}{\tau_i!}\gamma(\mathcal{I}(\tau_i))\right),
\end{aligned}
$$

where in the second sum we used the identity for the symmetry factors

$$\left(\mathcal{I}'(\tau_{i_0})\prod_{i\in I\setminus\{i_0\}}\mathcal{I}(\tau_i)^{\beta_i}\right)! = \tau_{i_0}!\prod_{i\in I\setminus\{i_0\}}\tau_i!\,\beta_i! = \tau_{i_0}!\left(\prod_{i=1}^{n-1}\mathcal{I}(\tau_i)\right)! = \tau_{i_0}!\,\frac{(n-1)!}{\delta(\tau_i)_{i\in I}}\prod_{i=1}^{n-1}\tau_i!$$

for $n=\sum_{i\in I}\beta_i$, which allowed us to re-index the sum with the factor $(\sum_{i\in I}\beta_i-1)!=(n-1)!$. ☐

## 2.4. Proof of Lemma 2.5

The following decomposition of Lemma 2.9 will be useful in the proof of Lemma 2.5 and, later, for the proof of Lemma 3.15. We will require to apply a form of Taylor's remainder theorem to a function $f\colon\mathbb{R}^A\to\mathbb{R}$ with $A$ a subset of trees of negative homogeneity, and therefore the identification of trees and multi-indices from (2.11) will be useful. To make sense of Lemma 2.8 with $\mathbb{N}^A$ instead of $\mathbb{N}^d$, we consider the order induced on $A$ by the homogeneity map.

**Lemma 2.8**. ([Hai14, Proposition A.1]) *Let $A\subset\mathbb{N}^d$ be such that $k\in A\Rightarrow k_<\subset A$ and define $\partial A=\{k\notin A\colon k-e_{\mathfrak{m}(k)}\in A\}$ with $\mathfrak{m}(k):=\min\{i\in\{1,\ldots,d\}\colon k_i\neq 0\}$. Then the identity*

$$f(x)=\sum_{k\in A}\frac{D^k f(0)}{k!}x^k+\sum_{k\in\partial A}\int_{\mathbb{R}^d}D^k f(y)\,\mathcal{Q}^k(x,\mathrm{d}y)$$

*holds for every $f\in C^{r_A}(\mathbb{R}^d)$ where $r_A:=\max_{k\in\partial A}\sum_{i=1}^d k(i)$.*

**LEMMA 2.9**. *Let $\gamma\in(1,2)$ and assume $\sigma\in C_b^k$ for $k=\lceil\gamma/\alpha\rceil$. Let $U_\gamma\in\mathcal{D}^\gamma$ a coherent modelled distribution and $\tau\in\mathcal{T}_\gamma$ a tree without a noise edge at the root. Consider the germ $H^\tau:\mathbb{R}^2\times\mathbb{R}^2\to\mathbb{R}$ defined by*

$$H^\tau(x,y):=\langle\tau,\hat{\sigma}(U_\gamma)(y)-\Gamma_{yx}\,\hat{\sigma}(U_\gamma)(x)\rangle.$$

*If $\tau$ has no polynomial decorations, then*

$$H^\tau(x,y)=H_1^\tau(x,y)+\sum_{\beta\in\partial A_\tau}H_\beta^\tau(x,y)+\sum_{n\in\mathbb{N}}\frac{1}{n!}\sum_{\substack{\tau_1,\ldots,\tau_n\in\mathcal{T}\\ \sum_{i=1}^n|\mathcal{I}(\tau_i)|<\gamma-|\tau|}}H^\tau_{(\tau_i)_{i=1}^n}(x,y),\tag{2.19}$$

*where*

$$\begin{aligned}
H_1^\tau(x,y) &:= \Upsilon[\tau\,\Xi](u(y))-\Upsilon[\tau\,\Xi](u(y)-\langle\mathbf{1},U_{\gamma-|\tau|}(y)-\Gamma_{yx}\,U_{\gamma-|\tau|}(x)\rangle),\\
H_\beta^\tau(x,y) &:= \int(\Upsilon[\tau\,\Xi])^{(\sum_\rho\beta(\rho))}(u(x)+\mathbb{1}_{\gamma-|\tau|>1}\,u'(x)\,(y-x)_1+\boldsymbol{U}_{\gamma-|\tau|;x,y}(w))\,\mathcal{Q}^\beta(\underline{1},\mathrm{d}w)\\
&\qquad\times\prod_{\rho\in\mathcal{T}}\left(\frac{\Upsilon_x[\rho]}{\rho!}\,\gamma_{yx}(\mathcal{I}(\rho))\right)^{\beta(\rho)},\\
H^\tau_{(\tau_i)_{i=1}^n}(x,y) &:= \{(\Upsilon[\tau\,\Xi])^{(n)}(u(x)+\mathbb{1}_{\gamma-|\tau|>1}u'(x)\,(y-x)_1)-(\Upsilon[\tau\,\Xi])^{(n)}(u(x))\\
&\qquad-\mathbb{1}_{\sum_{i=1}^n|\mathcal{I}(\tau_i)|<\gamma-|\tau|-1}\,(\Upsilon[\tau\,\Xi])^{(n+1)}(u(x))\,u'(x)\,(y-x)_1\}\\
&\qquad\times\prod_{i=1}^n\left(\frac{\Upsilon_x[\tau_i]}{\tau_i!}\,\gamma_{yx}(\mathcal{I}(\tau_i))\right),
\end{aligned}\tag{2.20}$$

*and all the derivatives of $\Upsilon[\tau\,\Xi]$ appearing are bounded.*

*If $\tau$ has polynomial decorations, then*

$$\begin{aligned}
H^\tau(x,y) &= u'(y)\left(H_1^\tau(x,y)+H_2^\tau(x,y)+\sum_{\beta\in\partial A_\tau}H_\beta^\tau(x,y)\right)\\
&\quad+\sum_{n\in\mathbb{N}}\frac{1}{n!}\sum_{\substack{\tau_1,\ldots,\tau_n\in\mathcal{T}\\ \sum_{i=1}^n|\mathcal{I}(\tau_i)|<\gamma-|\tau|}}H^\tau_{(\tau_i)_{i=1}^n}(x,y),
\end{aligned}\tag{2.21}$$

*where*

$$\begin{aligned}
H_1^\tau(x,y) &:= \Upsilon'[\tau\,\Xi](u(y))-\Upsilon'[\tau\,\Xi](u(y)-\langle\mathbf{1},U_{\gamma-|\tau|+1}(y)-\Gamma_{yx}\,U_{\gamma-|\tau|+1}(x)\rangle),\\
H_2^\tau(x,y) &:= \Upsilon'[\tau\,\Xi](u(y)-\langle\mathbf{1},U_{\gamma-|\tau|+1}(y)-\Gamma_{yx}\,U_{\gamma-|\tau|+1}(x)\rangle)\\
&\quad-\Upsilon'[\tau\,\Xi](u(y)-\langle\mathbf{1},U_{\gamma-|\tau|}(y)-\Gamma_{yx}\,U_{\gamma-|\tau|}(x)\rangle),\\
H_\beta^\tau(x,y) &:= \int(\Upsilon'[\tau\,\Xi])^{(\sum_\rho\beta(\rho))}(u(x)+\boldsymbol{U}_{\gamma-|\tau|;x,y}(w))\prod_{\rho\in\mathcal{T}}\left(\frac{\Upsilon_x[\rho]}{\rho!}\,\gamma_{yx}(\mathcal{I}(\rho))\right)^{\beta(\rho)}\mathcal{Q}^\beta(\underline{1},\mathrm{d}w),\\
H^\tau_{(\tau_i)_{i=1}^n}(x,y) &:= (\Upsilon'[\tau\,\Xi])^{(n)}(u(x))\prod_{i=1}^n\left(\frac{\Upsilon_x[\tau_i]}{\tau_i!}\,\gamma_{yx}(\mathcal{I}(\tau_i))\right)\\
&\qquad\times\left\{u'(y)-u'(x)-\sum_{\substack{\rho\in\mathcal{T}\\ 1<|\mathcal{I}(\rho)|<\gamma-|\tau|+1-\sum_{i=1}^n|\mathcal{I}(\tau_i)|}}\frac{\Upsilon_x[\rho]}{\rho!}\,\gamma_{yx}(\mathcal{I}'(\rho))\right\},
\end{aligned}$$

*and all the derivatives of $\Upsilon'[\tau\,\Xi]$ appearing are bounded.*

**Proof.** We start by considering $\tau$ without polynomial decorations, and without loss of generality we assume that $\gamma_\tau := \gamma - |\tau| > 1$ since otherwise the argument is the same with the notational simplification that the terms with $u'$ do not appear. Given $x, y \in B_\lambda(z)$ we define the function $\boldsymbol{U}_{\gamma_\tau; x,y}: [0,1]^{\mathcal{T}_{\gamma_\tau - 2}} \to \mathbb{R}$ as

$$z \mapsto \boldsymbol{U}_{\gamma_\tau; x,y}(z) := \sum_{\rho \in \mathcal{T}_{\gamma_\tau - 2}} z_\rho \frac{\Upsilon_x[\rho]}{\rho!} \gamma_{yx}(\mathcal{I}(\rho)). \tag{2.22}$$

Since $U_{\gamma_\tau}$ is function-like one can show as in [EW24, Lemma 2.19] that

$$\begin{aligned}\langle \mathbf{1}, U_{\gamma_\tau}(y) - \Gamma_{yx} U_{\gamma_\tau}(x) \rangle &= u(y) - u(x) - u'(x)\,(y-x)_1 - \sum_{\rho \in \mathcal{T}_{\gamma_\tau - 2}} \frac{\Upsilon_x[\rho]}{\rho!} \gamma_{yx}(\mathcal{I}(\rho))\\ &= u(y) - u(x) - u'(x)\,(y-x)_1 - \boldsymbol{U}_{\gamma_\tau; x,y}(\underline{1}),\end{aligned}$$

i.e., $\boldsymbol{U}_{\gamma_\tau; x,\cdot}(\underline{1})$, for $\underline{1} \in \{1\}^{\mathcal{T}_{\gamma_\tau - 2}}$, contains the non-polynomial part of $\Gamma_{yx} U_{\gamma_\tau}(x)$. By assumption $\tau \in \mathcal{T}_\gamma$ has no polynomial decorations, which by Lemma 2.7 implies that $\Upsilon[\tau\,\Xi]$ is only a function of one variable and we can consider the function $F_\tau: [0,1]^{\mathcal{T}_{\gamma_\tau - 2}} \to \mathbb{R}$ defined as

$$F_\tau(z) := \Upsilon[\tau\,\Xi](u(x) + u'(x)\,(y-x)_1 + \boldsymbol{U}_{\gamma_\tau; x,y}(z)). \tag{2.23}$$

In particular we have

$$F_\tau(\underline{1}) = \Upsilon[\tau\,\Xi](u(y) - \langle \mathbf{1}, U_{\gamma_\tau}(y) - \Gamma_{yx} U_{\gamma_\tau}(x) \rangle). \tag{2.24}$$

For $\rho \in \mathcal{T}_{\gamma_\tau - 2}$ let $\partial_\rho = \partial_{z_\rho}$ be the partial derivative corresponding to the coordinate $z_\rho$ of $z \in [0,1]^{\mathcal{T}_{\gamma_\tau - 2}}$. Since $\boldsymbol{U}_{\gamma_\tau; x,y}(z)$ is a linear function of $z$, we have by chain rule that

$$\begin{aligned}\partial_\rho F_\tau(z) &= (\Upsilon[\tau\,\Xi])^{(1)}(u(x) + u'(x)\,(y-x)_1 + \boldsymbol{U}_{\gamma_\tau; x,y}(z))\, \partial_\rho \boldsymbol{U}_{\gamma_\tau; x,y}(z)\\ &= (\Upsilon[\tau\,\Xi])^{(1)}(u(x) + u'(x)\,(y-x)_1 + \boldsymbol{U}_{\gamma_\tau; x,y}(z)) \frac{\Upsilon_x[\rho]}{\rho!} \gamma_{yx}(\mathcal{I}(\rho)),\end{aligned}$$

i.e., the effect of the derivative is to differentiate $\Upsilon[\tau\,\Xi]$ and multiply by a $\rho$-dependent constant. Iterating this we obtain for any multi-index $\beta \in \mathbb{N}^{\mathcal{T}_{\gamma_\tau - 2}}$ and $\partial^\beta = \prod_{\rho \in \mathcal{T}} \partial_\rho^{\beta(\rho)}$ that

$$\partial^\beta F_\tau(z) = (\Upsilon[\tau\,\Xi])^{(\sum_\rho \beta(\rho))}(u(x) + u'(x)\,(y-x)_1 + \boldsymbol{U}_{\gamma_\tau; x,y}(z)) \prod_{\rho \in \mathcal{T}} \left( \frac{\Upsilon_x[\rho]}{\rho!} \gamma_{yx}(\mathcal{I}(\rho)) \right)^{\beta(\rho)}. \tag{2.25}$$

Since $\beta! = \prod_\rho \beta(\rho)!$, then re-indexing the sum via the identification $\beta \leftrightarrow \mathfrak{a}(\beta) = \prod_\rho \mathcal{I}(\rho)^{\beta(\rho)}$ we obtain for any set of multi-indexes $A \subset \mathbb{N}^{\mathcal{T}_{\gamma_\tau - 2}}$ the identity:

$$\begin{aligned}&\sum_{\beta \in A} \frac{\partial^\beta F_\tau(0)}{\beta!}\\ &= \sum_{\beta \in A} \frac{1}{\beta!} (\Upsilon[\tau\,\Xi])^{(\sum_\rho \beta(\rho))}(u(x) + u'(x)\,(y-x)_1) \prod_{\rho \in \mathcal{T}} \left( \frac{\Upsilon_x[\rho]}{\rho!} \gamma_{yx}(\mathcal{I}(\rho)) \right)^{\beta(\rho)}\\ &= \sum_{\beta \in A} \frac{1}{\prod_\rho \beta(\rho)!} (\Upsilon[\tau\,\Xi])^{(\sum_\rho \beta(\rho))}(u(x) + u'(x)\,(y-x)_1) \prod_{\rho \in \mathcal{T}} \left( \frac{\Upsilon_x[\rho]}{\rho!} \gamma_{yx}(\mathcal{I}(\rho)) \right)^{\beta(\rho)}\\ &= \sum_{n \in \mathbb{N}} \sum_{(\tau_i, \beta_i)_{i=1}^n \subset \mathcal{T} \times \mathbb{N}} \frac{1}{\prod_{i=1}^n \beta_i!} \prod_{i=1}^n \left( \frac{\Upsilon_x[\tau_i]}{\tau_i!} \gamma_{yx}(\mathcal{I}(\tau_i)) \right)^{\beta_i} (\Upsilon[\tau\,\Xi])^{(\sum_i \beta_i)}(u(x) + u'(x)\,(y-x)_1)\, \mathbb{1}_A\\ &= \sum_{n \in \mathbb{N}} \sum_{(\tau_i, \beta_i)_{i=1}^n \subset \mathcal{T} \times \mathbb{N}} \frac{1}{\prod_{i=1}^n \tau_i!^{\beta_i} \beta_i!} (\Upsilon[\tau\,\Xi])^{(\sum_i \beta_i)}(u(x) + u'(x)\,(y-x)_1) \prod_{i=1}^n (\Upsilon_x[\tau_i]\, \gamma_{yx}(\mathcal{I}(\tau_i)))^{\beta_i}\, \mathbb{1}_A\\ &= \sum_{n \in \mathbb{N}} \sum_{(\tau_i, \beta_i)_{i=1}^n \subset \mathcal{T} \times \mathbb{N}} \frac{1}{(\prod_{i=1}^n \mathcal{I}(\tau_i)^{\beta_i})!} (\Upsilon[\tau\Xi])^{(\sum_i \beta_i)}(u(x) + u'(x)(y-x)_1) \prod_{i=1}^n (\Upsilon_x[\tau_i]\, \gamma_{yx}(\mathcal{I}(\tau_i)))^{\beta_i}\, \mathbb{1}_A,\end{aligned}$$

where in the last identity we used the definition of the symmetry factor (2.3) and the sum is restricted by $\mathbb{1}_A$ to those trees for which the associated multi-index belongs to $A$ (see (2.11)). Re-indexing the sum via the identification $\prod_{i=1}^n \mathcal{I}(\tau_i)^{\beta_i} \leftrightarrow \prod_{i=1}^n \mathcal{I}(\tau_i)$ leads to a factor $\delta((\mathcal{I}(\tau_i))_{i=1}^n)$ as in (2.4) that comes from the non-uniqueness of the latter representation, with which we obtain:

$$
\begin{aligned}
&\sum_{\beta\in A} \frac{\partial^\beta F_\tau(0)}{\beta!}\\
&= \sum_{n\in\mathbb{N}} \sum_{\tau_1,\ldots,\tau_n\in\mathcal{T}} \frac{1}{\delta((\mathcal{I}(\tau_i))_{i=1}^n)(\prod_{i=1}^n \mathcal{I}(\tau_i))!} (\Upsilon[\tau\Xi])^{(n)}(u(x)+u'(x)(y-x)_1)\prod_{i=1}^n (\Upsilon_x[\tau_i]\gamma_{yx}(\mathcal{I}(\tau_i)))\mathbb{1}_A\\
&= \sum_{n\in\mathbb{N}} \sum_{\tau_1,\ldots,\tau_n\in\mathcal{T}} \frac{1}{n!\prod_{i=1}^n \tau_i!} (\Upsilon[\tau\,\Xi])^{(n)}(u(x)+u'(x)\,(y-x)_1) \prod_{i=1}^n (\Upsilon_x[\tau_i]\,\gamma_{yx}(\mathcal{I}(\tau_i)))\,\mathbb{1}_A\\
&= \sum_{n\in\mathbb{N}} \frac{1}{n!} \sum_{\tau_1,\ldots,\tau_n\in\mathcal{T}} \prod_{i=1}^n \left(\frac{\Upsilon_x[\tau_i]}{\tau_i!}\right)\gamma_{yx}\left(\prod_{i=1}^n \mathcal{I}(\tau_i)\right)(\Upsilon[\tau\,\Xi])^{(n)}(u(x)+u'(x)\,(y-x)_1)\,\mathbb{1}_A. \qquad (2.26)
\end{aligned}
$$

Using Lemma 2.6 (with $\Gamma = \mathrm{Id}$) and (2.24) we can write

$$
\begin{aligned}
H^\tau(x,y) &= \langle \tau, \hat{\sigma}(U_\gamma)(y) - \Gamma_{yx}\,\hat{\sigma}(U_\gamma)(x)\rangle\\
&= \Upsilon[\tau\,\Xi](u(y)) - \langle \tau, \Gamma_{yx}\,\hat{\sigma}(U_{\gamma_\tau})(x)\rangle\\
&= \Upsilon[\tau\,\Xi](u(y)) - \Upsilon[\tau\,\Xi](u(y) - \langle \mathbf{1}, U_{\gamma_\tau}(y) - \Gamma_{yx}\,U_{\gamma_\tau}(x)\rangle) + F_\tau(\underline{1}) - \langle \tau, \Gamma_{yx}\,\hat{\sigma}(U_{\gamma_\tau})(x)\rangle\\
&= H_1^\tau(x,y) + F_\tau(\underline{1}) - \sum_{\beta\in A_\tau} \frac{\partial^\beta F_\tau(0)}{\beta!} + \sum_{\beta\in A_\tau} \frac{\partial^\beta F_\tau(0)}{\beta!} - \langle \tau, \Gamma_{yx}\,\hat{\sigma}(U_{\gamma_\tau})(x)\rangle. \qquad (2.27)
\end{aligned}
$$

By Lemma 2.8 and (2.25) we have that

$$
F_\tau(\underline{1}) - \sum_{\beta\in A_\tau} \frac{\partial^\beta F_\tau(0)}{\beta!} = \sum_{\beta\in\partial A_\tau} \int_{\mathbb{R}^{\mathcal{T}_{\gamma_\tau-2}}} \partial^\beta F(z)\,\mathcal{Q}^\beta(\underline{1},\mathrm{d}z) = \sum_{\beta\in\partial A_\tau} H_\beta^\tau(x,y).
$$

For $\beta\in\partial A_\tau$ the term $\partial^\beta F$ contains the term $(\Upsilon[\tau\,\Xi])^{(\sum_\rho \beta(\rho))}$, to see that this derivative is well-defined recall that by definition of $\partial A_\tau$ there exist $\rho_0\in\mathcal{T}_{\gamma_\tau-2}$ such that $\beta - e_{\rho_0}\in A_\tau$. Since every factor in $\mathfrak{a}(\beta)$, see (2.11), has homogeneity at least $\alpha$, one has

$$
\alpha\left(\sum_\rho \beta(\rho) - 1\right) \leqslant |\mathfrak{a}(\beta - e_{\rho_0})| < \gamma_\tau = \gamma - |\tau|.
$$

Because $\tau$ has no polynomial decoration, then $\alpha\,\mathfrak{C}(\tau) = |\tau|$ by (2.9), and since $\alpha\in(0,1)$ we can conclude that $\mathbb{N}\ni\mathfrak{C}(\tau) + \sum_\rho \beta(\rho) - 1 < \gamma/_\alpha$ and $\mathfrak{C}(\tau) + \sum_\rho \beta(\rho) \leqslant \lceil \gamma/_\alpha \rceil$ which makes $(\Upsilon[\tau\,\Xi])^{(\sum_\rho \beta(\rho))}$ well-defined and bounded by Lemma 2.7. To conclude this case we identify the remaining terms in (2.27). In (2.26) $\tau_1,\ldots,\tau_n\in\mathcal{T}$ satisfy $\sum_{i=1}^n |\mathcal{I}(\tau_i)| < \gamma_\tau$ iff the multi-index associated to $\prod_{i=1}^n \mathcal{I}(\tau_i)$ (i.e., $\mathfrak{a}(\beta) = \prod_{i=1}^n \mathcal{I}(\tau_i)$) belongs to $A_\tau$, and therefore using Lemma 2.4 we can write

$$
\sum_{\beta\in A_\tau} \frac{\partial^\beta F_\tau(0)}{\beta!} - \langle \tau, \Gamma_{yx}\,\hat{\sigma}(U_{\gamma_\tau})(x)\rangle = \sum_{n\in\mathbb{N}} \frac{1}{n!} \sum_{\substack{\tau_1,\ldots,\tau_n\in\mathcal{T}\\ \sum_{i=1}^n |\mathcal{I}(\tau_i)| < \gamma_\tau}} H^\tau_{(\tau_i)_{i=1}^n}(x,y).
$$

Now consider $\tau$ which contains a polynomial decoration, we proceed similarly as before by considering the function $F_\tau\colon [0,1]^{\mathcal{T}_{\gamma_\tau-2}} \to \mathbb{R}$ defined as

$$
F_\tau(z) := \Upsilon'[\tau\,\Xi](u(x) + \boldsymbol{U}_{\gamma_\tau;x,y}(z)), \qquad (2.28)
$$

where $\boldsymbol{U}$ is defined as in (2.22) and from where we obtain that

$$
F_\tau(\underline{1}) = \Upsilon'[\tau\,\Xi](u(y) - \langle \mathbf{1}, U_{\gamma_\tau}(y) - \Gamma_{yx}\,U_{\gamma_\tau}(x)\rangle).
$$

Since $\tau$ has a polynomial decoration $|\tau| \geqslant |\boldsymbol{X}| = 1$ and $\gamma_\tau < 1$, so the linear term $u'(x)\,(y-x)_1$ in (2.23) is not needed. Also observe that in contrast to (2.23) which is defined in terms of $\Upsilon$, (2.28) is defined with $\Upsilon'$ as defined in Lemma 2.7. Analogously to before we obtain

$$\partial^\beta F_\tau(z) = (\Upsilon'[\tau\,\Xi])^{(\sum_\rho \beta(\rho))}\left(u(x) + \sum_{\rho \in \mathcal{T}_{\gamma_\tau - 2}} z_\rho \frac{\Upsilon_x[\rho]}{\rho!}\,\gamma_{yx}(\mathcal{I}(\rho))\right) \prod_{\rho \in \mathcal{T}} \left(\frac{\Upsilon_x[\rho]}{\rho!}\,\gamma_{yx}(\mathcal{I}(\rho))\right)^{\beta(\rho)},$$

and

$$\sum_{\beta \in A_\tau} \frac{\partial^\beta F_\tau(0)}{\beta!} = \sum_{n \in \mathbb{N}} \frac{1}{n!}\,(\Upsilon'[\tau\,\Xi])^{(n)}(u(x)) \sum_{\tau_1,\ldots,\tau_n \in \mathcal{T}} \prod_{i=1}^{n} \left(\frac{\Upsilon_x[\tau_i]}{\tau_i!}\,\gamma_{yx}(\mathcal{I}(\tau_i))\right) \mathbb{1}_A. \tag{2.29}$$

From the definition of the set $A_\tau$ and Lemma 2.4 we obtain

$$\begin{aligned}
&\langle \tau, \Gamma_{yx}\,\hat{\sigma}(U_\gamma)\rangle(u(x), u'(y)) - u'(y) \sum_{\beta \in A_\tau} \frac{\partial^\beta F_\tau(0)}{\beta!} \\
= &\sum_{n \in \mathbb{N}} \frac{1}{n!}\,(\Upsilon'[\tau\,\Xi])^{(n)}(u(x)) \sum_{\substack{\tau_1,\ldots,\tau_n \in \mathcal{T} \\ \sum_{i=1}^n |\mathcal{I}(\tau_i)| < \gamma_\tau}} \prod_{i=1}^{n} \left(\frac{\Upsilon_x[\tau_i]}{\tau_i!}\,\gamma_{yx}(\mathcal{I}(\tau_i))\right) \\
&\times \sum_{\substack{\tau_0 \in \mathcal{T} \\ 0 < |\mathcal{I}'(\tau_0)| < \gamma_\tau - \sum_{i=1}^n |\mathcal{I}(\tau_i)|}} \frac{\Upsilon_x[\tau_0]}{\tau_0!}\,\gamma_{yx}(\mathcal{I}'(\tau_0))
\end{aligned}$$

Considering $\langle \tau, \Gamma\,\hat{\sigma}(U)\rangle$ as a function of $(u, u')$ (see Lemma 2.6) and using (2.28) and the splitting $\Upsilon = u'\,\Upsilon'$ from Lemma 2.7 we can write

$$\begin{aligned}
&H^\tau(x,y) \\
= &\Upsilon[\tau\,\Xi](u(y), u'(y)) - \langle \tau, \Gamma_{yx}\,\hat{\sigma}(U)\rangle(u(x), u'(x)) \\
= &\Upsilon[\tau\,\Xi](u(y), u'(y)) - \Upsilon[\tau\,\Xi](u(y) - \langle \mathbf{1}, U_{\gamma_\tau+1}(y) - \Gamma_{yx}\,U_{\gamma_\tau+1}(x)\rangle, u'(y)) \\
&+ \Upsilon[\tau\Xi](u(y) - \langle \mathbf{1}, U_{\gamma_\tau+1}(y) - \Gamma_{yx} U_{\gamma_\tau+1}(x)\rangle, u'(y)) - \Upsilon[\tau\Xi](u(y) - \langle \mathbf{1}, U_{\gamma_\tau}(y) - \Gamma_{yx} U_{\gamma_\tau}(x)\rangle, u'(y)) \\
&+ \Upsilon[\tau\,\Xi](u(y) - \langle \mathbf{1}, U_{\gamma_\tau}(y) - \Gamma_{yx}\,U_{\gamma_\tau}(x)\rangle, u'(y)) - u'(y) \sum_{\beta \in A_\tau} \frac{\partial^\beta F_\tau(0)}{\beta!} \\
&+ u'(y) \sum_{\beta \in A_\tau} \frac{\partial^\beta F_\tau(0)}{\beta!} - \langle \tau, \Gamma_{yx}\,\hat{\sigma}(U)\rangle(u(x), u'(x)) \\
= &u'(y) \left\{ H_1^\tau(x,y) + H_2^\tau(x,y) + F_\tau(\underline{1}) - \sum_{\beta \in A_\tau} \frac{\partial^\beta F_\tau(0)}{\beta!} \right\} \\
&+ u'(y) \sum_{\beta \in A_\tau} \frac{\partial^\beta F_\tau(0)}{\beta!} - \langle \tau, \Gamma_{yx}\,\hat{\sigma}(U)\rangle(u(x), u'(x)).
\end{aligned}$$

By Lemma 2.8 and (2.28) one has

$$F_\tau(\underline{1}) - \sum_{\beta \in A} \frac{\partial^\beta F_\tau(0)}{\beta!} = \sum_{\beta \in \partial A_\tau} \int_{\mathbb{R}^{\mathcal{T}_{\gamma_\tau - 2}}} \partial^\beta F(z)\,\mathcal{Q}^\beta(\underline{1}, \mathrm{d}z) = \sum_{\beta \in \partial A_\tau} H_\beta^\tau(x,y).$$

Analogously to the non-decoration case, for $\beta \in \partial A_\tau$ the term $(\Upsilon'[\tau\,\Xi])^{(\sum_\rho \beta(\rho))}$

The term $(\Upsilon'[\tau\,\Xi])^{(\sum_\rho \beta(\rho))}$ appearing in $\partial^\beta F$ for $\beta \in \partial A_\tau$ is well-defined, and bounded, because there exist $\rho_0 \in \mathcal{T}_{\gamma_\tau - 2}$ such that

$$\alpha\left(\sum_\rho \beta(\rho) - 1\right) \leqslant |\mathfrak{a}(\beta - e_{\rho_0})| < \gamma_\tau = \gamma - |\tau|,$$

and since in this case $\tau$ has a polynomial decoration, then $\alpha\,\mathfrak{C}(\tau)=|\tau|-1$, see (2.9), which allows us to conclude that $\mathbb{N}\ni\mathfrak{C}(\tau)+\sum_\rho\beta(\rho)<\gamma/_\alpha$ and $1+\mathfrak{C}(\tau)+\sum_\rho\beta(\rho)\leqslant\lceil\gamma/_\alpha\rceil$

At last, the highest order derivatives of $(\Upsilon'[\tau\,\Xi])^{(n)}$ appear when $n=\sum_\rho\beta(\rho)$ for some $\beta\in\partial A_\tau$. Analogously to the non-decoration case, there exists $\rho_0\in\mathcal{T}_{\gamma_\tau-2}$ such that $\beta-e_{\rho_0}\in A_\tau,$ and therefore

$$\alpha\left(\sum_\rho\beta(\rho)-1\right)\leqslant|\mathfrak{a}(\beta-e_{\rho_0})|<\gamma_\tau=\gamma-|\tau|.$$

Using that $\alpha\,\mathfrak{C}(\tau)=|\tau|-1$, see (2.9), and $\alpha\in(0,1)$ we conclude that $\mathbb{N}\ni\mathfrak{C}(\tau)+\sum_\rho\beta(\rho)<\gamma/_\alpha$ and $1+\mathfrak{C}(\tau)+\sum_\rho\beta(\rho)\leqslant\lceil\gamma/_\alpha\rceil$ which makes $(\Upsilon'[\tau\,\Xi])^{(\sum_\rho\beta(\rho))}$ well-defined by Lemma 2.7.

At last, by Lemma 2.4 and (2.29)

$$\begin{aligned}
& u'(y)\sum_{\beta\in A_\tau}\frac{\partial^\beta F_\tau(0)}{\beta!}-\langle\tau,\Gamma_{yx}\,\hat{\sigma}(U)\rangle(u(x),u'(x))\\
=\;& \langle\tau,\Gamma_{yx}\,\hat{\sigma}(U)\rangle(u(x),u'(y))-\langle\tau,\Gamma_{yx}\,\hat{\sigma}(U)\rangle(u(x),u'(x))\\
& -\left(\langle\tau,\Gamma_{yx}\,\hat{\sigma}(U)\rangle(u(x),u'(y))-u'(y)\sum_{\beta\in A_\tau}\frac{\partial^\beta F_\tau(0)}{\beta!}\right)\\
=\;& \{u'(y)-u'(x)\}\sum_{n\in\mathbb{N}}\frac{1}{n!}(\Upsilon'[\tau\,\Xi])^{(n)}(u(x))\sum_{\substack{\tau_1,\ldots,\tau_n\in\mathcal{T}\\ \sum_{i=1}^n|\mathcal{I}(\tau_i)|<\gamma_\tau}}\prod_{i=1}^n\left(\frac{\Upsilon_x[\tau_i]}{\tau_i!}\gamma_{yx}(\mathcal{I}(\tau_i))\right)\\
& -\sum_{n\in\mathbb{N}}\frac{1}{n!}(\Upsilon'[\tau\,\Xi])^{(n)}(u(x))\sum_{\substack{\tau_1,\ldots,\tau_n\in\mathcal{T}\\ \sum_{i=1}^n|\mathcal{I}(\tau_i)|<\gamma_\tau}}\prod_{i=1}^n\left(\frac{\Upsilon_x[\tau_i]}{\tau_i!}\gamma_{yx}(\mathcal{I}(\tau_i))\right)\\
& \quad\times\sum_{\substack{\tau_0\in\mathcal{T}\\ 0<|\mathcal{I}'(\tau_0)|<\gamma_\tau-\sum_{i=1}^n|\mathcal{I}(\tau_i)|}}\frac{\Upsilon_x[\tau_0]}{\tau_0!}\gamma_{yx}(\mathcal{I}'(\tau_0))\\
=\;& \sum_{n\in\mathbb{N}}\frac{1}{n!}\sum_{\substack{\tau_1,\ldots,\tau_n\in\mathcal{T}\\ \sum_{i=1}^n|\mathcal{I}(\tau_i)|<\gamma_\tau}}H^\tau_{(\tau_i)_{i=1}^n}(x,y).
\end{aligned}$$

□

We are now in place to prove Lemma 2.5.

**Proof of Lemma 2.5.** We start by considering $\tau$ without polynomial decorations, and without loss of generality we assume that $\gamma_\tau:=\gamma-|\tau|>1$ since in the other case the proof reduces to [BCMW22, Corollary 3.12]. Our starting point is the decomposition (2.19) from Lemma 2.9. By Lemma 2.7, $\Upsilon[\tau\,\Xi]$ is function of $\sigma$ and its derivatives up to degree $\mathfrak{C}(\tau)$, then $(\Upsilon[\tau\,\Xi])^{(1)}$ is bounded by our assumptions, and therefore for all $x,y\in B_\lambda(z)$ the germ $H_1^\tau$ can be bounded as

$$\begin{aligned}
|H_1^\tau(x,y)| &\lesssim \sup_{\lambda\in[0,1]}|(\Upsilon[\tau\,\Xi])^{(1)}(u(y)-\lambda\,\langle\mathbf{1},U_{\gamma_\tau}(y)-\Gamma_{yx}\,U_{\gamma_\tau}(x)\rangle)|\,|\langle\mathbf{1},U_{\gamma_\tau}(y)-\Gamma_{yx}\,U_{\gamma_\tau}(x)\rangle|\\
&\lesssim [U_{\gamma_\tau};\mathbf{1}]_{B_\lambda(z)}\,d(x,y)^{\gamma_\tau}.
\end{aligned}\tag{2.30}$$

For any $\beta\in\partial A_\tau$ we have by Lemma 2.7, definition (2.11) of $\mathfrak{a}(\beta)$ and the multiplicativity of $\gamma$:

$$\begin{aligned}
\left|\prod_{\rho\in\mathcal{T}_{\gamma_\tau-2}}\left(\frac{\Upsilon_x[\rho]}{\rho!}\gamma_{yx}(\mathcal{I}(\rho))\right)^{\beta(\rho)}\right| &\lesssim |\gamma_{yx}(\mathfrak{a}(\beta))|\prod_{\rho\in\mathcal{T}_{\gamma_\tau-2}}|\Upsilon_x[\rho]|^{\beta(\rho)}\\
&\lesssim d(x,y)^{|\mathfrak{a}(\beta)|}\,[\Gamma\,\mathfrak{a}(\beta);\mathbf{1}]_{B_\lambda(z)}\prod_{\substack{\rho\in\mathcal{T}_{\gamma_\tau-2}\\ \mathfrak{n}(\rho)\neq0}}|u'(x)|^{\beta(\rho)}\\
&\lesssim \|u'\|_{B_\lambda(z)}^{|\mathfrak{n}(\mathfrak{a}(\beta))|}\,[\Gamma\,\mathfrak{a}(\beta);\mathbf{1}]_{B_\lambda(z)}\,d(x,y)^{|\mathfrak{a}(\beta)|},
\end{aligned}\tag{2.31}$$

since $\mathfrak{n}(\mathfrak{a}(\beta)) = \sum_\rho \beta(\rho)\,\mathfrak{n}(\rho) \in \{0, e_{(0,1)}\}$. By definition of $\partial A_\tau$ we have that $|\mathfrak{a}(\beta)| \geqslant \gamma_\tau$ and by (2.23) we conclude that

$$\begin{aligned}|H_\beta(x,y)| &\lesssim \int_{\mathbb{R}^{\mathcal{T}_{\gamma_\tau-2}}} |\partial^\beta F(z)|\,\mathcal{Q}^\beta(\underline{1},\mathrm{d}z)\\ &\leqslant \left\| (\Upsilon[\tau\,\Xi])^{(\sum_\rho \beta(\rho))}(u(x)+u'(x)\,(y-x)_1+\boldsymbol{U}_{\gamma;x,y}(\cdot)) \prod_{\rho\in\mathcal{T}_{\gamma_\tau-2}} \left(\frac{\Upsilon_x[\rho]}{\rho!}\,\gamma_{yx}(\mathcal{I}(\rho))\right)^{\beta(\rho)} \right\|\\ &\lesssim \lambda^{|\mathfrak{a}(\beta)|-\gamma_\tau}\,\|u'\|_{B_\lambda(z)}^{|\mathfrak{n}(\mathfrak{a}(\beta))|}\,[\Gamma\,\mathfrak{a}(\beta);\mathbf{1}]_{B_\lambda(z)}\,d(x,y)^{\gamma_\tau} \end{aligned} \tag{2.32}$$

since $\Upsilon[\tau\,\Xi]$ has bounded derivative. For the last terms in (2.19) observe that $\tau_1,\ldots,\tau_n \in \mathcal{T}$ satisfy $\sum_{i=1}^n |\mathcal{I}(\tau_i)| < \gamma_\tau$ iff the multi-index $\beta$ such that $\mathfrak{a}(\beta) = \prod_{i=1}^n \mathcal{I}(\tau_i)$ belongs to $A_\tau$, and therefore in this case

$$\begin{aligned}\left|\prod_{i=1}^n \frac{\Upsilon_x[\tau_i]}{\tau_i!}\,\gamma_{yx}(\mathcal{I}(\tau_i))\right| &\lesssim |\gamma_{yx}(\mathfrak{a}(\beta))| \prod_{\rho\in\mathcal{T}_{\gamma_\tau-2}} |\Upsilon_x[\rho]|^{\beta(\rho)}\\ &\lesssim d(x,y)^{|\mathfrak{a}(\beta)|}\,[\Gamma\,\mathfrak{a}(\beta);\mathbf{1}]_{B_\lambda(z)} \prod_{\substack{\rho\in\mathcal{T}_{\gamma_\tau-2}\\ \mathfrak{n}(\rho)\neq 0}} |u'(x)|^{\beta(\rho)}\\ &\lesssim \|u'\|_{B_\lambda(z)}^{|\mathfrak{n}(\mathfrak{a}(\beta))|}\,[\Gamma\,\mathfrak{a}(\beta);\mathbf{1}]_{B_\lambda(z)}\,d(x,y)^{|\mathfrak{a}(\beta)|}.\end{aligned}$$

In contrast to the previous case now we have $\gamma_\tau > |\mathfrak{a}(\beta)|$ and therefore this term is not enough to obtain the required $d(x,y)^{\gamma_\tau}$ control on $H^\tau_{(\tau_i)_{i=1}^n}(x,y)$. The missing factor can be obtained from

$$\begin{aligned}&|(\Upsilon[\tau\Xi])^{(n)}(u(x)+u'(x)(y-x)_1) - (\Upsilon[\tau\Xi])^{(n)}(u(x)) - (\Upsilon[\tau\Xi])^{(n+1)}(u(x))\,u'(x)((y-x)_1)\,\mathbb{1}_{|\mathfrak{a}(\beta)|<\gamma_\tau-1}|\\ &\lesssim \begin{cases} \|u'\|_{B_\lambda(z)}\,d(x,y) & \text{if } |\mathfrak{a}(\beta)| \geqslant \gamma_\tau - 1\\ \|u'\|_{B_\lambda(z)}^2\,d(x,y)^2 & \text{if } |\mathfrak{a}(\beta)| < \gamma_\tau - 1 \end{cases}.\end{aligned}$$

For the analysis of Section 3.2 the quadratic term $\|u'\|$ is too bad, and since a factor $d(x,y)^{\gamma_\tau-|\mathfrak{a}(\beta)|}$ is needed we can interpolate ($\theta = \gamma_\tau - |\mathfrak{a}(\beta)| \in (0,1)$ in the first case and $\theta = \gamma_\tau - |\mathfrak{a}(\beta)| - 1 \in (0,1)$ for the second one) with the following coarser bound

$$\begin{aligned}&|(\Upsilon[\tau\Xi])^{(n)}(u(x)+u'(x)(y-x)_1) - (\Upsilon[\tau\Xi])^{(n)}(u(x)) - (\Upsilon[\tau\Xi])^{(n+1)}(u(x))\,u'(x)((y-x)_1)\,\mathbb{1}_{|\mathfrak{a}(\beta)|<\gamma_\tau-1}|\\ &\lesssim \begin{cases} 1 & \text{if } |\mathfrak{a}(\beta)| \geqslant \gamma_\tau - 1\\ \|u'\|_{B_\lambda(z)}\,d(x,y) & \text{if } |\mathfrak{a}(\beta)| < \gamma_\tau - 1 \end{cases},\end{aligned}$$

and conclude the bound

$$\begin{aligned}&|(\Upsilon[\tau\Xi])^{(n)}(u(x)+u'(x)(y-x)_1) - (\Upsilon[\tau\Xi])^{(n)}(u(x)) - (\Upsilon[\tau\Xi])^{(n+1)}(u(x))\,u'(x)((y-x)_1)\,\mathbb{1}_{|\mathfrak{a}(\beta)|<\gamma_\tau-1}|\\ &\lesssim \|u'\|_{B_\lambda(z)}^{\gamma_\tau-|\mathfrak{a}(\beta)|}\,d(x,y)^{\gamma_\tau-|\mathfrak{a}(\beta)|},\end{aligned}$$

which gives us the right power of $d(x,y)$ and a lower power in the $u'$ term. We conclude

$$|H^\tau_{(\tau_i)_{i=1}^n}(x,y)| \lesssim \|u'\|_{B_\lambda(z)}^{|\mathfrak{n}(\mathfrak{a}(\beta))|+\gamma_\tau-|\mathfrak{a}(\beta)|}\,[\Gamma\,\mathfrak{a}(\beta);\mathbf{1}]_{B_\lambda(z)}\,d(x,y)^{\gamma_\tau}. \tag{2.33}$$

Combining the bounds (2.30), (2.32), and (2.33) we obtain a bound for $[\hat{\sigma}(U_\gamma);\tau]_{B_\lambda(z)}$ and multiplying this bound by $\lambda^{\gamma_\tau}$ we obtain the stated result in the case where $\tau$ has no polynomial decorations.

Now we consider $\tau$ with polynomial decorations, for which we have $\gamma_\tau := \gamma - |\tau| \in (0,1)$. We bound each term from the decomposition of Lemma 2.9. For any $x,y \in B_\lambda(z)$ we have by Taylor's remainder

$$\begin{aligned}|u'(y)\,H_1^\tau(x,y)| &\lesssim \|u'\|_{B_\lambda(z)}\,\|(\Upsilon'[\tau\,\Xi])^{(1)}\|\,[U_{\gamma_\tau+1};\mathbf{1}]_{B_\lambda(z)}\,d(x,y)^{\gamma_\tau+1}\\ &\lesssim d(x,y)^{\gamma_\tau+1}\,\|u'\|_{B_\lambda(z)}\,[U_{\gamma_\tau+1};\mathbf{1}]_{B_\lambda(z)},\end{aligned} \tag{2.34}$$

which we can interpolate with the coarser bound

$$|u'(y)\, H_1^{\tau}(x,y)| \lesssim \|u'\|_{B_\lambda(z)}\, \|\Upsilon'[\tau\,\Xi]\| \lesssim \|u'\|_{B_\lambda(z)} \tag{2.35}$$

to obtain

$$|u'(y)\, H_1^{\tau}(x,y)| \lesssim \|u'\|_{B_\lambda(z)}\, [U_{\gamma_\tau+1};\mathbf{1}]_{B_\lambda(z)}^{\frac{\gamma_\tau}{\gamma_\tau+1}}\, d(x,y)^{\gamma_\tau}. \tag{2.36}$$

For the next term, we use the identity from [EW24, Lemma 2.19], and the fact that $\gamma_\tau \in (0,1)$ to write

$$\langle \mathbf{1}, \Gamma_{yx}\, U_{\gamma_\tau+1}(x)\rangle - \langle \mathbf{1}, \Gamma_{yx}\, U_{\gamma_\tau}(x)\rangle \;=\; u'(x)\,(y-x)_1 + \sum_{\substack{\rho\in\mathcal{T}\\ |\mathcal{I}(\rho)|\in[\gamma_\tau,\gamma_\tau+1)}} \frac{\Upsilon_x[\rho]}{\rho!}\, \gamma_{yx}(\mathcal{I}(\rho)), \tag{2.37}$$

which, by a Taylor remainder argument, implies the bound

$$\begin{aligned} &|u'(y)\, H_2^{\tau}(x,y)| \\ \leqslant\;& \|u'\|_{B_\lambda(z)}\, \|(\Upsilon'[\tau\,\Xi])^{(1)}\|\, |\langle \mathbf{1}, \Gamma_{yx}\, U_{\gamma_\tau+1}(x)\rangle - \langle \mathbf{1}, \Gamma_{yx}\, U_{\gamma_\tau}(x)\rangle| \\ \lesssim\;& \|u'\|_{B_\lambda(z)} \left( d(x,y)\, \|u'\|_{B_\lambda(z)} + \sum_{\substack{\rho\in\mathcal{T}\\ |\mathcal{I}(\rho)|\in[\gamma_\tau,\gamma_\tau+1)}} \|u'\|_{B_\lambda(z)}^{|\mathfrak{n}(\rho)|}\, [\Gamma\,\mathcal{I}(\rho);\mathbf{1}]\, d(x,y)^{|\mathcal{I}(\rho)|} \right). \end{aligned} \tag{2.38}$$

Analogously to (2.35), we also have the coarser bound $|u'(y)\, H_2^{\tau}(x,y)| \lesssim \|u'\|_{B_\lambda(z)}$. We can interpolate between this two bounds using the inequality:

$$x \leqslant a \wedge \sum_{i\in I} b_i \Longrightarrow x \lesssim C(\theta) \sum_{i\in I} a^{1-\theta_i}\, b_i^{\theta_i}, \qquad \theta\in[0,1]^I, \tag{2.39}$$

by setting $\theta_0 = \gamma_\tau$ for the term $d(x,y)\, \|u'\|_{B_\lambda(z)}$, $\theta_\rho = \gamma_\tau$ for trees with polynomial decoration and $\theta_\rho = 1$ for the ones without polynomial decoration. Applying these interpolation yields:

$$\begin{aligned} &|u'(y)\, H_2^{\tau}(x,y)| \\ \lesssim\;& \left( d(x,y)^{\gamma_\tau} + \sum_{\substack{\rho\in\mathcal{T},\mathfrak{n}(\rho)\neq 0\\ |\mathcal{I}(\rho)|\in[\gamma_\tau,\gamma_\tau+1)}} ([\Gamma\,\mathcal{I}(\rho);\mathbf{1}]\, d(x,y)^{|\mathcal{I}(\rho)|})^{\gamma_\tau} \right) \|u'\|_{B_\lambda(z)}^{\gamma_\tau+1} \\ &+ \|u'\|_{B_\lambda(z)} \sum_{\substack{\rho\in\mathcal{T},\mathfrak{n}(\rho)=0\\ |\mathcal{I}(\rho)|\in[\gamma_\tau,\gamma_\tau+1)}} [\Gamma\,\mathcal{I}(\rho);\mathbf{1}]\, d(x,y)^{|\mathcal{I}(\rho)|} \\ \lesssim\;& \|u'\|_{B_\lambda(z)}^{\gamma_\tau+1} \left( 1 + \sum_{\substack{\rho\in\mathcal{T},\mathfrak{n}(\rho)\neq 0\\ |\mathcal{I}(\rho)|\in[\gamma_\tau,\gamma_\tau+1)}} (\lambda^{(|\mathcal{I}(\rho)|-1)}\, [\Gamma\,\mathcal{I}(\rho);\mathbf{1}])^{\gamma_\tau} \right) d(x,y)^{\gamma_\tau} \\ &+ \|u'\|_{B_\lambda(z)} \left( \sum_{\substack{\rho\in\mathcal{T},\mathfrak{n}(\rho)=0\\ |\mathcal{I}(\rho)|\in[\gamma_\tau,\gamma_\tau+1)}} [\Gamma\,\mathcal{I}(\rho);\mathbf{1}]\, \lambda^{|\mathcal{I}(\rho)|-\gamma_\tau} \right) d(x,y)^{\gamma_\tau}. \end{aligned} \tag{2.40}$$

For each $\beta \in \partial A_\tau$ we have $|\mathfrak{a}(\beta)| - \gamma_\tau > 0$, and therefore by (2.28) and (2.31)

$$\begin{aligned} |u'(y)\, H_\beta^{\tau}(x,y)| \;\lesssim\;& \|u'\|_{B_\lambda(z)} \int_{\mathbb{R}^{\mathcal{T}_{\gamma_\tau-2}}} |\partial^\beta F(z)|\, \mathcal{Q}^\beta(\underline{1}, \mathrm{d}z) \\ \lesssim\;& \|u'\|_{B_\lambda(z)} \left\| (\Upsilon[\tau\,\Xi])^{(\sum_\rho \beta(\rho))} (u(x) + \boldsymbol{U}_{\gamma;x,y}(\cdot)) \prod_{\rho\in\mathcal{T}_{\gamma_\tau-2}} \left( \frac{\Upsilon_x[\rho]}{\rho!}\, \gamma_{yx}(\mathcal{I}(\rho)) \right)^{\beta(\rho)} \right\| \\ \lesssim\;& \lambda^{|\mathfrak{a}(\beta)|-\gamma_\tau}\, \|u'\|_{B_\lambda(z)}\, [\Gamma\,\mathfrak{a}(\beta);\mathbf{1}]_{B_\lambda(z)}\, d(x,y)^{\gamma_\tau}, \end{aligned} \tag{2.41}$$

where we used that $\partial A_\tau \subset \mathbb{N}^{\mathcal{T}_{\gamma_\tau - 2}}$ which in particular implies that $\mathfrak{n}(\mathfrak{a}(\beta)) = 0$, since for any $\rho \in \mathcal{T}_{\gamma_\tau - 2}$ we have $|\mathcal{I}(\rho)| < \gamma_\tau < \gamma - 1 < 1$, and which implies $\mathfrak{n}(\rho) = 0$. For the last terms in (2.21), let $\beta \in A_\tau$ be such that $\mathfrak{a}(\beta) = \prod_{i=1}^n \mathcal{I}(\tau_i)$. By [EW24, eq. (2.37) and Lemma 2.23] we have that

$$\begin{aligned}
&\left| u'(y) - u'(x) - \sum_{\substack{\rho \in \mathcal{T} \\ 1 < |\mathcal{I}(\rho)| < \gamma_\tau + 1 - \sum_{i=1}^n |\mathcal{I}(\tau_i)|}} \frac{\Upsilon[\rho]}{\rho!} \gamma_{yx}(\mathcal{I}'(\rho)) \right| \\
= \;& |\langle \boldsymbol{X}, U_{\gamma_\tau - |\mathfrak{a}(\beta)| + 1}(y) - \Gamma_{yx} U_{\gamma_\tau - |\mathfrak{a}(\beta)| + 1}(x) \rangle| \\
\lesssim \;& [U_{\gamma_\tau - |\mathfrak{a}(\beta)| + 1}; \boldsymbol{X}]_{B_\lambda(z)} \, d(x, y)^{\gamma_\tau - |\mathfrak{a}(\beta)| + 1 - |\boldsymbol{X}|} \\
= \;& [U_{\gamma_\tau - |\mathfrak{a}(\beta)| + 1}; \boldsymbol{X}]_{B_\lambda(z)} \, d(x, y)^{\gamma_\tau - |\mathfrak{a}(\beta)|}
\end{aligned}$$

i.e., this term models $u'$ to order $\gamma_\tau - |\mathfrak{a}(\beta)| > 0$ and contributes precisely the missing power of $d(x, y)$ in (2.31) (as before $\mathfrak{n}(\mathfrak{a}(\beta)) = 0$) which allows us to conclude

$$|H^\tau_{(\tau_i)_{i=1}^n}(x, y)| \lesssim [U_{\gamma_\tau - |\mathfrak{a}(\beta)| + 1}; \boldsymbol{X}]_{B_\lambda(z)} \, [\Gamma\, \mathfrak{a}(\beta); \mathbf{1}]_{B_\lambda(z)} \, d(x, y)^{\gamma_\tau} \tag{2.42}$$

Combining (2.36), (2.40), (2.41), and (2.42) and multiplying with $\lambda^{\gamma_\tau + 1}$ we conclude the result. □

**Remark 2.10.** In the previous proof, one could instead write

$$u'(y) \{H_1^\tau(x, y) + H_2^\tau(x, y)\} = u'(y) \{\Upsilon'[\tau\, \Xi](u(y)) - \Upsilon'[\tau\, \Xi](u(y) - \langle \mathbf{1}, U_{\gamma_\tau}(y) - \Gamma_{yx} U_{\gamma_\tau}(x) \rangle)\}$$

and bound this term directly as

$$|u'(y) \{H_1^\tau(x, y) + H_2^\tau(x, y)\}| \lesssim \|u'\|_{B_\lambda(z)} \, [U_{\gamma_\tau}; \mathbf{1}]_{B_\lambda(z)} \, d(x, y)^{\gamma_\tau},$$

which leads to a much simpler bound that does not require interpolations. However, the term

$$\|u'\| \, \lambda^{\gamma_\tau} \, [U_{\gamma_\tau}; \mathbf{1}] \lesssim \|u'\| \, \lambda^\gamma \, [U_\gamma; \mathbf{1}] + \cdots$$

would eventually become a problem in Section 3.1. More precisely, using (3.29) to bound $\|u'\|$ produces the term $\|u\|_{B_\lambda(z)}^{2 - \gamma^{-1}} (\lambda^\gamma \, [U_\gamma; \mathbf{1}]_{B_\lambda(z)})^{1 + \gamma^{-1}}$, which is superlinear in $\lambda^\gamma \, [U_\gamma; \mathbf{1}]$ and does not allow us to absorb it to the right-hand side (the sublinear power $\|u\|^{2 - \gamma^{-1}}$ gets absorbed using the smallness assumption on the model terms from Definition 3.1). See the Proof of Lemma 3.9 for details.

# 3. A Priori Estimates

## 3.1. Small scales estimates

In this section, we assume that we are given a coherent modelled distribution $U \in \mathcal{D}^\gamma$, for some $\gamma \in (2 - \alpha, 2)$, whose reconstruction $u := \mathcal{R}\, U$ (also its coefficient at $\mathbf{1}$ since $U$ is function-like) additionally solves the equation

$$(\partial_t - \partial_x^2)\, u = \beta\, u - u\, |u|^{m - 1} + \mathcal{R}(\hat{\sigma}(U_\gamma)\, \Xi), \tag{3.1}$$

which is the renormalised version of (1.3). The main result of this section, Theorem 3.2, is a control on high regularity modelled seminorms of the solution when restricted to small scales, and where this scale depends on the $L^\infty$-norm of the solution itself.

The precise smallness scale is given in the next definition:

**Definition 3.1**. *For $c \in (0, 1)$ and $t \in [0, 1)$ we define*

$$\lambda_t := (1 \vee \|u\|_{P_t})^{-\frac{(m - 1)}{2}} \wedge (c\, (1 + |||(\Pi, \Gamma)|||)^{-1} \, (1 \vee \|u\|_{P_t})^{-(\gamma - 1)})^{\frac{1}{\alpha}},$$

*where we recall that* $|||(\Pi,\Gamma)|||$ *is the homogeneous model norm as defined in* (2.14)*.*

**THEOREM 3.2**. *Let* $(\Pi;\Gamma)$ *be a weakly admissible model (Definition 2.1) and* $U_\gamma \in \mathcal{D}^\gamma$ *a coherent modelled distribution of the form* (2.16) *for some fixed* $\gamma \in (2-\alpha, 2)$ *such that* $u = \langle \mathbf{1}, U\rangle = \mathcal{R} U$ *solves* (3.1) *weakly in* $P := (0,1) \times (-1,1)$*. Fix* $t \in (0,1)$ *such that* $\|u\|_{P_t} \geqslant 1$ *and set* $\lambda_t > 0$ *as in Definition 3.1. Then, if* $t + 2\lambda_t < 1$ *and* $c = c(\alpha,\gamma) \in (0,1)$ *is small enough, we can conclude*

$$\sup_{\lambda \in (0,\lambda_t]} \lambda^\eta \sup_{z \in P_{t+2\lambda}} [U_\eta; \mathbf{1}]_{B_{\lambda/2}(z)} \lesssim_{\alpha,\gamma,m,|\beta|,\|\sigma\|_{C_b^k}} \|u\|_{P_t} \qquad \forall \eta \in (0,\gamma],$$

*where* $U_\eta = \mathcal{Q}_{<\eta} U_\gamma$ *is the projection of* $U_\gamma$ *to trees of homogeneity* $<\eta$ *and* $P_t = (t^2,1) \times (-1+t, 1-t)$ *is the set of points with parabolic distance* $t$ *away from the parabolic boundary of* $P$*.*

**COROLLARY 3.3**. *Under the same assumptions as Theorem 3.2 we have*

$$\begin{aligned}
\sup_{\lambda \in (0,\lambda_t]} \lambda \sup_{z \in P_{t+2\lambda}} \|u'\|_{B_{\lambda/2}(z)} &\lesssim \|u\|_{P_t} \\
\sup_{\lambda \in (0,\lambda_t]} \lambda^\eta \sup_{z \in P_{t+2\lambda}} [U_\eta; \boldsymbol{X}]_{B_{\lambda/2}(z)} &\lesssim \|u\|_{P_t} \qquad \forall \eta \in (1,\gamma] \\
\sup_{\lambda \in (0,\lambda_t]} \lambda^{\gamma - (|\tau| - |\mathfrak{n}(\tau)|)} \sup_{z \in P_{t+2\lambda}} [\hat{\sigma}(U_\gamma); \tau]_{B_{\lambda/2}(z)} &\lesssim \|u\|_{P_t}^{1 + \mathbb{1}_{\gamma - (|\tau| - |\mathfrak{n}(\tau)|) > 1} (\gamma - (|\tau| - |\mathfrak{n}(\tau)|))(1 - \gamma^{-1})},
\end{aligned}$$

*for all* $\tau \in \mathcal{T}_\gamma$ *with no noise edge at the root, and with the implicit proportionality constant depending on* $\alpha, \gamma, m, |\beta|$ *and* $\|\sigma\|_{C_b^k}$*.*

The main tool to prove this result is the following form of the Schauder estimate for germs:

**THEOREM 3.4**. *Fix* $\gamma \in (1,2)$ *and a finite set* $A \subset (0,\gamma)$*. Let* $U: \mathbb{R}^2 \times \mathbb{R}^2 \to \mathbb{R}$ *be a germ such that* $U(x,x) = 0$ *for all* $x \in \mathbb{R}^2$ *and*

$$[U]_{\gamma; B_1} := \sup_{x,y \in B_1} \frac{|U(x,y) - \nu(x)\,(y-x)_1|}{d(x,y)^\gamma} < +\infty,$$

*where* $\nu(x) := \partial U(x,\cdot)|_x \in \mathbb{R}$ *is the generalised spatial derivative of* $U$*. Assume that there exists a germ* $\Lambda: \mathbb{R} \times \mathbb{R} \to \mathbb{R}^2$ *such that* $\Lambda(x,x) = 0$ *for all* $x \in \mathbb{R}$ *and*

$$\begin{aligned}
[U]_{\gamma\text{-3pt}; B_1} &:= \sup_{x,y,z \in B_1} \frac{|U(x,z) - U(x,y) - U(y,z) + \Lambda(x,y)\,(z-y)_1|}{\sum_{\beta \in A} d(x,y)^\beta\, d(y,z)^{\gamma-\beta}} < +\infty. \\
[\Lambda]_{(\gamma-1)\text{-3pt}; B_1} &:= \sup_{x,y,z \in B_1} \frac{|\Lambda(x,z) - \Lambda(x,y) - \Lambda(y,z)|}{\sum_{\beta \in A \cap (1,\gamma)} d(x,y)^{\beta-1}\, d(y,z)^{\gamma-\beta}} < +\infty.
\end{aligned}$$

*Then, we have the estimate*

$$\lambda^\gamma [U]_{\gamma; B_{\lambda/2}(z)} \lesssim \lambda^\gamma \|(\partial_t - \partial_x^2) U\|_{\gamma-2; B_\lambda(z)} + \lambda^\gamma [U]_{\gamma\text{-3pt}; B_\lambda(z)} + \lambda^\gamma [\Lambda]_{(\gamma-1)\text{-3pt}; B_\lambda(z)} + \|U\|_{B_\lambda(z)}, \qquad (3.2)$$

*where* $\|U\|_B := \sup_{x,y \in B} |U(x,y)|$*,*

$$\|(\partial_t - \partial_x^2)\,U\|_{\gamma-2; B_\lambda(z)} := \sup_{x \in B_1} \sup_{\substack{r \in (0,1) \\ B_r(x) \subset B_\lambda(z)}} |\langle (\partial_t - \partial_x^2)\, U(x,\cdot), \psi_x^r \rangle|\, r^{-(\gamma-2)},$$

*and* $\psi$ *is any fixed test function.*

The precise version we use is taken from [EW24, Theorem 3.1]. Since we are considering the local heat operator $(\partial_t - \partial_x^2)$ instead of the non-local fractional heat operator considered there, the last $L^\infty$ term turns into a local one. See [MW20b, Lemma 2.11] for a similar statement.

With a slight abuse of notation, we denote by $U_\gamma: \mathbb{R}^2 \times \mathbb{R}^2 \to \mathbb{R}$ the germ induced by the modelled distribution $U_\gamma \in \mathcal{D}^\gamma$ (see (2.16)) via the definition

$$U_\gamma(x,y) := u(y) - u(x) - \sum_{\tau \in \mathcal{T}_{\gamma-2}} \frac{\Upsilon_x[\tau]}{\tau!} (\Pi_x \mathcal{I}(\tau))(y) = u(y) - \Pi_x(U_\gamma(x)) + u'(x)\,(y-x)_1. \tag{3.3}$$

First, we compute the term $\mathscr{L} U_\gamma$ explicitly and, using the Reconstruction Theorem, bound this term. We postpone the proofs of all the lemmas to Section 3.6.

**Lemma 3.5**. *Let $(\Pi;\Gamma)$ be an admissible model and $U_\gamma \in \mathcal{D}^\gamma$ a coherent modelled distribution of the form (2.16) for some fixed $\gamma \in (2-\alpha, 2)$ such that $u = \mathcal{R} U_\gamma$ solves (3.1) on $B_1$. Let $U_\gamma$ be the germ defined by (3.3), then on $B_1 \subset \mathbb{R}^{1+1}$*

$$\mathscr{L} U_\gamma(x,\cdot) = -u\,|u|^{m-1} + \mathcal{R}(\hat{\sigma}(U_\gamma)\,\Xi) - \Pi_x(\hat{\sigma}(U_\gamma)\,\Xi) + \sum_{\tau \in \mathcal{T}_{\gamma-\alpha,\gamma}} \frac{\Upsilon_x[\tau\,\Xi]}{\tau!} \Pi_x(\tau\,\Xi),$$

*and for any $B_\lambda(z) \subset (0,1) \times (-1,1)$*

$$\begin{aligned}\lambda^\gamma \,\|\mathscr{L} U_\gamma\|_{\gamma-2;B_\lambda(z)} \;\lesssim\; & \lambda^2 \,\|u\|_{B_\lambda(z)}^m + \sum_{\tau \in \mathcal{T}_\gamma} \lambda^{\gamma-(|\tau|-|\mathfrak{n}(\tau)|)} \,[\hat{\sigma}(U_\gamma);\tau]_{B_\lambda(z)}\, \lambda^{\alpha+|\tau|-|\mathfrak{n}(\tau)|}\,[\Pi;\tau\,\Xi]_{B_\lambda(z)} \\ & + \sum_{\tau \in \mathcal{T}_{\gamma-\alpha,\gamma}} (\lambda\,\|u'\|_{B_\lambda(z)})^{|\mathfrak{n}(\tau)|}\,\lambda^{|\mathcal{I}(\tau\Xi)|-|\mathfrak{n}(\tau)|}\,[\Pi;\tau\,\Xi]_{B_\lambda(z)}.\end{aligned}$$

Lemma 2.5 gives us bounds on the modelled seminorms of $\hat{\sigma}(U_\gamma)$. These bounds contain linear and nonlinear terms in the generalised derivative $u'$, which we control with the following lemma.

**Lemma 3.6**. *If $c = c(\alpha,\gamma) \in (0,1)$ from Definition 3.1 is small enough, then for all $0 < r < \lambda \leqslant \lambda_t$:*

$$r\,\|u'\|_{B_\lambda(z)} \lesssim r^\gamma\,[U_\gamma;\mathbf{1}]_{B_\lambda(z)} + \|u\|_{B_\lambda(z)} + 1, \tag{3.4}$$

*with the implicit proportionality constant depending only on $\gamma$ and $\alpha$. Moreover, if the bound*

$$\|u\|_{B_\lambda(z)} + 1 < \lambda^\gamma\,[U_\gamma;\mathbf{1}]_{B_\lambda(z)}, \tag{3.5}$$

*holds, then we can conclude for any $\eta \in (1,\gamma]$ the bound*

$$(\lambda\,\|u'\|_{B_\lambda(z)})^\eta \lesssim \lambda^\gamma\,[U_\gamma;\mathbf{1}]_{B_\lambda(z)} + \|u\|_{B_\lambda(z)} + 1 + \lambda^\gamma\,[U_\gamma;\mathbf{1}]_{B_\lambda(z)}\,\|u\|_{B_\lambda(z)}^{\eta(1-\gamma^{-1})}. \tag{3.6}$$

**Remark 3.7.** If bound (3.5) does not hold, then $\lambda^\gamma\,[U_\gamma;\mathbf{1}]_{B_\lambda(z)} \leqslant \|u\|_{B_\lambda(z)} + 1$, which, in combination with Lemma 3.23 (which does not use (3.6) from Lemma 3.6), implies the result of Theorem 3.2. Therefore, we will assume without loss of generality that (3.5) holds.

**Corollary 3.8**. *If $c = c(\alpha,\gamma) \in (0,1)$ in Definition 3.1 is small enough, then for all $0<\lambda\leqslant\lambda_t$:*

$$\|U_\gamma\|_{B_\lambda(z)} \;\lesssim\; c\,\lambda^\gamma\,[U_\gamma;\mathbf{1}]_{B_\lambda(z)} + \|u\|_{B_\lambda(z)} + 1.$$

After post-processing Lemma 2.5 with the previous results, we obtain the following:

**Lemma 3.9**. *Let $c = c(\alpha,\gamma) \in (0,1)$ in Definition 3.1 be small enough and $0 < \lambda \leqslant \lambda_t$. Then we have for $\tau \in \mathcal{T}_\gamma^+$ without polynomial decorations:*

$$\begin{aligned}\lambda^{\gamma-|\tau|}\,[\hat{\sigma}(U_\gamma);\tau]_{B_\lambda(z)} \;\lesssim\; & \lambda^\gamma\,[U_\gamma;\mathbf{1}]_{B_\lambda(z)} + \|u\|_{B_\lambda(z)} + 1 \\ & + \mathbb{1}_{\gamma-|\tau|>1}\,\lambda^\gamma\,[U_\gamma;\mathbf{1}]_{B_\lambda(z)}\,\|u\|_{B_\lambda(z)}^{(\gamma-|\tau|)(1-\gamma^{-1})},\end{aligned} \tag{3.7}$$

*and for $\tau \in \mathcal{T}_\gamma^+$ with polynomial decorations:*

$$\lambda^{\gamma-(|\tau|-1)}\,[\hat{\sigma}(U_\gamma);\tau]_{B_\lambda(z)} \lesssim \lambda^\gamma\,[U_\gamma;\mathbf{1}]_{B_\lambda(z)} + \|u\|_{B_\lambda(z)} + 1 + \lambda^\gamma\,[U_\gamma]_{\gamma\text{-3pt};B_\lambda(z)} + (\lambda^\gamma\,[U_\gamma;\mathbf{1}]_{B_\lambda(z)} + \|u\|_{B_\lambda(z)})\,\|u\|_{B_\lambda(z)}^{(\gamma-|\tau|+1)(1-\gamma^{-1})}. \tag{3.8}$$

For the next result, we start to use, in addition, the smallness of $\lambda_t$, coming from Definition 3.1, defined in terms of inverse powers of $\|u\|_{P_t}$, which allows us to reduce superlinear terms in $u$ (the one coming from the superlinear damping $u\,|u|^{m-1}$ and the last terms in Lemma 3.9, with the latter being the most relevant ones).

**LEMMA 3.10**. *Let $c = c(\alpha,\gamma) \in (0,1)$ in Definition 3.1 be small enough and $0 < \lambda \leqslant \lambda_t$, then*

$$\lambda^\gamma\,\|\mathscr{L} U_\gamma\|_{\gamma-2s;B_\lambda(z)} \lesssim c\,\lambda^\gamma\,[U_\gamma;\mathbf{1}]_{B_\lambda(z)} + \lambda^\gamma\,[U_\gamma]_{\gamma\text{-3pt};B_\lambda(z)} + \|u\|_{B_\lambda(z)} + 1. \tag{3.9}$$

A control on the three-point continuity seminorm can be obtained as in [EW24, Lemma 4.10], which combined with the identity $[U_\gamma;\mathcal{I}(\tau\,\Xi)]_{B_\lambda(z)} = [\hat{\sigma}(U_{\gamma-\alpha});\tau]_{B_\lambda(z)}$ consequence of Lemma 2.6

$$\begin{aligned}
&\lambda^\gamma\,[U_\gamma]_{\gamma\text{-3pt};B_\lambda(z)}\\
\lesssim& \sum_{\rho\in\mathcal{T}_{\gamma-2}} \lambda^\gamma\,[U_\gamma;\mathcal{I}(\rho)]_{B_\lambda(z)}\,[\Gamma\,\mathcal{I}(\rho);\mathbf{1}]_{B_\lambda(z)}\\
=& \sum_{\tau\in\mathcal{T}_{\gamma-\alpha}^+} \lambda^{\gamma-|\mathcal{I}(\tau\Xi)|+|\mathfrak{n}(\mathcal{I}(\tau\Xi))|}\,[U_\gamma;\mathcal{I}(\tau\,\Xi)]_{B_\lambda(z)}\,\lambda^{|\mathcal{I}(\tau\Xi)|-|\mathfrak{n}(\mathcal{I}(\tau\Xi))|}\,[\Gamma\,\mathcal{I}(\tau\,\Xi);\mathbf{1}]_{B_\lambda(z)}\\
=& \sum_{\tau\in\mathcal{T}_{\gamma-\alpha}^+} \lambda^{\gamma-\alpha-(|\tau|-|\mathfrak{n}(\tau)|)}\,[\hat{\sigma}(U_{\gamma-\alpha});\tau]_{B_\lambda(z)}\,\lambda^{\alpha+|\tau|-|\mathfrak{n}(\tau)|}\,[\Gamma\,\mathcal{I}(\tau\,\Xi);\mathbf{1}]_{B_\lambda(z)},
\end{aligned} \tag{3.10}$$

where in the second line we used the re-indexing $\rho \leftrightarrow \tau\,\Xi$ coming from the identity of trees (2.8), and that $\alpha = |\mathcal{I}(\Xi)| = 2 + |\Xi|$. After some post-processing we can conclude the following:

**LEMMA 3.11**. *Let $c = c(\alpha,\gamma) \in (0,1)$ in Definition 3.1 be small enough and $0 < \lambda \leqslant \lambda_t$, then*

$$\lambda^\gamma\,[\Lambda_\gamma]_{(\gamma-1)\text{-3pt};B_\lambda(z)} \lesssim c\,\lambda^\gamma\,[U_\gamma;\mathbf{1}]_{B_\lambda(z)} + \|u\|_{B_\lambda(z)} + 1, \tag{3.11}$$

*Moreover, under an additional smallness condition on $c = c(\alpha,\gamma) \in (0,1)$ we can also conclude*

$$\lambda^\gamma\,[U_\gamma]_{\gamma\text{-3pt};B_\lambda(z)} \lesssim c\,\lambda^\gamma\,[U_\gamma;\mathbf{1}]_{B_\lambda(z)} + \|u\|_{B_\lambda(z)} + 1. \tag{3.12}$$

**Remark 3.12.** Since $\gamma \in (2-\alpha, 2)$ then $\gamma - \alpha \in (2 - 2\,\alpha, 2)$, and in particular we have to deal with gradients when working with $\hat{\sigma}(U_{\gamma-\alpha})$ if $2 - 2\,\alpha \geqslant 1 \Longleftrightarrow \alpha \leqslant \frac{1}{2}$, which includes the case of $(1+1)$-space-time white noise $\left(\alpha = \frac{1}{2}^-\right)$. In this case, from Lemma 3.9 we see that the terms with decoration are themselves bounded by the three-point continuity seminorm we are trying to bound. However, the small constant in (3.10) coming from the model norms and Definition 3.1 allows us to absorb the three-point continuity seminorm to the right-hand side and close the estimate.

With the previous lemmas, we can prove Theorem 3.2.

**Proof of Theorem 3.2.** From the Schauder estimate in the form of Theorem 3.4 and Lemma 3.10 we obtain for all $\lambda \in (0,\lambda_t]$ and $z \in P_{t+2\lambda}$

$$\lambda^\gamma\,[U_\gamma]_{\gamma;B_{\lambda/2}(z)} \lesssim c\,\lambda^\gamma\,[U_\gamma]_{\gamma;B_\lambda(z)} + \|u\|_{P_t} + 1 + \lambda^\gamma\,[U_\gamma]_{\gamma\text{-3pt};B_\lambda(z)} + \lambda^\gamma\,[\Lambda_\gamma]_{(\gamma-1)\text{-3pt};B_\lambda(z)} + \|U_\gamma\|_{B_\lambda(z)}.$$

This bound can be improved by Lemma 3.11 and Corollary 3.8 to obtain

$$\lambda^\gamma\,[U_\gamma]_{\gamma;B_{\lambda/2}(z)} \lesssim c\,\lambda^\gamma\,[U_\gamma]_{\gamma;B_\lambda(z)} + \|u\|_{P_t} + 1.$$

From here the proof follows analogously to [EW24, Theorem 4.2] using the absorption lemma [EW24, Lemma 4.25]. □

## 3.2. Large scales estimates

We consider a regularisation $(u)_\lambda(x) := \langle u, \psi_x^\lambda \rangle$ of the solution to (3.1) at some scale $\lambda > 0$ to be determined. This regularisation solves for any $z \in (\lambda^2, 1) \times (-1+\lambda, 1-\lambda)$ the equation

$$\begin{aligned}((\partial_t - \partial_x^2)\,(u)_\lambda)(z) &= \beta\,(u)_\lambda - (u)_\lambda\,|(u)_\lambda|^{m-1}(z) + [(u)_\lambda\,|(u)_\lambda|^{m-1} - (u\,|u|^{m-1})_\lambda](z)\\ &\quad + \langle \mathcal{R}(\hat{\sigma}(U_\gamma)\,\Xi) - \Pi_z(\hat{\sigma}(U_\gamma)(z)\,\Xi), \psi_z^\lambda \rangle + \langle \Pi_z(\hat{\sigma}(U_\gamma)(z)\,\Xi), \psi_z^\lambda \rangle. \end{aligned} \tag{3.13}$$

If the damping is strong enough, then we can replicate the argument from [MW20b] to conclude the following:

**Theorem 3.13**. *Fix $\alpha \in (0,1)$, $\beta \in \mathbb{R}$ and $m > \frac{2-\alpha}{\alpha}$. Let $\gamma \in \left(2-\alpha, 2 \wedge \left(1 + \frac{\alpha\,(m-1)}{2}\right)\right]$. Set $k := \lceil \gamma / \alpha \rceil$ and assume that $\sigma \in C_b^k(\mathbb{R})$. Let $(\Pi, \Gamma)$ be a weakly admissible model on the regularity structure described in Section 2.1, and let $U \in \mathcal{D}^\gamma$ be a coherent modelled distribution of the form (2.16) whose reconstruction $u := \mathcal{R}U$ solves weakly in $P := (0,1) \times (-1,1)$ the equation*

$$(\partial_t - \partial_x^2)\,u = \beta\,u - u|u|^{m-1} + \mathcal{R}\,(\hat{\sigma}(U)\,\Xi).$$

*Then*

$$\|u\|_{L^\infty(P_t)} \lesssim_{\alpha, |\beta|, m, \gamma, \|\sigma\|_{C_b^k}} \max\left\{ t^{-\frac{2}{m-1}}, (1 + |||(\Pi,\Gamma)|||)^{\frac{2}{\alpha(m-1) - 2(\gamma-1)}} \right\}, \qquad t \in (0,1),$$

*where $P_t := (t^2, 1) \times (-1+t, 1-t)$ and $|||(\Pi,\Gamma)|||$ is defined in (2.14). In particular, the implicit constant is independent of the model and of any initial or boundary data.*

As discussed in the introduction, the previous argument does not cover the situation of space-time white noise with cubic damping. Following [CDW26], the goal is to extract the highest super-linear in $u$ terms as transport terms of the form $(\partial_{e_1} u_\lambda)\,b_\lambda$ which lead to improved estimates on $(\partial_t - \partial_x^2)\,(u)_\lambda - (\partial_{e_1} u_\lambda)\,b_\lambda$. We first decompose the model term in (3.13). The restriction on the scales imposed in (3.14) and (3.16) ensures that every other term in this decomposition is bounded by the contribution from the noise symbol which has the worst scaling.

**Lemma 3.14**. *Under the same assumptions as Theorem 3.2. Assume that $m \leqslant (2-\alpha)\,\alpha^{-1}$. Since $\gamma > 2-\alpha$ we then have $\gamma - 1 - \frac{\alpha\,(m-1)}{2} > 0$, and therefore $\lambda_t$ from Definition 3.1 is given by $\lambda_t = \left(c\,(1 + |||(\Pi,\Gamma)|||)^{-1}\,\|u\|_{P_t}^{-(\gamma-1)}\right)^{\frac{1}{\alpha}}$. Then for any*

$$\lambda \in \left(0, \lambda_t\,\|u\|_{P_t}^{-\frac{1}{\gamma}}\right] \subset (0, \lambda_t) \tag{3.14}$$

*and every $z \in P_{t+\lambda_t}$*

$$|\langle \Pi_z(\hat{\sigma}(U_\gamma)(z)\,\Xi), \varphi_z^\lambda \rangle - (\partial_{e_1} u_\lambda)(z)\,b_\lambda^\Pi(z)| \lesssim_{\alpha,\gamma} \lambda^{\alpha-2}\,|||(\Pi,\Gamma)|||, \tag{3.15}$$

*where $b_\lambda^\Pi \colon \mathbb{R}^2 \to \mathbb{R}$ is given by (3.43).*

We have an analogous decomposition for the reconstruction term in (3.13) which relies on the multi-scale decomposition (3.51):

**Lemma 3.15**. *Under the same assumptions as Lemma 3.14 for any*

$$\lambda \in \left(0, \lambda_t\,\|u\|_{P_t}^{-\frac{2}{\gamma+1}}\right] \subset (0, \lambda_t) \tag{3.16}$$

*and every* $z\in P_{t+\lambda_t}$

$$|\langle \mathcal{R}(\hat{\sigma}(U_\gamma)\,\Xi) - \Pi_z(\hat{\sigma}(U_\gamma)(z)\,\Xi), \varphi_z^\lambda\rangle - (\partial_{e_1}u_\lambda)(z)\, b_\lambda^{\mathcal{R}}(z)| \lesssim_{\alpha,\gamma} \lambda^{\alpha-2}\, |||(\Pi,\Gamma)|||,$$

*where* $b_\lambda^{\mathcal{R}}\colon \mathbb{R}^2\to\mathbb{R}$ *is given by* (3.62)*.*

The proofs of Lemmas 3.14 and 3.15 are postponed to Section 3.7.

**COROLLARY 3.16**. *Under the same assumptions as Lemma 3.14 we set* $\lambda := \lambda_t\, \|u\|_{P_t}^{-\frac{2}{\gamma+1}} \in (0,\lambda_t)$*. Then, for* $z\in P_{t+\lambda_t}$*, we have*

$$\begin{aligned}
|\langle \Pi_z(\hat{\sigma}(U_\gamma)(z)\,\Xi), \varphi_z^\lambda\rangle - (\partial_{e_1}u_\lambda)(z)\, b_\lambda^{\Pi}(z)| &\lesssim \|u\|_{P_t}^{\left(\frac{2-\alpha}{\alpha}\right)\left(\gamma-1+\frac{2\alpha}{\gamma+1}\right)} (1+|||(\Pi,\Gamma)|||)^{\frac{2}{\alpha}},\\
|\langle \mathcal{R}(\hat{\sigma}(U_\gamma)\,\Xi) - \Pi_z(\hat{\sigma}(U_\gamma)(z)\,\Xi), \varphi_z^\lambda\rangle - (\partial_{e_1}u_\lambda)(z)\, b_\lambda^{\mathcal{R}}(z)| &\lesssim \|u\|_{P_t}^{\left(\frac{2-\alpha}{\alpha}\right)\left(\gamma-1+\frac{2\alpha}{\gamma+1}\right)} (1+|||(\Pi,\Gamma)|||)^{\frac{2}{\alpha}}.
\end{aligned}$$

*Moreover, if* $\gamma\in(2-\alpha,2)$ *is sufficiently close to* $2-\alpha$*, then* $\gamma-1+\frac{2\alpha}{\gamma+1}<1$*.*

**Proof.** Since $\gamma\in(2-\alpha,2)\subset(1,2)$ then $\gamma+1<2\gamma \Longrightarrow \frac{1}{\gamma}<\frac{2}{\gamma+1}$, and therefore

$$\|u\|_{P_t}^{-\frac{2}{\gamma+1}} \leqslant \|u\|_{P_t}^{-\frac{1}{\gamma}},$$

and therefore both Lemmas 3.14 and 3.15 can be applied with this choice of $\lambda=\lambda_t\,\|u\|_{P_t}^{-\frac{2}{\gamma+1}}$. Since

$$\begin{aligned}
\lambda^{\alpha-2}\, |||(\Pi,\Gamma)||| &= \left(c\,(1+|||(\Pi,\Gamma)|||)^{-\frac{1}{\alpha}}\, \|u\|_{P_t}^{-\frac{(\gamma-1)}{\alpha}-\frac{2}{\gamma+1}}\right)^{\alpha-2} |||(\Pi,\Gamma)|||\\
&= \|u\|_{P_t}^{\left(\frac{2-\alpha}{\alpha}\right)\left(\gamma-1+\frac{2\alpha}{\gamma+1}\right)} (1+|||(\Pi,\Gamma)|||)^{\frac{2-\alpha}{\alpha}+1}\\
&= \|u\|_{P_t}^{\left(\frac{2-\alpha}{\alpha}\right)\left(\gamma-1+\frac{2\alpha}{\gamma+1}\right)} (1+|||(\Pi,\Gamma)|||)^{\frac{2}{\alpha}}.
\end{aligned}$$

On the other hand, the map $f\colon(2-\alpha,2)\to\mathbb{R}$ given by $f(\gamma)=\gamma-1+\frac{2\alpha}{\gamma+1}$ is increasing since $f'(\gamma)=1-\frac{2\alpha}{(\gamma+1)^2}>1-\frac{2\alpha}{(3-\alpha)^2}$ and

$$\inf_{\gamma\in(2-\alpha,2)} f'(\gamma)>0 \Longleftrightarrow 1>\frac{2\alpha}{(3-\alpha)^2} \Longleftrightarrow (3-\alpha)^2>2\alpha \Longleftrightarrow \alpha^2-8\alpha+9>0,$$

and since $\alpha^2-8\alpha+9>0-8+9=1>0$ then we conclude that $f$ is increasing, and therefore

$$\lim_{\gamma\downarrow 2-\alpha} f(\gamma) = \lim_{\gamma\downarrow 2-\alpha} \gamma-1+\frac{2\alpha}{\gamma+1} = 1-\alpha+\frac{2\alpha}{3-\alpha} = 1-\alpha\left(1-\frac{2}{3-\alpha}\right).$$

Since $\alpha\in(0,1)$, then $2<3-\alpha$ and $\frac{2}{3-\alpha}<1$, and we conclude that $\lim_{\gamma\downarrow 2-\alpha} f(\gamma)<1$ which concludes the proof. □

Motivated by Lemmas 3.14 and 3.15 we can rewrite (3.13) as

$$\begin{aligned}
((\partial_t-\partial_x^2)\,(u)_\lambda)(z) &= \beta(u)_\lambda-(u)_\lambda|(u)_\lambda|^{m-1}(z)+(\partial_{e_1}u_\lambda)(z)b_\lambda(z)+[(u)_\lambda|(u)_\lambda|^{m-1}-(u|u|^{m-1})_\lambda](z)\\
&\quad+\langle \mathcal{R}(\hat{\sigma}(U_\gamma)\,\Xi) - \Pi_z(\hat{\sigma}(U_\gamma)(z)\,\Xi), \varphi_z^\lambda\rangle - (\partial_{e_1}u_\lambda)(z)\, b_\lambda^{\mathcal{R}}(z)\\
&\quad+\langle \Pi_z(\hat{\sigma}(U_\gamma)(z)\,\Xi), \varphi_z^\lambda\rangle - (\partial_{e_1}u_\lambda)(z)\, b_\lambda^{\Pi}(z), \qquad (3.17)
\end{aligned}$$

where $b_\lambda(z)=b_\lambda^{\Pi}(z)+b_\lambda^{\mathcal{R}}(z)$. The following modification of [MW20b, Lemma 2.7] takes into account the presence of the transport term in the equation and allows us to conclude an $L^\infty$-estimate on the solution that does not depend on this transport term. However, the price for removing this transport term is that our estimate needs to be global in space and therefore we restrict ourselves to periodic boundary conditions.

**Lemma 3.17**. *Let* $u: [0,1]\times\mathbb{T}\to\mathbb{R}$ *be a smooth function which satisfies pointwise in* $(0,1)\times\mathbb{T}$

$$((\partial_t-\partial_x^2)\,u)(z)=-u(z)\,|u(z)|^{m-1}+(\partial_x u)(z)\,b(z)+g(z)$$

*for some bounded function* $g: [0,1]\times\mathbb{T}\to\mathbb{R}$ *and* $m>1$. *Then for all* $(t,x)\in(0,1]\times\mathbb{T}$ *it holds that*

$$|u(t,x)|\lesssim\max\left\{t^{-\frac{1}{m-1}},\|g\|^{\frac{1}{m}}\right\}.$$

**Proof.** The proof follows the same barrier argument as [EW24, Lemma 4.17], so we only mention the differences. We consider the temporal barrier defined as

$$\eta(t)=\lambda\left(B(t)+\lambda\,\|g\|^{\frac{1}{m}}\right)^{-1}$$

for $B(t)=((m-1)\,t)^{-\frac{1}{m-1}}$. One can easily see that

$$\left(\frac{\partial_t\eta}{\eta}\right)(t)=\left(B(t)+\lambda\,\|g\|^{\frac{1}{m}}\right)^{-1}B(t)^m,$$

and therefore

$$\eta^{m-1}\,\frac{\partial_t\eta}{\eta}=\lambda^{m-1}\left(\frac{B(t)}{B(t)+\lambda\,\|g\|^{\frac{1}{m}}}\right)^m\leqslant\lambda^{m-1},$$

and in particular if $\lambda\leqslant 2^{-\frac{1}{m-1}}$ then $\frac{\partial_t\eta}{\eta}\leqslant\frac{1}{2\,\eta^{m-1}}$ which is the analogue of condition (4.70) in [EW24, Lemma 4.17] for a general damping $m>1$. Moreover, one trivially has the bound $\eta\leqslant\|g\|^{-\frac{1}{m}}$. Instead of the (fractional) heat operator we consider the differential operator $\mathscr{L}=(\partial_t-\partial_x^2-b(\cdot)\,\partial_x)$. Since $B(s)\to+\infty$ as $s\downarrow 0$, we have $\eta(s)\to 0$. The continuity of $u$ on $[0,1]\times\mathbb{T}$ therefore implies that

$$u(s,x)\,\eta(s)\longrightarrow 0$$

uniformly in $x\in\mathbb{T}$ as $s\downarrow 0$. Thus $u\,\eta$ can be continuously extended by setting $(u\,\eta)(0,\cdot)=0$. Suppose now that the positive maximum of $u\,\eta$ on $[0,1]\times\mathbb{T}$ is non-zero. Since $u\,\eta$ vanishes at $t=0$, this maximum is attained at some $x\in(0,1]\times\mathbb{T}$. Since the torus $\mathbb{T}$ has no boundary, the spatial component of $\arg\max_{y\in[0,1]\times\mathbb{T}}(u\,\eta)\,(y)$ is attained at the interior and therefore $\partial_x(u\,\eta)(x)=0$. Then

$$0\leqslant\mathscr{L}(u\,\eta)(x)=(\eta\,\mathscr{L}\,u+u\,\partial_t\,\eta\,)(x)=-\eta\,(u\,|u|^{m-1}-g(\cdot))(x)+(u\,\partial_t\,\eta)(x),$$

independently of the sign of $b$, and therefore

$$(u\,|u|^{m-1}-g)(x)\leqslant\left(u\,\frac{\partial_t\eta}{\eta}\right)(x),$$

which is the analogue of [EW24, eq. (4.69)]. One then concludes analogously that

$$|u(t,\cdot)|\lesssim\eta^{-1}(t)\sim B(t)+\|g\|^{\frac{1}{m}}\sim\max\left\{B(t),\|g\|^{\frac{1}{m}}\right\}\sim\max\left\{t^{-\frac{1}{m-1}},\|g\|^{\frac{1}{m}}\right\}.\qquad\square$$

We can conclude the following estimate on the regularisation $(u)_\lambda$, whose proof follows immediately from (3.17) and Lemma 3.17, after absorbing the lower-order term $\beta(u)_\lambda$ using

$$|\beta|\,|(u)_\lambda|\leqslant\frac{1}{4}\,|(u)_\lambda|^m+C(m)\,|\beta|^{\frac{m}{m-1}}.$$

**LEMMA 3.18**. *Let $m>1$, $t\in[0,1)$, and $0<\lambda\leqslant R'<R$ be such that $t+R<1$. Then*

$$\begin{aligned}\|(u)_\lambda\|_{P_{t+R}} \lesssim \max\Big\{&(R-R')^{-\frac{2}{m-1}}, |\beta|^{\frac{1}{m-1}}, \|(u\,|u|^{m-1})_\lambda-(u)_\lambda\,|(u)_\lambda|^{m-1}\|_{P_{t+R'}}^{\frac{1}{m}},\\ &\sup_{z\in P_{t+R'}}|\langle\mathcal{R}(\hat\sigma(U_\gamma)\,\Xi)-\Pi_z(\hat\sigma(U_\gamma)(z)\,\Xi),\varphi_z^\lambda\rangle-(\partial_{e_1}u_\lambda)(z)\,b_\lambda^{\mathcal{R}}(z)|^{\frac{1}{m}},\\ &\sup_{z\in P_{t+R'}}|\langle\Pi_z(\hat\sigma(U_\gamma)(z)\,\Xi),\varphi_z^\lambda\rangle-(\partial_{e_1}u_\lambda)(z)\,b_\lambda^{\Pi}(z)|^{\frac{1}{m}}\Big\}.\end{aligned}$$

The commutator term and removing the regularisation are dealt with as in [EW24, Lemma 4.19].

**LEMMA 3.19**. *Assume the hypotheses and conclusion of Theorem 3.2. Let $0<4\lambda\leqslant\lambda_t$. Then*

$$\begin{aligned}\|u-(u)_\lambda\|_{P_{t+\lambda_t}} &\lesssim \left(\frac{\lambda}{\lambda_t}\right)^\alpha\|u\|_{P_t}\\ \|(u|u|^{m-1})_\lambda-(u)_\lambda|(u)_\lambda|^{m-1}\|_{P_{t+\lambda_t}} &\lesssim \left(\frac{\lambda}{\lambda_t}\right)^\alpha\|u\|_{P_t}^m.\end{aligned}$$

*In particular, the same estimates hold with $P_{t+\lambda_t}$ replaced by $P_{t+R}$ whenever $R\geqslant\lambda_t$ and $t+R<1$.*

**Proof.** Let $z\in P_{t+\lambda_t}$. We apply Corollary 3.3 with regularity exponent $\eta=\alpha$ at scale ${}^{\lambda_t}\!/_2$. Since $z\in P_{t+2(\lambda_t/2)}$ and $\operatorname{supp}\varphi_z^\lambda\subset B_\lambda(z)\subset B_{\lambda_t/4}(z)$, identity (3.3) gives, for every $z'\in\operatorname{supp}\varphi_z^\lambda$,

$$|u(z)-u(z')|\leqslant[u]_{\alpha;B_\lambda(z)}\,\lambda^\alpha\lesssim\left(\frac{\lambda}{\lambda_t}\right)^\alpha\|u\|_{P_t}.$$

Integrating this estimate against $\varphi_z^\lambda$ proves the first claim. Moreover,

$$|a\,|a|^{m-1}-b\,|b|^{m-1}|\lesssim_m\max\{|a|,|b|\}^{m-1}\,|a-b|.$$

For $z'\in\operatorname{supp}\varphi_z^\lambda$, the first estimate also implies

$$|u(z')-(u)_\lambda(z)|\lesssim\left(\frac{\lambda}{\lambda_t}\right)^\alpha\|u\|_{P_t}.$$

Consequently,

$$|(u|u|^{m-1})_\lambda(z)-(u)_\lambda(z)\,|(u)_\lambda(z)|^{m-1}|\lesssim\|u\|_{P_t}^{m-1}\int|\varphi_z^\lambda(z')|\,|u(z')-(u)_\lambda(z)|\,\mathrm{d}z'\lesssim\left(\frac{\lambda}{\lambda_t}\right)^\alpha\|u\|_{P_t}^m,$$

which proves the second claim. □

## 3.3. Small times estimates

The coming-down-from-infinity estimate controls the solution away from the initial time. We complement it with the following local estimate.

**THEOREM 3.20**. *Under the same assumptions as Theorem 1.1 there exists $T_0\in(0,1]$, depending only on $\alpha,\gamma,m,|\beta|,\|\sigma\|_{C_b^k},\|u_0\|_{L^\infty(\mathbb{T})}$ and $|||(\Pi,\Gamma)|||$ such that*

$$\|u\|_{[0,T_0]\times\mathbb{T}}\lesssim(\|u(0,\,\cdot\,)\|_{L^\infty(\mathbb{T})}\vee1).$$

*Moreover, $T_0$ can be chosen uniformly when the initial-condition norm and the model norm remain bounded, with all the other parameters fixed.*

We split the solution $u:[0,1]\times\mathbb{T}\to\mathbb{R}$ as $u(t,\cdot)=e^{t\partial_x^2}u_0+w(t,\cdot)$ for $t\in[0,1]$. Then $w:[0,1]\times\mathbb{T}\to\mathbb{R}$ satisfies

$$\begin{cases} (\partial_t-\partial_x^2)\,w=\beta\,u-u\,|u|^{m-1}+\mathcal{R}(\hat{\sigma}(U)\,\Xi) & \text{in } (0,1)\times\mathbb{T} \\ w(0,\cdot)=0 \end{cases}.$$

Due to the 0 initial condition, $w$ can be continuously extended to negative times, i.e., $w:\mathbb{R}\times\mathbb{T}\to\mathbb{R}$ by just setting $w(t,\cdot)=0$ for all $t\leqslant 0$. With this we have the estimate

$$|u(t,x)|\leqslant\|P_t\,u_0\|_{\mathbb{T}}+|w(t,x)|\leqslant\|u_0\|_{\mathbb{T}}+|w(t,x)-w(0,x)|\leqslant\|u_0\|_{\mathbb{T}}+t^{\frac{\alpha}{2}}\,[w]_{\alpha;(-\infty,t]\times\mathbb{T}}$$

and therefore

$$\|u\|_{[0,T]\times\mathbb{T}}\leqslant\|u_0\|_{\mathbb{T}}+T^{\frac{\alpha}{2}}\,[w]_{\alpha;(-\infty,T]\times\mathbb{T}}. \tag{3.18}$$

The idea is to use (global) Schauder estimates (see e.g. [MW18, Lemma B.1]) to obtain

$$T^{\frac{\alpha}{2}}\,[w]_{\alpha;(-\infty,T]\times\mathbb{T}}\lesssim T^{\frac{\alpha}{2}}\,\|(\partial_t-\partial_x^2)\,w\|_{\alpha-2;(-\infty,T]\times\mathbb{T}}. \tag{3.19}$$

The distribution $(\partial_t-\partial_x^2)\,w=\beta\,u-u\,|u|^{m-1}+\mathcal{R}(\hat{\sigma}(U)\,\Xi)$ is initially only defined as a distribution on positive times. This can be extended to negative times with an explicit estimate of its $\mathcal{C}^{\alpha-2}$-norm. This is the content of [Hai14, Proposition 6.9], and by keeping track of the precise seminorms in that proof one can conclude the following:

**Lemma 3.21**. *Let $f\in\mathcal{C}^{\alpha-2}((-\infty,1)\times\mathbb{T})$ be the extension of $\beta\,u-u\,|u|^{m-1}+\mathcal{R}(\hat{\sigma}(U)\,\Xi)$ to negative times, then for any $T\in(0,1)$ satisfying*

$$\sqrt{T}\leqslant(1\vee\|u\|_{[0,T]\times\mathbb{T}})^{-\frac{m-1}{2}}\wedge(c\,(1+|\!|\!|(\Pi,\Gamma)|\!|\!|)^{-1}\,(1\vee\|u\|_{[0,T]\times\mathbb{T}})^{-(\gamma-1)})^{1/\alpha},$$

*with $c\in(0,1)$ chosen as in Theorem 3.2. Then*

$$\begin{aligned} &\|f\|_{\alpha-2;(-\infty,T]\times\mathbb{T}} \\ \lesssim\;& T^{1-\alpha/2}\,(1+\|u\|_{[0,T]\times\mathbb{T}}^m)+\sum_{\substack{\tau\in\mathcal{T}_\gamma\\ \tau\Xi\in\mathcal{T}}}T^{\frac{1}{2}(|\tau|-|\mathfrak{n}(\tau)|)}\left(1+\|u\|_{[0,T]\times\mathbb{T}}^{|\mathfrak{n}(\tau)|}\right)[\Pi;\tau\,\Xi] \\ &+\sum_{\substack{\tau\in\mathcal{T}_\gamma\\ \tau\Xi\in\mathcal{T}}}T^{\frac{1}{2}(|\tau|-|\mathfrak{n}(\tau)|)}[\Pi;\tau\Xi]\left(1+\|u\|_{[0,T]\times\mathbb{T}}^{1+\mathbb{1}_{\gamma-(|\tau|-|\mathfrak{n}(\tau)|)>1}(\gamma-(|\tau|-|\mathfrak{n}(\tau)|))(1-\gamma^{-1})}\right). \end{aligned}$$

**Proof.** To simplify notation we may assume that $\|u\|_{[0,T]\times\mathbb{T}}\geqslant 1$. The extension $f\in\mathcal{C}^{\alpha-2}$ of $\beta\,u-u\,|u|^{m-1}+\mathcal{R}(\hat{\sigma}(U)\,\Xi)$ is well-defined by [Hai14, Proposition 6.9]. Fix $z\in(-\infty,T]\times\mathbb{T}$ and $\lambda\in(0,1]$, the goal is to estimate $|\langle f,\psi_z^\lambda\rangle|\,\lambda^{-(\alpha-2)}$. First, we observe that, doing a change of variable in the spatial coordinate:

$$\|\psi_z^\lambda\|_{L^1}=\lambda^{-2}\int_{0\vee(z_0-\lambda^2)}^{z_0}\int_{\mathbb{R}}|\psi(\lambda^{-2}\,(s-z_0),y)|\,\mathrm{d}y\mathrm{d}s\leqslant\lambda^{-2}\,(z_0\wedge\lambda^2)\sup_{s\in\mathbb{R}}\int_{\mathbb{R}}|\psi(s,y)|\,\mathrm{d}y\lesssim\lambda^{-2}\,T\wedge 1,$$

and therefore, since the function $\lambda\mapsto\lambda^{-\alpha}\,T\wedge\lambda^{2-\alpha}$ maximises at $\lambda^{-\alpha}\,T=\lambda^{2-\alpha}\Longleftrightarrow\lambda=\sqrt{T}\in(0,1)$,

$$|\langle u\,|u|^{m-1},\varphi_z^\lambda\rangle|\,\lambda^{-(\alpha-2)}\leqslant\lambda^{2-\alpha}\,\|u\|_{[0,T]\times\mathbb{T}}^m\,\|\varphi_z^\lambda\|_{L^1}\leqslant(\lambda^{-\alpha}\,T\wedge\lambda^{2-\alpha})\,\|u\|_{[0,T]\times\mathbb{T}}^m\leqslant T^{1-\alpha/2}\,\|u\|_{[0,T]\times\mathbb{T}}^m.$$

For the linear term $\beta\,u$ we have

$$\lambda^{2-\alpha}|\langle\beta\,u,\varphi_z^\lambda\rangle|\lesssim|\beta|\,(\lambda^{-\alpha}\,T\wedge\lambda^{2-\alpha})\|u\|_{[0,T]\times\mathbb{T}}\lesssim_{|\beta|}T^{1-\alpha/2}\,(1+\|u\|_{[0,T]\times\mathbb{T}}^m).$$

For the reconstruction term $\mathcal{R}(\hat{\sigma}(U)\,\Xi)$ we only consider the case $z_0 \leqslant 2\,\lambda^2$, since otherwise everything is localised away from $t=0$ and the argument follows by a similar argument without the need of a partition of unity. Following the notation of [Hai14, Proposition 6.9] we consider the partition of unity

$$1 = \sum_{n\in\mathbb{N}} \sum_{w\in\Xi_P^n} \varphi_{w,n}(y), \qquad \forall\, y\in(\mathbb{R}\setminus\{0\})\times\mathbb{R}$$

and rewrite the test function $\psi_z^\lambda$ as

$$\psi_z^\lambda = \sum_{n\geqslant n_0} \sum_{w\in\Xi_P^n} \psi_z^\lambda\,\varphi_{w,n} = \sum_{n\geqslant n_0} \sum_{w\in\Xi_P^n} \lambda^{-|\mathfrak{s}|}\,2^{-n|\mathfrak{s}|}\,\chi_{n,zw} = \lambda^{-|\mathfrak{s}|} \sum_{n\geqslant n_0} 2^{-n|\mathfrak{s}|} \sum_{w\in\Xi_P^n} \sum_{j=1}^{M} \chi_{n,zw}^{(j)},$$

where in our case the parabolic scaling is $|\mathfrak{s}| = |(2,1)| = 3$. We start by writing

$$\begin{aligned} |\langle \mathcal{R}(\hat{\sigma}(U)\,\Xi), \psi_z^\lambda\rangle| &\leqslant \lambda^{-|\mathfrak{s}|} \sum_{n\geqslant n_0} 2^{-n|\mathfrak{s}|} \sum_{w\in\Xi_P^n} \sum_{j=1}^{M} \big|\big\langle \mathcal{R}(\hat{\sigma}(U)\,\Xi), \chi_{n,zw}^{(j)}\big\rangle\big| \\ &\leqslant \lambda^{-|\mathfrak{s}|} \sum_{n\geqslant n_0} 2^{-n|\mathfrak{s}|} \sum_{w\in\Xi_P^n} \sum_{j=1}^{M} \big|\big\langle \mathcal{R}(\hat{\sigma}(U)\,\Xi) - \Pi_{z_j}(\hat{\sigma}(U_\gamma)\,\Xi), \chi_{n,zw}^{(j)}\big\rangle\big| \\ &\quad + \lambda^{-|\mathfrak{s}|} \sum_{n\geqslant n_0} 2^{-n|\mathfrak{s}|} \sum_{w\in\Xi_P^n} \sum_{j=1}^{M} \big|\big\langle \Pi_{z_j}(\hat{\sigma}(U_\gamma)\,\Xi), \chi_{n,zw}^{(j)}\big\rangle\big|, \end{aligned}$$

where $z_j$ is such that $\mathrm{supp}\big(\chi_{n,zw}^{(j)}\big) \subset B_{\frac{1}{2}d(z_j,P)}(z_j)$ and $d(z_j,P) = \sqrt{|(z_j)_0|} \sim 2^{-n}$. Moreover, since $f\equiv 0$ in $(-\infty,0)\times\mathbb{T}$ we assume that $(z_j)_0 > 0$. We start by controlling the model term as

$$\begin{aligned} \big|\big\langle \Pi_{z_j}(\hat{\sigma}(U_\gamma)\,\Xi), \chi_{n,zw}^{(j)}\big\rangle\big| &\lesssim \sum_{\substack{\tau\in\mathcal{T}_\gamma\\ \tau\Xi\in\mathcal{T}}} |\Upsilon_{z_j}[\tau]|\,\big|\big\langle \Pi_{z_j}(\tau\,\Xi), \chi_{n,zw}^{(j)}\big\rangle\big| \\ &\lesssim \sum_{\substack{\tau\in\mathcal{T}_\gamma\\ \tau\Xi\in\mathcal{T}}} |u'(z_j)|^{|\mathfrak{n}(\tau)|}\,[\Pi;\tau\,\Xi]\,2^{-n|\tau\Xi|} \\ &\lesssim \sum_{\substack{\tau\in\mathcal{T}_\gamma\\ \tau\Xi\in\mathcal{T}}} \big((z_j)_0^{-1/2}\,\|u\|_{[0,T]\times\mathbb{T}}\big)^{|\mathfrak{n}(\tau)|}\,[\Pi;\tau\,\Xi]\,2^{-n|\tau\Xi|} \\ &\lesssim \sum_{\substack{\tau\in\mathcal{T}_\gamma\\ \tau\Xi\in\mathcal{T}}} (z_j)_0^{-\frac{1}{2}|\tau|}\,2^{-n|\tau\Xi|}\,z_0^{\frac{1}{2}(|\tau|-|\mathfrak{n}(\tau)|)}\,\|u\|_{[0,T]\times\mathbb{T}}^{|\mathfrak{n}(\tau)|}\,[\Pi;\tau\,\Xi] \\ &\lesssim \sum_{\substack{\tau\in\mathcal{T}_\gamma\\ \tau\Xi\in\mathcal{T}}} (z_j)_0^{-\frac{1}{2}|\tau|}\,2^{-n|\tau\Xi|}\,T^{\frac{1}{2}(|\tau|-|\mathfrak{n}(\tau)|)}\,\|u\|_{[0,T]\times\mathbb{T}}^{|\mathfrak{n}(\tau)|}\,[\Pi;\tau\,\Xi] \end{aligned}$$

where we used that $|u'(z)| \lesssim (z_j)_0^{-1/2}\,\|u\|_{[0,T]\times\mathbb{T}}$ by a trivial modification of Corollary 3.3 where one replaces $\|u\|_{P_t}$ by $\|u\|_{[0,T]\times\mathbb{T}}$. Since $(z_j)_0^{-\frac{1}{2}|\tau|} \sim 2^{n|\tau|}$ and $\#\Xi_P^n \sim (2^n\,\lambda)^d$ for $d=1$, therefore

$$\begin{aligned} \lambda^{-|\mathfrak{s}|} \sum_{n,w,j} 2^{-n|\mathfrak{s}|}\,(z_j)_0^{-\frac{1}{2}|\tau|}\,2^{-n|\tau\Xi|} &\sim \lambda^{-|\mathfrak{s}|} \sum_{n,w,j} 2^{-n|\mathfrak{s}|}\,2^{n|\tau|}\,2^{-n|\tau\Xi|} \\ &\lesssim \lambda^{-|\mathfrak{s}|} \sum_{n\geqslant n_0} 2^{-n|\mathfrak{s}|}\,2^{n|\tau|}\,2^{-n|\tau\Xi|}\,(2^n\,\lambda)^d \\ &= \lambda^{d-|\mathfrak{s}|} \sum_{n\geqslant n_0} 2^{-n(|\mathfrak{s}|+|\Xi|-d)} \\ &\sim \lambda^{d-|\mathfrak{s}|}\,2^{-n_0(|\Xi|+|\mathfrak{s}|-d)} \sim \lambda^{|\Xi|} = \lambda^{\alpha-2}, \end{aligned}$$

since $|\Xi|+|\mathfrak{s}|-d=\alpha>0$, which makes the geometric sum summable, and $n_0$ was such that $2^{-n_0}\sim\lambda$. We conclude that

$$\lambda^{-|\mathfrak{s}|}\sum_{n\geqslant n_0}2^{-n|\mathfrak{s}|}\sum_{w\in\Xi_P^n}\sum_{j=1}^{M}\left|\left\langle\Pi_{z_j}(\hat{\sigma}(U_\gamma)\,\Xi),\chi_{n,zw}^{(j)}\right\rangle\right|\lesssim\lambda^{\alpha-2}\sum_{\tau\in\mathcal{T}_\gamma^+}T^{\frac{1}{2}(|\tau|-|\mathfrak{n}(\tau)|)}\,\|u\|_{[0,T]\times\mathbb{T}}^{|\mathfrak{n}(\tau)|}\,[\Pi;\tau\,\Xi].$$

For the reconstruction term we have

$$\left|\left\langle\mathcal{R}(\hat{\sigma}(U)\,\Xi)-\Pi_{z_j}(\hat{\sigma}(U_\gamma)\,\Xi),\chi_{n,zw}^{(j)}\right\rangle\right|\lesssim2^{-n(\gamma+|\Xi|)}\sum_{\tau\in\mathcal{T}_\gamma^+}[\hat{\sigma}(U_\gamma);\tau]_{B_{\frac{1}{2}d(z_j,P)}(z_j)}\,[\Pi;\tau\,\Xi].$$

By the same modification of Corollary 3.3 we have

$$\begin{aligned}[\hat{\sigma}(U_\gamma);\tau]_{B_{\frac{1}{2}d(z_j,P)}(z_j)} &\lesssim (z_j)_0^{-\frac{1}{2}(\gamma-(|\tau|-|\mathfrak{n}(\tau)|))}\,\|u\|_{[0,T]\times\mathbb{T}}^{1+\mathbb{1}_{\gamma-(|\tau|-|\mathfrak{n}(\tau)|)>1}(\gamma-(|\tau|-|\mathfrak{n}(\tau)|))(1-\gamma^{-1})}\\ &\lesssim T^{\frac{1}{2}(|\tau|-|\mathfrak{n}(\tau)|)}(z_j)_0^{-\frac{\gamma}{2}}\,\|u\|_{[0,T]\times\mathbb{T}}^{1+\mathbb{1}_{\gamma-(|\tau|-|\mathfrak{n}(\tau)|)>1}(\gamma-(|\tau|-|\mathfrak{n}(\tau)|))(1-\gamma^{-1})}\end{aligned}$$

and analogously as with the model term we obtain

$$\lambda^{-|\mathfrak{s}|}\sum_{n,w,j}2^{-n|\mathfrak{s}|}\,(z_j)_0^{-\frac{\gamma}{2}}\,2^{-n(\gamma+|\Xi|)}\lesssim\lambda^{d-|\mathfrak{s}|}\sum_{n\geqslant n_0}2^{-n(|\mathfrak{s}|+|\Xi|-d)}\sim\lambda^{\alpha-2},$$

and therefore

$$\begin{aligned}&\lambda^{-|\mathfrak{s}|}\sum_{n\geqslant n_0}2^{-n|\mathfrak{s}|}\sum_{w\in\Xi_P^n}\sum_{j=1}^{M}\left|\left\langle\mathcal{R}(\hat{\sigma}(U)\,\Xi)-\Pi_{z_j}(\hat{\sigma}(U_\gamma)\,\Xi),\chi_{n,zw}^{(j)}\right\rangle\right|\\ \lesssim\;&\lambda^{\alpha-2}\sum_{\tau\in\mathcal{T}_\gamma^+}T^{\frac{1}{2}(|\tau|-|\mathfrak{n}(\tau)|)}\,[\Pi;\tau\,\Xi]\,\|u\|_{[0,T]\times\mathbb{T}}^{1+\mathbb{1}_{\gamma-(|\tau|-|\mathfrak{n}(\tau)|)>1}(\gamma-(|\tau|-|\mathfrak{n}(\tau)|))(1-\gamma^{-1})},\end{aligned}$$

which concludes the result. □

**Proof of Theorem 3.20.** The proof follows the same short-time bootstrap and absorption argument as in the proof of [CDW26, Theorem 5]. Choose $T_0\in(0,1]$ sufficiently small so that, for every $T\in(0,T_0]$ satisfying

$$\|u\|_{[0,T]\times\mathbb{T}}\leqslant3\,(1+\|u_0\|_{L^\infty(\mathbb{T})}),$$

the scale condition of Lemma 3.21 holds and, together with (3.19), its estimate gives

$$T^{\alpha/2}\,[w]_{\alpha;(-\infty,T]\times\mathbb{T}}\leqslant1+\|u_0\|_{L^\infty(\mathbb{T})}.$$

Such a choice is possible since all powers of $T$ arising from Lemma 3.21 are strictly positive. Hence, by (3.18),

$$\|u\|_{[0,T]\times\mathbb{T}}\leqslant2\,(1+\|u_0\|_{L^\infty(\mathbb{T})}).$$

The usual continuity argument closes the bootstrap and proves the result. The choice of $T_0$ is uniform when the initial-condition and model norms remain bounded. □

## 3.4. Proof of a priori estimates under a strong damping

**Proof of Theorem 3.13.** Following [EW24, Remark 4.15] we may assume w.l.o.g. that $\|u\|_{P_t}\geqslant1$ for all $t\in(0,1)$. Since $m>\frac{2-\alpha}{\alpha}$ then by (1.7) we can fix $\gamma\in(2-\alpha,2)$ close enough to $2-\alpha$ such that $\frac{m-1}{2}-\frac{\gamma-1}{\alpha}>0$. In this case we have that

$$\lambda_t=\|u\|_{P_t}^{-\frac{(m-1)}{2}}\Longleftrightarrow c^{-1}\,(1+|||(\Pi,\Gamma)|||)\leqslant\|u\|_{P_t}^{\frac{\alpha(m-1)}{2}-(\gamma-1)}. \tag{3.20}$$

If (3.20) does not hold, then

$$\|u\|_{P_t} \lesssim (1 + \|(\Pi, \Gamma)\|)^{\frac{2}{\alpha(m-1)-2(\gamma-1)}},$$

which is part of the stated bound in the result, and therefore we will assume that (3.20) holds. By [MW20a, Theorem 4.2] applied to (3.13) with $\lambda = k^{-1}\lambda_t$ for some $k \geqslant 4$ we can conclude analogously to Lemma 3.18 that

$$\begin{aligned} \|u\|_{P_{t+R}} &\lesssim \max\{\|u-(u)_\lambda\|, \|(u)_\lambda\|_{P_{t+R}}\} \\ &\lesssim \max\Big\{(k-1)^{-\alpha}\,\|u\|_{P_t}, (R-\lambda_t)^{-\frac{2}{m-1}}, k^{-\frac{\alpha}{m}}\,\|u\|_{P_t}, \sup_{z\in P_{t+R'}} |\langle \Pi_x(\hat{\sigma}(U_\gamma)(z)\,\Xi), \varphi_z^\lambda\rangle|^{\frac{1}{m}}, \\ &\qquad \sup_{z\in P_{t+R'}} |\langle \mathcal{R}(\hat{\sigma}(U_\gamma)\,\Xi) - \Pi_z(\hat{\sigma}(U_\gamma)(z)\,\Xi), \varphi_z^\lambda\rangle|^{\frac{1}{m}}\Big\} \end{aligned}$$

By (2.17) and (3.25) we have

$$\begin{aligned} \lambda^2\,|\langle \Pi_x(\hat{\sigma}(U_\gamma)(z)\,\Xi), \varphi_z^\lambda\rangle| &\lesssim \lambda^2 \sum_{\tau\in\mathcal{T}_\gamma^+} |u'(x)|^{|\mathfrak{n}(\tau)|}\,\lambda^{|\tau\Xi|}\,[\Pi;\tau\,\Xi] \\ &\lesssim_k \sum_{\tau\in\mathcal{T}_\gamma^+} \|u\|_{P_t}^{|\mathfrak{n}(\tau)|}\,\lambda^{2+|\tau\Xi|-|\mathfrak{n}(\tau)|}\,[\Pi;\tau\,\Xi] \\ &\lesssim c\,\|u\|_{P_t}, \end{aligned}$$

and therefore

$$|\langle \Pi_x(\hat{\sigma}(U_\gamma)(z)\,\Xi), \varphi_z^\lambda\rangle| \lesssim c\,\lambda^{-2}\,\|u\|_{P_t} \sim_k c\,\|u\|_{P_t}^m.$$

Moreover, by (3.28) combined with (3.33) and (3.36) (which itself needs Lemma 3.11) we have that

$$\lambda^2\,|\langle \mathcal{R}(\hat{\sigma}(U_\gamma)\,\Xi) - \Pi_z(\hat{\sigma}(U_\gamma)(z)\,\Xi), \varphi_z^\lambda\rangle| \lesssim c\,\|u\|_{P_t}$$

and therefore

$$\|\mathcal{R}(\hat{\sigma}(U_\gamma)\,\Xi) - \Pi_z(\hat{\sigma}(U_\gamma)(z)\,\Xi)\| \lesssim c\,\lambda^{-2}\,\|u\|_{P_t} \sim_k c\,\lambda_t^{-2}\,\|u\|_{P_t} = c\,\|u\|_{P_t}^m.$$

We conclude that

$$\|u\|_{P_{t+R}} \lesssim \max\Big\{(R-\lambda_t)^{-\frac{2}{m-1}}, \max\Big\{(k-1)^{-\alpha}, k^{-\frac{\alpha}{m}}, c^{\frac{1}{m}}\,C(k)\Big\}\,\|u\|_{P_t}\Big\}.$$

Choosing $k = k(\alpha, \gamma) \in \mathbb{Z}^+$ big enough and imposing a smallness condition on $c \in (0,1)$ depending on $C(k) = C(\alpha, \gamma)$ we can conclude that

$$\|u\|_{P_{t+R}} \lesssim \max\Big\{(R-\lambda_t)^{-\frac{2}{m-1}}, \frac{1}{2}\,\|u\|_{P_t}\Big\}.$$

This is the analogue of [MW20b, eq. (4.32)] and one can prove analogously that $\|u\|_{P_t} \lesssim t^{-\frac{2}{m-1}}$. □

## 3.5. Proof of a priori estimates under a weaker damping

We can now proceed with the proof of the main result Theorem 1.1.

**Proof of Theorem 1.1.** Following [EW24, Remark 4.15] we can assume w.l.o.g. that $\|u\|_{P_t} \geqslant 1$ for all $t \in (0,1)$. As in the proof of Corollary 3.16 we have

$$\varepsilon_\alpha := \lim_{\gamma\downarrow 2-\alpha}\left(\gamma - 1 + \frac{2\,\alpha}{\gamma+1}\right) = 1 - \alpha\left(1 - \frac{2}{3-\alpha}\right) \in (0,1). \tag{3.21}$$

The function $f_\alpha(\gamma) := \gamma - 1 + \frac{2\alpha}{\gamma+1}$ is continuous and increasing on the interval $(2-\alpha, 2)$, and satisfies $\lim_{\gamma \downarrow 2-\alpha} f_\alpha(\gamma) = \varepsilon_\alpha$. Since $m > (2-\alpha)\,\alpha^{-1}\varepsilon_\alpha$, we can fix $\gamma > 2-\alpha$ sufficiently close to $2-\alpha$ and depending on $(\alpha, m)$, such that $m > \frac{(2-\alpha)}{\alpha} f_\alpha(\gamma)$. Fix $c = c(\alpha,\gamma) \in (0,1)$ small enough such that Theorem 3.2 holds, later we will impose additional smallness conditions on it. Recalling Definition 3.1 we have

$$\begin{aligned}\lambda_t = \big(c\,(1 + |||(\Pi,\Gamma)|||)^{-1}\,\|u\|_{P_t}^{-(\gamma-1)}\big)^{\frac{1}{\alpha}} &\iff \big(c\,(1 + |||(\Pi,\Gamma)|||)^{-1}\,\|u\|_{P_t}^{-(\gamma-1)}\big)^{\frac{1}{\alpha}} \leqslant \|u\|_{P_t}^{-\frac{(m-1)}{2}} \\ &\iff c\,(1 + |||(\Pi,\Gamma)|||)^{-1} \leqslant \|u\|_{P_t}^{\gamma-1-\frac{\alpha(m-1)}{2}}. \end{aligned} \tag{3.22}$$

Since $m \leqslant (2-\alpha)\,\alpha^{-1}$ and $\gamma > 2-\alpha$, we have $\frac{\alpha\,(m-1)}{2} \leqslant 1-\alpha < \gamma - 1$, and $\gamma - 1 - \frac{\alpha\,(m-1)}{2} > 0$. Moreover, since $c \in (0,1)$, $c\,(1 + \|(\Pi,\Gamma)\|)^{-1} \leqslant 1$, whereas $\|u\|_{P_t}^{\gamma-1-\frac{\alpha(m-1)}{2}} \geqslant 1$. Consequently (3.22) holds and therefore

$$\lambda_t = \big(c\,(1 + \|(\Pi,\Gamma)\|)^{-1}\,\|u\|_{P_t}^{-(\gamma-1)}\big)^{\frac{1}{\alpha}},$$

which is part of the assumptions in Lemmas 3.14 and 3.15 and Corollary 3.16. We now fix $k \geqslant 4$ and use the scale

$$\lambda := k^{-1}\,\lambda_t \|u\|_{P_t}^{-2/(\gamma+1)}.$$

In particular, $4\,\lambda \leqslant \lambda_t$, and $\lambda$ lies in the ranges of both Lemmas 3.14 and 3.15. Applying Lemma 3.19 at the scale $\lambda$, we obtain

$$\|u - u_\lambda\|_{P_{t+R}} \lesssim k^{-\alpha} \|u\|_{P_t}^{1-\frac{2\alpha}{\gamma+1}} \leqslant k^{-\alpha}\,\|u\|_{P_t}, \tag{3.23}$$

and

$$|u_\lambda |u_\lambda|^{m-1} - (u|u|^{m-1})_\lambda|_{P_{t+R}}^{1/m} \lesssim k^{-\alpha/m} \|u\|_{P_t}^{1-\frac{2\alpha}{m(\gamma+1)}} \leqslant k^{-\alpha/m}\,\|u\|_{P_t}.$$

By Corollary 3.16 applied at this smaller admissible scale, the $m$-th roots of the two transport remainders are bounded by

$$C(k)\,(1 + \|(\Pi,\Gamma)\|)^{\frac{2}{\alpha m}} \|u\|_{P_t}^{\frac{2-\alpha}{\alpha m}\left(\gamma - 1 + \frac{2\alpha}{\gamma+1}\right)}.$$

We now apply Lemma 3.18 to $u_\lambda$, with $R' = \lambda_t$, and use (3.23) to remove the regularisation. For every $R > \lambda_t$, this gives

$$\begin{aligned} &\|u\|_{P_{t+R}} \\ \lesssim\;& \max\{\|u - (u)_\lambda\|_{P_{t+R}}, \|(u)_\lambda\|_{P_{t+R}}\} \\ \lesssim\;& \max\left\{k^{-\alpha}\|u\|_{P_t}, k^{-\frac{\alpha}{m}}\|u\|_{P_t}, C(k)(1 + |||(\Pi,\Gamma)|||)^{\frac{2}{\alpha m}} \|u\|_{P_t}^{\frac{(2-\alpha)}{\alpha m}\left(\gamma-1+\frac{2\alpha}{\gamma+1}\right)}, |\beta|^{\frac{1}{m-1}}, (R-\lambda_t)^{-\frac{2}{m-1}}\right\}. \end{aligned}$$

Choosing $k = k(\alpha,\gamma) \in \mathbb{Z}$ big enough depending on the implicit proportionality constant we can guarantee that

$$\|u\|_{P_{t+R}} \leqslant \max\left\{\frac{1}{2}\,\|u\|_{P_t}, C(\alpha,\gamma)\,(1 + |||(\Pi,\Gamma)|||)^{\frac{2}{\alpha m}}\,\|u\|_{P_t}^{\frac{(2-\alpha)}{\alpha m}\left(\gamma-1+\frac{2\alpha}{\gamma+1}\right)}, |\beta|^{\frac{1}{m-1}}, (R-\lambda_t)^{-\frac{2}{m-1}}\right\}.$$

For simplicity we write $\theta = \frac{(2-\alpha)}{\alpha\,m}\left(\gamma - 1 + \frac{2\,\alpha}{\gamma+1}\right)$, then $\theta \in (0,1)$ by the choice of $\gamma \in (2-\alpha, 2)$, and therefore

$$\|u\|_{P_t} \geqslant \left(2\,C(\alpha,\gamma)\,(1 + |||(\Pi,\Gamma)|||)^{\frac{2}{\alpha m}}\right)^{\frac{1}{1-\theta}} \iff C(\alpha,\gamma)\,(1 + |||(\Pi,\Gamma)|||)^{\frac{2}{\alpha m}}\,\|u\|_{P_t}^{\theta} \leqslant \frac{1}{2}\,\|u\|_{P_t}. \tag{3.24}$$

If (3.24) does not hold, then we can conclude the bound

$$\|u\|_{P_t} \leqslant \left(2\,C(\alpha,\gamma)\,(1 + |||(\Pi,\Gamma)|||)^{\frac{2}{\alpha m}}\right)^{\frac{1}{1-\theta}} \sim (1 + |||(\Pi,\Gamma)|||)^{\frac{2}{\alpha m(1-\theta)}},$$

which is part of the stated bound in (1.5). On the other hand, if (3.24) holds, then

$$\|u\|_{P_{t+R}} \leqslant \max\left\{C\,(R-\lambda_t)^{-\frac{2}{m-1}}, \frac{1}{2}\,\|u\|_{P_t}\right\}.$$

This is the analogue of [MW20b, eq. (4.32)] and from here one can proceed as in there to conclude $\|u\|_{P_t} \lesssim t^{-\frac{2}{m-1}}$, or equivalently,

$$\|u\|_{(t,1]\times\mathbb{T}} \lesssim t^{-\frac{1}{m-1}}.$$

To incorporate the initial solution one uses Theorem 3.20. □

## 3.6. Proofs of Lemmas from Section 3.1

**Remark 3.22.** It can be easily seen that for any $\tau \in \mathcal{T}$ without a noise and with trivial polynomial decorations its homogeneity only depends on its number of noises $\mathfrak{C}(\tau)$. More precisely, for $\tau \in \mathcal{T}$ without a noise at the root one has $|\tau| - |\mathfrak{n}(\tau)| = \alpha\,\mathfrak{C}(\tau)$. For any $\lambda \in (0, \lambda_t]$ and $\tau \in \mathcal{T}_{\gamma-2}$ one has that $2 + |\tau| - |\mathfrak{n}(\tau)| = |\mathcal{I}(\tau)| - |\mathfrak{n}(\mathcal{I}(\tau))| = \alpha\,\mathfrak{C}(\mathcal{I}(\tau)) = \alpha\,\mathfrak{C}(\tau)$, and therefore

$$\begin{aligned} \lambda^{2+|\tau|-|\mathfrak{n}(\tau)|}\,[\Pi;\tau] \vee \lambda^{|\mathcal{I}(\tau)|-|\mathfrak{n}(\tau)|}\,[\Gamma\,\mathcal{I}(\tau);\boldsymbol{X}^k] &\leqslant \lambda_t^{\alpha\,\mathfrak{C}(\tau)}\,|||(\Pi,\Gamma)|||^{\mathfrak{C}(\tau)} \\ &\leqslant \left(c\,|||(\Pi,\Gamma)|||^{-1}\,\|u\|_{P_t}^{-(\gamma-1)}\right)^{\mathfrak{C}(\tau)}\,|||(\Pi,\Gamma)|||^{\mathfrak{C}(\tau)} \\ &\leqslant c^{\mathfrak{C}(\tau)}\,\|u\|_{P_t}^{-(\gamma-1)\,\mathfrak{C}(\tau)} \\ &\leqslant c, \end{aligned} \tag{3.25}$$

since $\mathfrak{C}(\tau) \geqslant 1, c \in (0,1), \gamma > 1$ and $\|u\|_{P_t} \geqslant 1$. More generally, given $\tau \in \mathcal{T}^+$ let $\beta \in \mathbb{N}^{\mathcal{T}}$ be such that $\tau = \boldsymbol{X}^m\,\mathfrak{a}(\beta)$ for some $m \in \mathbb{N}^2$ (see (2.11)), then if $\beta(\rho) = 0$ for all $\rho \notin \mathcal{T}_{\gamma-2}$, the previous bound (3.25) and multiplicativity in $\mathcal{T}^+$ implies that

$$\lambda^{|\tau|-|\mathfrak{n}(\tau)|}\,[\Gamma\,\tau;\boldsymbol{X}^k] = \prod_{\rho\in\mathcal{T}} (\lambda^{|\mathcal{I}(\rho)|-|\mathfrak{n}(\rho)|}\,[\Gamma\,\mathcal{I}(\rho);\boldsymbol{X}^k])^{\beta(\rho)} \leqslant c^{\sum_\rho \beta(\rho)} \leqslant c. \tag{3.26}$$

Observe that the previous holds even though $|\tau|$ might be bigger than $\gamma$.

On the other hand, given $\tau \in \mathcal{T}_\gamma^+$ we have that $\tau\,\Xi \in \mathcal{T}_{\gamma+\alpha-2}$, and therefore

$$\begin{aligned} &\|u\|_{P_t}^{(\gamma-(|\tau|-|\mathfrak{n}(\tau)|))(1-\gamma^{-1})}\,\lambda_t^{\alpha+|\tau|-|\mathfrak{n}(\tau)|}\,|||(\Pi,\Gamma)|||^{\mathfrak{C}(\tau\Xi)} \\ = {}& \|u\|_{P_t}^{(\gamma-(2-\alpha+|\tau\Xi|-|\mathfrak{n}(\tau\Xi)|))(1-\gamma^{-1})}\,\lambda_t^{2+|\tau\Xi|-|\mathfrak{n}(\tau\Xi)|}\,|||(\Pi,\Gamma)|||^{\mathfrak{C}(\tau\Xi)} \\ = {}& \|u\|_{P_t}^{(\gamma-\alpha(\mathfrak{C}(\tau\Xi)-1))(1-\gamma^{-1})}\,\lambda_t^{\alpha\,\mathfrak{C}(\tau\Xi)}\,|||(\Pi,\Gamma)|||^{\mathfrak{C}(\tau\Xi)} \\ \leqslant {}& \|u\|_{P_t}^{\gamma-1}\,c^{\mathfrak{C}(\tau\Xi)}\,\|u\|_{P_t}^{-(\gamma-1)\,\mathfrak{C}(\tau)} \\ \leqslant {}& c \end{aligned} \tag{3.27}$$

where we used that $\mathfrak{C}(\tau) \geqslant 1$ and the function $x \mapsto (\gamma - \alpha\,(x-1))\,(1-\gamma^{-1})$ is decreasing.

**Proof of Lemma 3.5.** By (2.17) only trees $\tau \in \mathcal{T}_\gamma$ without a noise at the root contribute; by subcriticality, such trees satisfy $|\tau| \geqslant 0$. Applying $\Pi_x$ and the heat operator $(\partial_t - \partial_x^2)$ to (2.16), and using the weak admissibility condition from Definition 2.1 we have for any $x \in \mathbb{R}^2$

$$\begin{aligned} (\partial_t - \Delta)\,\Pi_x(U_\gamma(x)) &= (\partial_t - \partial_x^2)\left(u(x) + u'(x)\,(\cdot - x)_1 + \sum_{\tau\in\mathcal{T}_{\gamma-2}} \frac{\Upsilon_x[\tau]}{\tau!}\,\Pi_x\mathcal{I}(\tau)\right) \\ &= \sum_{\tau\in\mathcal{T}_{\gamma-2}} \frac{\Upsilon_x[\tau]}{\tau!}\,\Pi_x\,\tau = \sum_{\tau\in\mathcal{T}_{\gamma-\alpha}^+} \frac{\Upsilon_x[\tau\,\Xi]}{\tau!}\,\Pi_x(\tau\,\Xi) = \Pi_x(\hat{\sigma}(U_{\gamma-\alpha}(x))\,\Xi) \\ &= \Pi_x\,(\hat{\sigma}(U_\gamma(x))\,\Xi) - \sum_{\tau\in\mathcal{T}_{\gamma-\alpha,\gamma}^+} \frac{\Upsilon_x[\tau\,\Xi]}{\tau!}\,\Pi_x(\tau\,\Xi). \end{aligned}$$

Combined with the assumption of $u$ solving the PDE (3.1) and the identity (3.3), we conclude the first part. Applying the seminorm $\|\cdot\|_{\gamma-2}$ and the triangle inequality on it, we have that

$$\begin{aligned}\lambda^{\gamma}\|\mathscr{L}U_{\gamma}\|_{\gamma-2;B_{\lambda}(z)} &\lesssim \lambda^{\gamma}\|u\,|u|^{m-1}\|_{\gamma-2;B_{\lambda}(z)}+\lambda^{\gamma}\|\mathcal{R}(\hat{\sigma}(U_{\gamma})\,\Xi)-\Pi.(\hat{\sigma}(U_{\gamma})\,\Xi)\|_{\gamma-2;B_{\lambda}(z)}\\ &\quad+\lambda^{\gamma}\sum_{\tau\in\mathcal{T}^{+}_{\gamma-\alpha,\gamma}}\|\Upsilon.[\tau\,\Xi]\|_{B_{\lambda}(z)}\|\Pi.(\tau\,\Xi)\|_{\gamma-2;B_{\lambda}(z)},\end{aligned}$$

where the first term is interpreted as a constant germ. Using for $\gamma-2<\eta$ the elementary bound $\|\cdot\|_{\gamma-2;B_{\lambda}(z)}\leqslant\lambda^{-(\gamma-2-\eta)}\|\cdot\|_{\eta,B_{\lambda}(z)}$, we can control the low $(\gamma-2)$-regularity norm with the $L^{\infty}$-norm $(\eta=0)$ for the first two terms and a $|\tau|$-regularity for each of the last terms to conclude:

$$\begin{aligned}\lambda^{\gamma}\|\mathscr{L}U_{\gamma}\|_{\gamma-2;B_{\lambda}(z)} &\lesssim \lambda^{2}\|u\|^{m}_{B_{\lambda}(z)}+\lambda^{2}\|\mathcal{R}(\hat{\sigma}(U_{\gamma})\,\Xi)-\Pi.(\hat{\sigma}(U_{\gamma})\,\Xi)\|_{B_{\lambda}(z)}\\ &\quad+\lambda^{\gamma}\sum_{\tau\in\mathcal{T}^{+}_{\gamma-\alpha,\gamma}}\|\Upsilon.[\tau\,\Xi]\|_{B_{\lambda}(z)}\,\lambda^{-(\gamma-2-|\tau\Xi|)}\,\|\Pi.(\tau\,\Xi)\|_{|\tau|;B_{\lambda}(z)}.\end{aligned}$$

For the second term, since $\hat{\sigma}(U_{\gamma})\in\mathcal{D}^{\gamma}$, then $\hat{\sigma}(U_{\gamma})\,\Xi\in\mathcal{D}^{\gamma+\alpha-2}$ by [FH20, Theorem 14.5] and by the Reconstruction theorem [FH20, Theorem 13.6] (see [MW20b, Theorem 2.8] for the precise form used here) we conclude using the re-indexing coming from the identity (2.8):

$$\begin{aligned}\lambda^{2}\|\mathcal{R}(\hat{\sigma}(U_{\gamma})\,\Xi)-\Pi.(\hat{\sigma}(U_{\gamma})\,\Xi)\|_{B_{\lambda}(z)} &\lesssim \lambda^{2+\gamma+\alpha-2}\sum_{\rho\in\mathcal{T}_{\gamma+\alpha-2}}[\hat{\sigma}(U_{\gamma})\,\Xi;\rho]_{B_{\lambda}(z)}\,[\Pi;\rho] \qquad (3.28)\\ &= \lambda^{\gamma+\alpha}\sum_{\tau\in\mathcal{T}^{+}_{\gamma}}[\hat{\sigma}(U_{\gamma})\,\Xi;\tau\,\Xi]_{B_{\lambda}(z)}\,[\Pi;\tau\,\Xi]\\ &= \sum_{\tau\in\mathcal{T}^{+}_{\gamma}}\lambda^{\gamma-(|\tau|-|\mathfrak{n}(\tau)|)}[\hat{\sigma}(U_{\gamma});\tau]_{B_{\lambda}(z)}\lambda^{\alpha+|\tau|-|\mathfrak{n}(\tau)|}[\Pi;\tau\Xi]\\ &= \sum_{\tau\in\mathcal{T}^{+}_{\gamma}}\lambda^{\gamma-(|\tau|-|\mathfrak{n}(\tau)|)}\,[\hat{\sigma}(U_{\gamma});\tau]_{B_{\lambda}(z)}\,\lambda^{\alpha+|\tau|-|\mathfrak{n}(\tau)|}\,[\Pi;\tau\,\Xi].\end{aligned}$$

At last, for each $\tau\in\mathcal{T}^{+}_{\gamma-\alpha,\gamma}$, we have the identity

$$\lambda^{\gamma}\,\lambda^{-(\gamma-2-|\tau\Xi|)}\,\|\Pi.(\tau\,\Xi)\|_{|\tau|;B_{\lambda}(z)}=\lambda^{|\mathcal{I}(\tau\Xi)|}\,[\Pi;\tau\,\Xi]_{B_{\lambda}(z)},$$

which, combined with $\|\Upsilon.[\tau\,\Xi]\|\lesssim\|u'\|^{|\mathfrak{n}(\tau)|}$, consequence of Lemma 2.7, concludes the result. $\square$

**Proof of Lemma 3.6.** Fix $x\in B_{\lambda}(z)$ and $r\in(0,\lambda)$. At least one of $x_1+r$ or $x_1-r$ belongs to the interval $(x_1-\lambda,x_1+\lambda)$, so we define $y:=(x_0,x_1\pm r)\in\mathbb{R}\times\mathbb{R}$ such that $y\in B_{\lambda}(z)$ and $d(x,y)=|y_1-x_1|=r$. By the smallness conditions on $\lambda_t$ as in (3.25) from Remark 3.22 and (3.3)

$$\begin{aligned}r\,|u'(x)| &= |u'(x)\,(y-x)_1|\\ &\leqslant |U_{\beta}(x,y)-u'(x)\,(y-x)_1|+|U_{\beta}(x,y)|\\ &\lesssim d(x,y)^{\beta}\,[U_{\beta};\mathbf{1}]_{B_{\lambda}(z)}+|u(x)|+|u(y)|+\sum_{\tau\in\mathcal{T}_{\beta-2}}|\Upsilon_x[\tau]|\,[\Gamma\,\mathcal{I}(\tau);\mathbf{1}]\,d(x,y)^{|\mathcal{I}(\tau)|}\\ &\lesssim r^{\beta}\,[U_{\beta};\mathbf{1}]_{B_{\lambda}(z)}+\|u\|_{B_{\lambda}(z)}+\sum_{\tau\in\mathcal{T}_{\beta-2}}(r\,|u'(x)|)^{|\mathfrak{n}(\tau)|}\,[\Gamma\,\mathcal{I}(\tau);\mathbf{1}]\,\lambda^{|\mathcal{I}(\tau)|-|\mathfrak{n}(\tau)|}\\ &\lesssim r^{\beta}\,[U_{\beta};\mathbf{1}]_{B_{\lambda}(z)}+\|u\|_{B_{\lambda}(z)}+c\,r\,|u'(x)|+1,\end{aligned}$$

where the proportionality constant depends only on the set $\{\rho!\colon\rho\in\mathcal{T}^{+}_{\beta}\}$ and can be chosen uniformly for $\eta\in(1,\gamma]$ depending only on $\gamma$ and the regularity of the noise $\alpha-2=|\Xi|$. Taking supremum over $x\in B_{\lambda}(z)$ we obtain

$$r\,\|u'\|_{B_{\lambda}(z)}\lesssim r^{\eta}\,[U_{\eta};\mathbf{1}]_{B_{\lambda}(z)}+\|u\|_{B_{\lambda}(z)}+c\,r\,\|u'\|_{B_{\lambda}(z)}+1.$$

If $c=c(\gamma,\alpha)\in(0,1)$ from Definition 3.1 is small enough we can absorb the term $c\,r\,\|u'\|_{B_{\lambda}(z)}$ into the left hand side and obtain (3.4). For the second part define

$$r:=\left(\frac{\|u\|_{B_{\lambda}(z)}+1}{[U_{\gamma};\mathbf{1}]_{B_{\lambda}(z)}}\right)^{\frac{1}{\gamma}},$$

which is precisely chosen to balance both terms in (3.4). Then assuming the bound (3.5) means precisely that $r<\lambda$ which allows us to apply (3.4) for this particular value and obtain

$$\|u'\|_{B_\lambda(z)} \lesssim r^{-1}(\|u\|_{B_\lambda(z)}+1) = [U_\gamma;\mathbf{1}]_{B_\lambda(z)}^{\frac{1}{\gamma}}(\|u\|_{B_\lambda(z)}+1)^{1-\frac{1}{\gamma}}, \tag{3.29}$$

and therefore for any $\eta\in(1,\gamma]$

$$\begin{aligned}(\lambda\,\|u'\|_{B_\lambda(z)})^\eta &\lesssim (\lambda^\gamma[U_\gamma;\mathbf{1}]_{B_\lambda(z)})^{\frac{\eta}{\gamma}}(\|u\|_{B_\lambda(z)}+1)^{\eta(1-\gamma^{-1})}\\ &\lesssim (\lambda^\gamma[U_\gamma;\mathbf{1}]_{B_\lambda(z)}+1)\Big(\|u\|_{B_\lambda(z)}^{\eta(1-\gamma^{-1})}+1\Big)\\ &\lesssim \lambda^\gamma[U_\gamma;\mathbf{1}]_{B_\lambda(z)}+1+\|u\|_{B_\lambda(z)}^{\eta(1-\gamma^{-1})}+\lambda^\gamma[U_\gamma;\mathbf{1}]_{B_\lambda(z)}\|u\|_{B_\lambda(z)}^{\eta(1-\gamma^{-1})}\\ &\lesssim \lambda^\gamma[U_\gamma;\mathbf{1}]_{B_\lambda(z)}+1+\|u\|_{B_\lambda(z)}+\lambda^\gamma[U_\gamma;\mathbf{1}]_{B_\lambda(z)}\|u\|_{B_\lambda(z)}^{\eta(1-\gamma^{-1})},\end{aligned}$$

where we used that $\eta(1-\gamma^{-1})\leqslant 1$ since for $\gamma\in(2-\alpha,2)\subset(1,2)$ we have $1-\gamma^{-1}>0$ and $1<\frac{1}{\gamma-1}$, and therefore $\eta(1-\gamma^{-1})\leqslant 1\Longleftrightarrow\eta\leqslant(1-\gamma^{-1})^{-1}=\frac{\gamma}{\gamma-1}$ which holds for all $\eta\leqslant\gamma$ since $\gamma<\frac{\gamma}{\gamma-1}$. ☐

**Proof of Corollary 3.8.** We have by (3.3), (3.25) from Remark 3.22 and Lemma 3.6:

$$\begin{aligned}\|U_\gamma\|_{B_\lambda(z)} &\lesssim \|u\|_{B_\lambda(z)}+\sum_{\tau\in\mathcal{T}_{\gamma-2}}\|\Upsilon_\cdot[\tau]\|_{B_\lambda(z)}\lambda^{|\mathcal{I}(\tau)|}[\Gamma\,\mathcal{I}(\tau);\mathbf{1}]\\ &\lesssim \|u\|_{B_\lambda(z)}+\lambda^{|\mathcal{I}(\tau)|}[\Gamma\,\mathcal{I}(\tau);\mathbf{1}]+\lambda\,\|u'\|_{B_\lambda(z)}\lambda^{|\mathcal{I}(\tau)|-1}[\Gamma\,\mathcal{I}(\tau);\mathbf{1}]\\ &\lesssim \|u\|_{B_\lambda(z)}+c\,\lambda\,\|u'\|_{B_\lambda(z)}+1\\ &\lesssim c\,\lambda^\gamma[U_\gamma;\mathbf{1}]_{B_\lambda(z)}+\|u\|_{B_\lambda(z)}+1,\end{aligned}$$

which concludes the result. ☐

The following lemma allows us to control lower-order modelled seminorms.

**LEMMA 3.23**. *If $c=c(\alpha,\gamma)\in(0,1)$ in Definition 3.1 is small enough, then for any $\lambda>0$, $z\in\mathbb{R}^{1+1}$ and $\eta\in(0,\gamma]$*

$$\lambda^\eta[U_\eta;\mathbf{1}]_{B_\lambda(z)}\lesssim\lambda^\gamma[U_\gamma;\mathbf{1}]_{B_\lambda(z)}+\|u\|_{B_\lambda(z)}+1,$$

*and for any $\eta\in(1,\gamma]$*

$$\lambda^\eta[U_\beta;\boldsymbol{X}]_{B_\lambda(z)}\lesssim\lambda^\gamma[U_\gamma;\boldsymbol{X}]_{B_\lambda(z)}+\|u\|_{B_\lambda(z)}+1.$$

**Proof.** We have the identity

$$\begin{aligned}&U_\eta(x,y)-\mathbb{1}_{\eta>1}u'(x)(y-x)_1\\ =\;& U_\gamma(x,y)-u'(x)(y-x)_1+\mathbb{1}_{\eta\leqslant1}u'(x)(y-x)_1+\sum_{\tau\in\mathcal{T}_{\eta-2,\gamma-2}}\frac{\Upsilon_x[\tau]}{\tau!}\Pi_x(\mathcal{I}(\tau))(y),\end{aligned}$$

from where we obtain for all $x,y\in B_\lambda(z)$

$$\begin{aligned}&|U_\eta(x,y)-\mathbb{1}_{\eta>1}u'(x)(y-x)_1|\\ \lesssim\;& d(x,y)^\gamma[U_\gamma;\mathbf{1}]_{B_\lambda(z)}+\mathbb{1}_{\eta\leqslant1}|u'(x)|\,d(x,y)+\sum_{\tau\in\mathcal{T}_{\eta-2,\gamma-2}}|\Upsilon_x[\tau]|\,[\Gamma\,\mathcal{I}(\tau);\mathbf{1}]\,d(x,y)^{|\mathcal{I}(\tau)|}\\ \leqslant\;& d(x,y)^\eta\Bigg(\lambda^{\gamma-\eta}[U_\gamma;\mathbf{1}]_{B_\lambda(z)}+\mathbb{1}_{\eta\leqslant1}\lambda^{1-\eta}|u'(x)|+\sum_{\tau\in\mathcal{T}_{\eta-2,\gamma-2}}|u'(x)|^{|\mathfrak{n}(\tau)|}[\Gamma\mathcal{I}(\tau);\mathbf{1}]\lambda^{|\mathcal{I}(\tau)|-\eta}\Bigg),\end{aligned}$$

and using (3.25) from Remark 3.22 we obtain

$$\lambda^\eta[U_\eta;\mathbf{1}]_{B_\lambda(z)}\lesssim\lambda^\gamma[U_\gamma;\mathbf{1}]_{B_\lambda(z)}+\lambda\,\|u'\|_{B_\lambda(z)}+1,$$

and conclude the first part using (3.4) from Lemma 3.6. For the second part, the identity

$$\begin{aligned}\langle \boldsymbol{X}, U_\eta(y) - \Gamma_{yx} U_\eta(x)\rangle &= u'(y) - u'(x) - \sum_{\substack{\tau \in \mathcal{T}_{\eta-2}\\ |\mathcal{I}(\tau)|>1}} \frac{\Upsilon_x[\tau]}{\tau!} \gamma_{yx}(\mathcal{I}'(\tau))\\ &= \langle \boldsymbol{X}, U_\gamma(y) - \Gamma_{yx} U_\gamma(x)\rangle + \sum_{\tau \in \mathcal{T}_{\eta-2,\gamma-2}} \frac{\Upsilon_x[\tau]}{\tau!} \gamma_{yx}(\mathcal{I}'(\tau)),\end{aligned}$$

(observe that for $\tau \in \mathcal{T}_{\eta-2,\gamma-2}$ one has $|\mathcal{I}'(\tau)| = 1 + |\tau| > \eta - 1 > 0$) which implies the bound

$$\begin{aligned}&\langle \boldsymbol{X}, U_\eta(y) - \Gamma_{yx} U_\eta(x)\rangle\\ \lesssim\; & d(x,y)^{\gamma-1} [U_\gamma; \boldsymbol{X}]_{B_\lambda(z)} + \sum_{\tau \in \mathcal{T}_{\beta-2,\gamma-2}} \|u'\|^{|\mathfrak{n}(\tau)|} d(x,y)^{|\mathcal{I}(\tau\Xi)|-1} [\Gamma \mathcal{I}(\tau); \boldsymbol{X}]\\ \lesssim\; & d(x,y)^{\eta-1} \lambda^{\gamma-\eta} [U_\gamma; \boldsymbol{X}]_{B_\lambda(z)} + d(x,y)^{\eta-1} \sum_{\tau \in \mathcal{T}_{\beta-2,\gamma-2}} (\lambda \|u'\|)^{|\mathfrak{n}(\tau)|} \lambda^{|\mathcal{I}(\tau)|-|\mathfrak{n}(\tau)|-\eta} [\Gamma \mathcal{I}(\tau); \boldsymbol{X}]\end{aligned}$$

and using again (3.25) from Remark 3.22 we obtain

$$\lambda^\eta [U_\eta; \boldsymbol{X}]_{B_\lambda(z)} \lesssim \lambda^\gamma [U_\gamma; \boldsymbol{X}]_{B_\lambda(z)} + \lambda \|u'\|_{B_\lambda(z)} + 1,$$

and conclude by (3.4) from Lemma 3.6. □

**Proof of Lemma 3.9.** We consider first the case of $\tau$ with no polynomial decorations. Using the smallness of $\lambda_t$ from Definition 3.1 in the form (3.26) of Remark 3.22 (observe that trees in the boundary $\mathfrak{a}(\partial A_\tau)$ fall in the last category in there) we can rewrite the conclusion of Lemma 2.5 as

$$\lambda^{\gamma-|\tau|} [\hat{\sigma}(U_\gamma); \tau]_{B_\lambda(z)} \lesssim \lambda^{\gamma-|\tau|} [U_{\gamma-|\tau|}; \mathbf{1}]_{B_\lambda(z)} + \lambda \|u'\|_{B_\lambda(z)} + 1 + \mathbb{1}_{\gamma-|\tau|>1} (\lambda \|u'\|_{B_\lambda(z)})^{\gamma-|\tau|}, \tag{3.30}$$

where we used that the sublinear terms in $\lambda \|u'\|_{B_\lambda(z)}$ coming from the last term in Lemma 2.5 can be bounded either by the linear term $\lambda \|u'\|_{B_\lambda(z)}$ or the constant 1, and the corresponding superlinear terms, since $|\rho| - |\mathfrak{n}(\rho)| > 0$ for $\rho \in \mathfrak{a}(A_\tau)$, by either $(\lambda \|u'\|_{B_\lambda(z)})^{\gamma-|\tau|}$ or the constant 1. Lemma 3.23 immediately implies the result for $\gamma - |\tau| \leqslant 1$, while for the case $\gamma - |\tau| > 1$ we are only left with the last term, the superlinear in $u$ term, which by (3.6) in Lemma 3.6 is bounded as

$$(\lambda \|u'\|_{B_\lambda(z)})^{\gamma-|\tau|} \lesssim \lambda^\gamma [U_\gamma; \mathbf{1}]_{B_\lambda(z)} + \|u\|_{B_\lambda(z)} + 1 + \lambda^\gamma [U_\gamma; \mathbf{1}]_{B_\lambda(z)} \|u\|_{B_\lambda(z)}^{(\gamma-|\tau|)(1-\gamma^{-1})},$$

and which concludes the case without decorations.

Now we consider the case of $\tau$ with polynomial decorations. We use Remark 3.22 to post-process the corresponding bounds from Lemma 2.5

$$\begin{aligned}\lambda^{\gamma-(|\tau|-1)} [\hat{\sigma}(U_\gamma); \tau]_{B_\lambda(z)} \lesssim\; & \lambda \|u'\|_{B_\lambda(z)} (\lambda^{\gamma-|\tau|+1} [U_{\gamma-|\tau|+1}; \mathbf{1}]_{B_\lambda(z)})^{\frac{\gamma-|\tau|}{\gamma-|\tau|+1}}\\ & + (\lambda \|u'\|_{B_\lambda(z)})^{\gamma-|\tau|+1} + \lambda \|u'\|_{B_\lambda(z)} + 1\\ & + \max_{\rho \in \mathfrak{a}(A_\tau)} \lambda^{\gamma-|\tau|-|\rho|+1} [U_{\gamma-|\tau|-|\rho|+1}; \boldsymbol{X}]_{B_\lambda(z)}.\end{aligned} \tag{3.31}$$

For the second term in (3.31) we use (3.6) from Lemma 3.6 to conclude

$$(\lambda \|u'\|_{B_\lambda(z)})^{\gamma-|\tau|+1} \lesssim \lambda^\gamma [U_\gamma; \mathbf{1}]_{B_\lambda(z)} + \|u\|_{B_\lambda(z)} + 1 + \lambda^\gamma [U_\gamma; \mathbf{1}]_{B_\lambda(z)} \|u\|_{B_\lambda(z)}^{(\gamma-|\tau|+1)(1-\gamma^{-1})}.$$

For the last terms in (3.31) we use Lemma 3.23 and [EW24, Lemma 3.9] to obtain

$$\begin{aligned}\lambda^{\gamma-|\tau|-|\rho|+1} [U_{\gamma-|\tau|-|\rho|+1}; \boldsymbol{X}]_{B_\lambda(z)} &\lesssim \lambda^\gamma [U_\gamma; \boldsymbol{X}]_{B_\lambda(z)} + \|u\|_{B_\lambda(z)} + 1\\ &\lesssim \lambda^\gamma [U_\gamma; \mathbf{1}]_{B_\lambda(z)} + \lambda^\gamma [U_\gamma]_{\gamma\text{-3pt}; B_\lambda(z)} + \|u\|_{B_\lambda(z)} + 1.\end{aligned}$$

At last, for the first term in (3.31) we obtain (3.29) and by Lemma 3.23:

$$
\begin{aligned}
& \lambda \|u'\|_{B_\lambda(z)} (\lambda^{\gamma-|\tau|+1} [U_{\gamma-|\tau|+1}; \mathbf{1}]_{B_\lambda(z)})^{\frac{\gamma-|\tau|}{\gamma-|\tau|+1}} \\
\lesssim\ & \lambda [U_\gamma; \mathbf{1}]_{B_\lambda(z)}^{\frac{1}{\gamma}} (\|u\|_{B_\lambda(z)} + 1)^{1-\frac{1}{\gamma}} (\lambda^{\gamma-|\tau|+1} [U_{\gamma-|\tau|+1}; \mathbf{1}]_{B_\lambda(z)})^{\frac{\gamma-|\tau|}{\gamma-|\tau|+1}} \\
\lesssim\ & (\lambda^\gamma [U_\gamma; \mathbf{1}]_{B_\lambda(z)})^{\frac{1}{\gamma}} \left(\|u\|_{B_\lambda(z)}^{1-\gamma^{-1}} + 1\right) (\lambda^\gamma [U_\gamma; \mathbf{1}]_{B_\lambda(z)} + \|u\|_{B_\lambda(z)} + 1)^{\frac{\gamma-|\tau|}{\gamma-|\tau|+1}} \\
\lesssim\ & (\lambda^\gamma [U_\gamma; \mathbf{1}]_{B_\lambda(z)})^{\gamma^{-1}} \left(\|u\|_{B_\lambda(z)}^{1-\gamma^{-1}} + 1\right) \left( (\lambda^\gamma [U_\gamma; \mathbf{1}]_{B_\lambda(z)})^{\frac{\gamma-|\tau|}{\gamma-|\tau|+1}} + \|u\|_{B_\lambda(z)}^{\frac{\gamma-|\tau|}{\gamma-|\tau|+1}} + 1 \right) \\
\lesssim\ & (\lambda^\gamma [U_\gamma; \mathbf{1}]_{B_\lambda(z)})^{\frac{1}{\gamma}+\frac{\gamma-|\tau|}{\gamma-|\tau|+1}} \left(\|u\|_{B_\lambda(z)}^{1-\gamma^{-1}} + 1\right) + (\lambda^\gamma [U_\gamma; \mathbf{1}]_{B_\lambda(z)})^{\gamma^{-1}} \left( \|u\|_{B_\lambda(z)}^{1-\frac{1}{\gamma}+\frac{\gamma-|\tau|}{\gamma-|\tau|+1}} + 1 \right). \qquad (3.32)
\end{aligned}
$$

$$
\begin{aligned}
\lesssim\ & \lambda^\gamma [U_\gamma; \mathbf{1}]_{B_\lambda(z)} \|u\|_{B_\lambda(z)}^{1-\gamma^{-1}} + \lambda^\gamma [U_\gamma; \mathbf{1}]_{B_\lambda(z)} + 1 + (\lambda^\gamma [U_\gamma; \mathbf{1}]_{B_\lambda(z)})^{\gamma^{-1}} \|u\|_{B_\lambda(z)}^{1-\frac{1}{\gamma}+\frac{\gamma-|\tau|}{\gamma-|\tau|+1}} \\
\lesssim\ & \lambda^\gamma [U_\gamma; \mathbf{1}]_{B_\lambda(z)} \|u\|_{B_\lambda(z)}^{(\gamma-|\tau|+1)(1-\gamma^{-1})} + \lambda^\gamma [U_\gamma; \mathbf{1}]_{B_\lambda(z)} + 1 + (\lambda^\gamma [U_\gamma; \mathbf{1}]_{B_\lambda(z)})^{\gamma^{-1}} \|u\|_{B_\lambda(z)}^{1-\frac{1}{\gamma}+\frac{\gamma-|\tau|}{\gamma-|\tau|+1}}.
\end{aligned}
$$

Since $|\tau| \geqslant 1$ then $1 \leqslant \gamma - |\tau| + 1 \leqslant \gamma < 2$ and therefore

$$
\frac{1}{\gamma} + \frac{\gamma-|\tau|}{\gamma-|\tau|+1} \leqslant \frac{1}{\gamma-|\tau|+1} + \frac{\gamma-|\tau|}{\gamma-|\tau|+1} = 1,
$$

which makes the exponent in the modelled norm of first term in (3.32) sublinear and we obtain, and since $1 - \gamma^{-1} \leqslant (\gamma - |\tau| + 1)(1 - \gamma^{-1})$ then

$$
\begin{aligned}
(\lambda^\gamma [U_\gamma; \mathbf{1}]_{B_\lambda(z)})^{\frac{1}{\gamma}+\frac{\gamma-|\tau|}{\gamma-|\tau|+1}} \left(\|u\|_{B_\lambda(z)}^{1-\gamma^{-1}} + 1\right) &\lesssim \lambda^\gamma [U_\gamma; \mathbf{1}]_{B_\lambda(z)} \|u\|_{B_\lambda(z)}^{1-\gamma^{-1}} + \lambda^\gamma [U_\gamma; \mathbf{1}]_{B_\lambda(z)} + 1 \\
&\lesssim \lambda^\gamma [U_\gamma; \mathbf{1}]_{B_\lambda(z)} \|u\|_{B_\lambda(z)}^{(\gamma-|\tau|+1)(1-\gamma^{-1})} + \lambda^\gamma [U_\gamma; \mathbf{1}]_{B_\lambda(z)} + 1.
\end{aligned}
$$

Let $a = \lambda^\gamma [U_\gamma; \mathbf{1}]_{B_\lambda(z)}$, $b = \|u\|_{B_\lambda(z)}$ and $\theta = 1 - \frac{1}{\gamma} + \frac{\gamma-|\tau|}{\gamma-|\tau|+1} - (\gamma - |\tau| + 1) \frac{(1-\gamma^{-1})}{\gamma}$, then using Young's inequality the last term in (3.32) is bounded as

$$
\begin{aligned}
a^{\gamma^{-1}} b^{1-\frac{1}{\gamma}+\frac{\gamma-|\tau|}{\gamma-|\tau|+1}} &= (a\, b^{(\gamma-|\tau|+1)(1-\gamma^{-1})})^{\gamma^{-1}} b^\theta \\
&\lesssim a\, b^{(\gamma-|\tau|+1)(1-\gamma^{-1})} + b^{\theta \frac{\gamma}{\gamma-1}} \\
&\lesssim a\, b^{(\gamma-|\tau|+1)(1-\gamma^{-1})} + b + 1,
\end{aligned}
$$

and therefore the term is bound as stated in the result. To justify the last line we used that

$$
\begin{aligned}
\theta \frac{\gamma}{\gamma-1} &= \left(1 - \frac{1}{\gamma} + \frac{\gamma-|\tau|}{\gamma-|\tau|+1} - (\gamma-|\tau|+1) \frac{1-\gamma^{-1}}{\gamma}\right) \frac{\gamma}{\gamma-1} \\
&= 1 + \frac{\gamma(\gamma-|\tau|)}{(\gamma-|\tau|+1)(\gamma-1)} - \frac{\gamma-|\tau|+1}{\gamma},
\end{aligned}
$$

that $1 < \gamma - |\tau| + 1 < \gamma < 2$, and that the map $x \mapsto \frac{x^2}{x-1}$ is decreasing in the interval, and therefore

$$
\frac{\gamma^2}{\gamma-1} \leqslant \frac{(\gamma-|\tau|+1)^2}{\gamma-|\tau|} \Longrightarrow \frac{\gamma(\gamma-|\tau|)}{(\gamma-|\tau|+1)(\gamma-1)} \leqslant \frac{\gamma-|\tau|+1}{\gamma}
$$

which makes the exponent $\theta \frac{\gamma}{\gamma-1} \leqslant 1$ sublinear. □

**Proof of Lemma 3.10.** By Definition 3.1 $\lambda_t \leqslant \|u\|_{P_t}^{-\frac{(m-1)}{2}}$, and therefore the first term from Lemma 3.5 can be bounded for any $0 < \lambda \leqslant \lambda_t$ as

$$
\lambda^2 \|u\|_{B_\lambda(z)}^m \lesssim \lambda_t^2 \|u\|_{P_t}^m \leqslant \|u\|_{P_t}.
$$

For the terms coming from the reconstruction in Lemma 3.5 we use Lemma 3.9. First consider $\tau \in \mathcal{T}_\gamma^+$ with $\mathfrak{n}(\tau) = 0$. If $\gamma - |\tau| \leqslant 1$ then by (3.25) from Remark 3.22 ($2 + |\tau\,\Xi| = \alpha + |\tau|$)

$$\lambda^{\gamma - |\tau|}\,[\hat{\sigma}(U_\gamma); \tau]_{B_\lambda(z)}\,\lambda^{\alpha + |\tau|}\,[\Pi; \tau\,\Xi] \lesssim c\,\lambda^\gamma\,[U_\gamma; \mathbf{1}]_{B_\lambda(z)} + \|u\|_{B_\lambda(z)} + 1. \tag{3.33}$$

If $\gamma - |\tau| > 1$, the only difference to before is the last term coming from (3.7) for which we use that, by definition (2.14), $[\Pi; \tau\,\Xi] \leqslant |||(\Pi, \Gamma)|||^{\mathfrak{c}(\tau\Xi)}$ and then by (3.27) from Remark 3.22:

$$\|u\|_{P_t}^{(\gamma - |\tau|)(1 - \gamma^{-1})}\,\lambda_t^{\alpha + |\tau|}\,[\Pi; \tau\,\Xi] \leqslant c. \tag{3.34}$$

We conclude that

$$\lambda^\gamma\,[U_\gamma; \mathbf{1}]_{B_\lambda(z)} \|u\|_{P_t}^{(\gamma - |\tau|)(1 - \gamma^{-1})}\,\lambda_t^{\alpha + |\tau|}\,[\Pi; \tau\,\Xi] \leqslant c\,\lambda^\gamma\,[U_\gamma; \mathbf{1}]_{B_\lambda(z)} \tag{3.35}$$

and (3.33) also holds if $\gamma - |\tau| > 1$. If $\mathfrak{n}(\tau) \neq 0$, then by (3.8) and the smallness of $\lambda_t$ in the form of (3.25) from Remark 3.22:

$$\begin{aligned} & \lambda^{\gamma - (|\tau| - 1)}\,[\hat{\sigma}(U_\gamma); \tau]_{B_\lambda(z)}\,\lambda^{\alpha + |\tau| - 1}\,[\Pi; \tau\,\Xi] \\ \lesssim\ & c\,\lambda^\gamma\,[U_\gamma; \mathbf{1}]_{B_\lambda(z)} + \|u\|_{B_\lambda(z)} + 1 + \lambda^\gamma\,[U_\gamma]_{\gamma\text{-3pt}; B_\lambda(z)} \\ & + (\lambda^\gamma\,[U_\gamma; \mathbf{1}]_{B_\lambda(z)} + \|u\|_{B_\lambda(z)})\,\|u\|_P^{(\gamma - |\tau| + 1)(1 - \gamma^{-1})}\,\lambda_t^{\alpha + |\tau| - 1}\,[\Pi; \tau\,\Xi]_{B_\lambda(z)}. \end{aligned} \tag{3.36}$$

In addition, by (3.27) from Remark 3.22:

$$\|u\|_{P_t}^{(\gamma - |\tau| + 1)(1 - \gamma^{-1})}\,\lambda_t^{\alpha + |\tau| - 1}\,[\Pi; \tau\,\Xi] \leqslant c.$$

and we conclude that (3.36) is bounded by (3.9). For the last terms from Lemma 3.5 we have for any $\tau \in \mathcal{T}_{\gamma - \alpha, \gamma}^+$ that the bound

$$(\lambda\,\|u'\|_{B_\lambda(z)})^{|\mathfrak{n}(\tau)|}\,\lambda^{|\mathcal{I}(\tau\Xi)| - |\mathfrak{n}(\tau)|}\,[\Pi; \tau\,\Xi]_{B_\lambda(z)} \leqslant c\,(\lambda\,\|u'\|_{B_\lambda(z)})^{|\mathfrak{n}(\tau)|} \leqslant c\,\lambda^\gamma\,[U_\gamma; \mathbf{1}]_{B_\lambda(z)} + \|u\|_{B_\lambda(z)} + 1$$

follows from (3.25) from Remark 3.22 in combination with (3.4) from Lemma 3.6 when $\mathfrak{n}(\tau) \neq 0$. □

**Proof of Lemma 3.11.** We prove first (3.12) since it is the most complex one. Our starting point is the bound (3.10). Fix $\tau \in \mathcal{T}_{\gamma - \alpha}^+$ without polynomial decorations. An explicit computation leads to the identity

$$(\gamma - |\tau|)\,(1 - \gamma^{-1}) - (\gamma - \alpha - |\tau|)\,(1 - (\gamma - \alpha)^{-1}) = \alpha \left(1 - \frac{|\tau|}{\gamma\,(\gamma - \alpha)}\right). \tag{3.37}$$

Since $|\tau| \in [0, \gamma - \alpha)$ then $1 - \frac{|\tau|}{\gamma\,(\gamma - \alpha)} \geqslant 1 - \gamma^{-1} > 0$, which combined with the conditions $\alpha \in (0, 1)$ and $\gamma \in (1, 2)$ allows us to conclude that (3.37) is positive, and therefore

$$(\gamma - \alpha - |\tau|)\,(1 - (\gamma - \alpha)^{-1}) < (\gamma - |\tau|)\,(1 - \gamma^{-1}) \leqslant \gamma\,(1 - \gamma^{-1}) = \gamma - 1 < 1. \tag{3.38}$$

By Lemma 3.9 (with $\gamma - \alpha$ in the role of $\gamma$ in there) and Lemma 3.23 (to go back from $\gamma - \alpha$ to $\gamma$ seminorms) we have if $\mathfrak{n}(\tau) = 0$:

$$\begin{aligned} & \lambda^{\gamma - \alpha - |\tau|}\,[\hat{\sigma}(U_{\gamma - \alpha}); \tau]_{B_\lambda(z)} \\ \lesssim\ & \lambda^{\gamma - \alpha}\,[U_{\gamma - \alpha}; \mathbf{1}]_{B_\lambda(z)} + \|u\|_{B_\lambda(z)} + 1 + \mathbb{1}_{\gamma - \alpha - |\tau| > 1}\,\lambda^{\gamma - \alpha}\,[U_{\gamma - \alpha}; \mathbf{1}]_{B_\lambda(z)}\,\|u\|_{B_\lambda(z)}^{(\gamma - \alpha - |\tau|)(1 - (\gamma - \alpha)^{-1})} \\ \lesssim\ & \lambda^\gamma\,[U_\gamma; \mathbf{1}]_{B_\lambda(z)} + \|u\|_{B_\lambda(z)} + 1 \\ & + \mathbb{1}_{\gamma - \alpha - |\tau| > 1}\,(\lambda^\gamma\,[U_\gamma; \mathbf{1}]_{B_\lambda(z)} + \|u\|_{B_\lambda(z)} + 1)\,\|u\|_{B_\lambda(z)}^{(\gamma - \alpha - |\tau|)(1 - (\gamma - \alpha)^{-1})} \\ \lesssim\ & \lambda^\gamma\,[U_\gamma; \mathbf{1}]_{B_\lambda(z)} + \|u\|_{B_\lambda(z)} + 1 + \mathbb{1}_{\gamma - \alpha - |\tau| > 1}\,(\lambda^\gamma\,[U_\gamma; \mathbf{1}]_{B_\lambda(z)} + \|u\|_{B_\lambda(z)})\,\|u\|_{B_\lambda(z)}^{(\gamma - \alpha - |\tau|)(1 - (\gamma - \alpha)^{-1})}, \end{aligned}$$

where we used Lemma 3.23, and (3.38) to remove the sublinear in $\|u\|$ terms. If $\gamma-\alpha-|\tau|\leqslant 1$ then

$$\lambda^{\gamma-\alpha-|\tau|}[\hat{\sigma}(U_{\gamma-\alpha});\tau]_{B_\lambda(z)}\lambda^{\alpha+|\tau|-|\mathfrak{n}(\tau)|}[\Gamma\mathcal{I}(\tau\Xi);\mathbf{1}]\lesssim\lambda^\gamma[U_\gamma;\mathbf{1}]_{B_\lambda(z)}+\|u\|_{B_\lambda(z)}+1, \tag{3.39}$$

which follows from (3.25) in Remark 3.22 ($\alpha+|\tau|=|\mathcal{I}(\tau\Xi)|$). If $\gamma-\alpha-|\tau|>1$, then by (3.38)

$$\begin{aligned}
&(\lambda^\gamma[U_\gamma;\mathbf{1}]_{B_\lambda(z)}+\|u\|_{B_\lambda(z)})\|u\|_{B_\lambda(z)}^{(\gamma-\alpha-|\tau|)(1-(\gamma-\alpha)^{-1})}\lambda^{\alpha+|\tau|-|\mathfrak{n}(\tau)|}[\Gamma\mathcal{I}(\tau\Xi);\mathbf{1}]\\
\leqslant\;&(\lambda^\gamma[U_\gamma;\mathbf{1}]_{B_\lambda(z)}+\|u\|_{B_\lambda(z)})\|u\|_{B_\lambda(z)}^{(\gamma-|\tau|)(1-\gamma^{-1})}\lambda^{\alpha+|\tau|-|\mathfrak{n}(\tau)|}[\Gamma\mathcal{I}(\tau\Xi);\mathbf{1}]\\
\leqslant\;&c\,\lambda^\gamma[U_\gamma;\mathbf{1}]_{B_\lambda(z)}+\|u\|_{B_\lambda(z)},
\end{aligned}$$

where the second bound follows as in (3.34), and therefore (3.39) also holds in this case. For the case with decorations, $\mathfrak{n}(\tau)\neq 0$, we obtain, analogously to the proof of Lemma 3.9, that

$$\begin{aligned}
&\lambda^{\gamma-\alpha-(|\tau|-1)}[\hat{\sigma}(U_{\gamma-\alpha});\tau]_{B_\lambda(z)}\lambda^{\alpha+|\tau|-|\mathfrak{n}(\tau)|}[\Gamma\mathcal{I}(\tau\Xi);\mathbf{1}]\\
\lesssim\;&c\,\lambda^\gamma[U_\gamma;\mathbf{1}]_{B_\lambda(z)}+\|u\|_{B_\lambda(z)}+1+c\,\lambda^\gamma[U_\gamma]_{\gamma\text{-3pt};B_\lambda(z)}\\
&+(\lambda^\gamma[U_\gamma;\mathbf{1}]_{B_\lambda(z)}+\|u\|_{B_\lambda(z)})\|u\|_{B_\lambda(z)}^{(\gamma-\alpha-|\tau|+1)(1-(\gamma-\alpha)^{-1})}\lambda^{\alpha+|\tau|-|\mathfrak{n}(\tau)|}[\Gamma\mathcal{I}(\tau\Xi);\mathbf{1}]\\
\lesssim\;&c\,\lambda^\gamma[U_\gamma;\mathbf{1}]_{B_\lambda(z)}+\|u\|_{B_\lambda(z)}+1+c\,\lambda^\gamma[U_\gamma]_{\gamma\text{-3pt};B_\lambda(z)}
\end{aligned} \tag{3.40}$$

where we used that, analogously to the first part

$$(\gamma-\alpha-|\tau|+1)(1-(\gamma-\alpha)^{-1})<(\gamma-|\tau|+1)(1-\gamma^{-1}),$$

and therefore by (3.27) from Remark 3.22

$$\|u\|_{B_\lambda(z)}^{(\gamma-\alpha-|\tau|+1)(1-(\gamma-\alpha)^{-1})}\lambda^{\alpha+|\tau|-|\mathfrak{n}(\tau)|}[\Gamma\mathcal{I}(\tau\Xi);\mathbf{1}]\leqslant c.$$

Combining the bounds (3.39) and (3.40) in (3.10) we conclude that

$$\lambda^\gamma[U_\gamma]_{\gamma\text{-3pt};B_\lambda(z)}\lesssim c\,\lambda^\gamma[U_\gamma;\mathbf{1}]_{B_\lambda(z)}+\|u\|_{B_\lambda(z)}+1+c\,\lambda^\gamma[U_\gamma]_{\gamma\text{-3pt};B_\lambda(z)}.$$

By imposing another smallness condition on $c\in(0,1)$ depending on the implicit constant (which only depends on $\gamma$ and $\alpha$) we can absorb this term onto the right hand-side and conclude this part of the result.

We proceed to prove (3.11). Analogously to (3.10), by [EW24, Lemma 4.10]:

$$\begin{aligned}
&\lambda^\gamma[\Lambda_\gamma]_{(\gamma-1)\text{-3pt};B_\lambda(z)}\\
\lesssim\;&\lambda^\gamma\sum_{\tau\in\mathcal{T}_{-1,\gamma-2}}[U_\gamma;\mathcal{I}(\tau)]_{B_\lambda(z)}[\Gamma\mathcal{I}(\tau);\boldsymbol{X}]_{B_\lambda(z)}\\
=\;&\sum_{\tau\in\mathcal{T}^+_{1-\alpha,\gamma-\alpha}}\lambda^{\gamma-|\mathcal{I}(\tau\Xi)|+|\mathfrak{n}(\mathcal{I}(\tau\Xi))|}[U_\gamma;\mathcal{I}(\tau\Xi)]_{B_\lambda(z)}\lambda^{|\mathcal{I}(\tau\Xi)|-|\mathfrak{n}(\mathcal{I}(\tau\Xi))|}[\Gamma\mathcal{I}(\tau\Xi);\boldsymbol{X}]_{B_\lambda(z)}\\
=\;&\sum_{\tau\in\mathcal{T}^+_{1-\alpha,\gamma-\alpha}}\lambda^{\gamma-\alpha-(|\tau|-|\mathfrak{n}(\tau)|)}[\hat{\sigma}(U_{\gamma-\alpha});\tau]_{B_\lambda(z)}\lambda^{\alpha+|\tau|-|\mathfrak{n}(\tau)|}[\Gamma\mathcal{I}(\tau\Xi);\boldsymbol{X}]_{B_\lambda(z)}.
\end{aligned}$$

From here the result follows analogously as for the 3-point seminorm of $U_\gamma$ with the only difference being the term $[\Gamma\mathcal{I}(\tau);\boldsymbol{X}]$ replacing the term $[\Gamma\mathcal{I}(\tau);\mathbf{1}]$. ☐

**Proof of Corollary 3.3.** By Lemma 3.6

$$\lambda\|u'\|_{B_\lambda(z)}\lesssim\lambda^\gamma[U_\gamma;\mathbf{1}]_{B_\lambda(z)}+\|u\|_{B_\lambda(z)}+1\lesssim\|u\|_{P_t}.$$

By Lemma 3.11 we have

$$\lambda^{\gamma}\,[U_{\gamma}]_{\gamma\text{-3pt};B_{\lambda}(z)} \lesssim c\,\lambda^{\gamma}\,[U_{\gamma};\mathbf{1}]_{B_{\lambda}(z)} + \|u\|_{B_{\lambda}(z)} + 1 \lesssim \|u\|_{P_t}.$$

By Lemma 3.23 we have

$$\lambda^{\eta}\,[U_{\eta};\mathbf{1}]_{B_{\lambda}(z)} \lesssim \lambda^{\gamma}\,[U_{\gamma};\mathbf{1}]_{B_{\lambda}(z)} + \|u\|_{B_{\lambda}(z)} + 1 \lesssim \|u\|_{P_t},$$

and

$$\lambda^{\eta}\,[U_{\eta};\boldsymbol{X}]_{B_{\lambda}(z)} \lesssim \lambda^{\gamma}\,[U_{\gamma};\boldsymbol{X}]_{B_{\lambda}(z)} + \|u\|_{B_{\lambda}(z)} + 1 \lesssim \lambda^{\gamma}\,[U_{\gamma};\mathbf{1}]_{B_{\lambda}(z)} + \lambda^{\gamma}\,[U_{\gamma}]_{\gamma\text{-3pt};B_{\lambda}(z)} + \|u\|_{P_t} \lesssim \|u\|_{P_t},$$

where we used Lemma 3.11 in the form

$$\lambda^{\gamma}\,[U_{\gamma};\boldsymbol{X}]_{B_{\lambda}(z)} \lesssim \lambda^{\gamma}\,[U_{\gamma};\mathbf{1}]_{B_{\lambda}(z)} + \lambda^{\gamma}\,[U_{\gamma}]_{\gamma\text{-3pt};B_{\lambda}(z)}.$$

By Lemma 3.9 we have for $\mathfrak{n}(\tau)=0$

$$\begin{aligned}\lambda^{\gamma-|\tau|}\,[\hat{\sigma}(U_{\gamma});\tau]_{B_{\lambda}(z)} &\lesssim \lambda^{\gamma}\,[U_{\gamma};\mathbf{1}]_{B_{\lambda}(z)} + \|u\|_{B_{\lambda}(z)} + 1 + \mathbb{1}_{\gamma-|\tau|>1}\,\lambda^{\gamma}\,[U_{\gamma};\mathbf{1}]_{B_{\lambda}(z)}\,\|u\|_{B_{\lambda}(z)}^{(\gamma-|\tau|)(1-\gamma^{-1})},\\ &\lesssim \|u\|_{P_t} + \mathbb{1}_{\gamma-|\tau|>1}\,\|u\|_{P_t}^{1+(\gamma-|\tau|)(1-\gamma^{-1})}\end{aligned}$$

and for $\mathfrak{n}(\tau)\neq 0$

$$\begin{aligned}\lambda^{\gamma-(|\tau|-1)}\,[\hat{\sigma}(U_{\gamma});\tau]_{B_{\lambda}(z)} &\lesssim \lambda^{\gamma}\,[U_{\gamma};\mathbf{1}]_{B_{\lambda}(z)} + \|u\|_{B_{\lambda}(z)} + 1 + \lambda^{\gamma}\,[U_{\gamma}]_{\gamma\text{-3pt};B_{\lambda}(z)}\\ &\quad + (\lambda^{\gamma}\,[U_{\gamma};\mathbf{1}]_{B_{\lambda}(z)} + \|u\|_{B_{\lambda}(z)})\,\|u\|_{B_{\lambda}(z)}^{(\gamma-|\tau|+1)(1-\gamma^{-1})}\\ &\lesssim \|u\|_{B_{\lambda}(z)}^{1+(\gamma-|\tau|+1)(1-\gamma^{-1})}.\end{aligned}$$

□

## 3.7. Proofs of Lemmas from Section 3.2

The following result gives us an identity between the generalised derivative and the real derivative of the regularisation.

**Lemma 3.24**. *For any $\gamma\in(1,2)$ we have the identity*

$$u'(x) = (\partial_{e_1}u_{\lambda})(x) + \langle U_{\gamma}(x,\cdot) - u'(x)\,(\cdot - x)_1, \partial_{e_1}\varphi_x^{\lambda}\rangle - \sum_{\tau\in\mathcal{T}_{\gamma-2}} \frac{\Upsilon_x[\tau]}{\tau!}\,\langle \Pi_x\,\mathcal{I}'(\tau), \varphi_x^{\lambda}\rangle. \tag{3.41}$$

**Proof.** Recall the definition

$$u_{\lambda}(x) := \int \varphi^{\lambda}(x-y)\,u(y)\,\mathrm{d}y = \int \tilde{\varphi}^{\lambda}(y-x)\,u(y)\,\mathrm{d}y = \langle u, \tilde{\varphi}_x^{\lambda}\rangle,$$

where $\tilde{\varphi}(x_0,x_1) = \varphi(-x_0,-x_1) = \varphi(-x_0,x_1)$ since $\varphi$ is symmetric in space, and observe that

$$\partial_{e_1}u_{\lambda} = \langle u, \partial_{x_1}\tilde{\varphi}_x^{\lambda}\rangle = \langle u, \partial_{x_1}(\lambda^{-|\mathfrak{s}|}\,\tilde{\varphi}(\lambda.(\cdot - x)))\rangle = -\langle u, \partial_{e_1}(\lambda^{-|\mathfrak{s}|}\,\tilde{\varphi}(\lambda.(\cdot - x)))\rangle = -\langle u, \partial_{e_1}\tilde{\varphi}_x^{\lambda}\rangle,$$

where $\partial_{x_1}\tilde{\varphi}_x^{\lambda}$ indicates the partial derivative in the base-point $x$ and $\partial_{e_1}\tilde{\varphi}_x^{\lambda}$ the partial derivative in the implicit argument and $\lambda.x = (\lambda^2 x_0, \lambda x_1)$. On the other hand, by (3.3) we have

$$U_{\gamma}(x,y) - u'(x)\,(y-x)_1 \;=\; u(y) - u(x) - u'(x)\,(y-x)_1 - \sum_{\tau\in\mathcal{T}_{\gamma-2}} \frac{\Upsilon_x[\tau]}{\tau!}\,\Pi_x(\mathcal{I}(\tau))(y).$$

and therefore

$$
\begin{aligned}
(\partial_{e_1} u_\lambda)(x) &= -\langle U_\gamma(x,\cdot) - u'(x)\,(\cdot - x)_1, \partial_{e_1}\tilde{\varphi}_x^\lambda\rangle \\
&\quad -u(x)\,\langle 1, \partial_{e_1}\tilde{\varphi}_x^\lambda\rangle - u'(x)\,\langle(\cdot - x)_1, \partial_{e_1}\tilde{\varphi}_x^\lambda\rangle - \sum_{\tau\in\mathcal{T}_{\gamma-2}} \frac{\Upsilon_x[\tau]}{\tau!}\,\langle \Pi_x\,\mathcal{I}(\tau), \partial_{e_1}\tilde{\varphi}_x^\lambda\rangle \\
&= -\langle U_\gamma(x,\cdot) - u'(x)\,(\cdot - x)_1, \partial_{e_1}\tilde{\varphi}_x^\lambda\rangle \\
&\quad +u(x)\,\langle \partial_{e_1} 1, \tilde{\varphi}_x^\lambda\rangle + u'(x)\,\langle \partial_{e_1}(\cdot - x)_1, \tilde{\varphi}_x^\lambda\rangle + \sum_{\tau\in\mathcal{T}_{\gamma-2}} \frac{\Upsilon_x[\tau]}{\tau!}\,\langle \partial_{e_1}\Pi_x\,\mathcal{I}(\tau), \tilde{\varphi}_x^\lambda\rangle \\
&= -\langle U_\gamma(x,\cdot) - u'(x)\,(\cdot - x)_1, \partial_{e_1}\tilde{\varphi}_x^\lambda\rangle + u'(x) + \sum_{\tau\in\mathcal{T}_{\gamma-2}} \frac{\Upsilon_x[\tau]}{\tau!}\,\langle \Pi_x\,\mathcal{I}'(\tau), \tilde{\varphi}_x^\lambda\rangle
\end{aligned}
$$

which concludes the identity. □

**Proof of Lemma 3.14.** By assumption

$$\lambda_t = c^{\frac{1}{\alpha}}\,(1 + |||(\Pi,\Gamma)|||)^{-\frac{1}{\alpha}}\,\|u\|_{P_t}^{-\frac{1}{\alpha}(\gamma-1)},$$

and therefore

$$\lambda_t^{-1}\,\|u\|_{P_t} \sim (1 + |||(\Pi,\Gamma)|||)^{\frac{1}{\alpha}}\,\|u\|_{P_t}^{\frac{1}{\alpha}(\gamma-1+\alpha)}. \tag{3.42}$$

Since $z \in P_{t+\lambda_t}$, Corollary 3.3, applied at scale $\lambda_t$, gives

$$|u'(z)| \leqslant \|u'\|_{B_{\lambda_t}(z)} \lesssim \lambda_t^{-1}\|u\|_{P_t},$$

and, since $\lambda \leqslant \lambda_t$,

$$[U_\gamma;1]_{B_\lambda(z)} \leqslant [U_\gamma;1]_{B_{\lambda_t}(z)} \lesssim \lambda_t^{-\gamma}\,\|u\|_{P_t}.$$

Moreover, since every model symbol appearing below contains at least one noise, we will repeatedly use

$$\|(\Pi,\Gamma)\|^q \leqslant \|(\Pi,\Gamma)\|\,(1 + \|(\Pi,\Gamma)\|)^{q-1}, \qquad q \geq 1.$$

We proceed by identifying the function $b_\lambda^\Pi$. By (2.17) and Lemma 2.7 we have

$$\langle \Pi_z(\hat{\sigma}(U_\gamma)(z)\,\Xi), \varphi_z^\lambda\rangle = \sum_{\tau\in\mathcal{T}_\gamma^+} \frac{\Upsilon_z[\tau\,\Xi]}{\tau!}\,\langle \Pi_z(\tau\,\Xi), \varphi_z^\lambda\rangle.$$

Using the decomposition (3.41) of the generalise gradient we write

$$
\begin{aligned}
\Upsilon_z[\boldsymbol{X}\,\Xi]\,\langle \Pi_z(\boldsymbol{X}\,\Xi), \varphi_z^\lambda\rangle &= u'(z)\,\Upsilon_z'[\boldsymbol{X}\,\Xi]\,\langle \Pi_z(\boldsymbol{X}\,\Xi), \varphi_z^\lambda\rangle \\
&= (\partial_{e_1} u_\lambda)(z)\,\Upsilon_z'[\boldsymbol{X}\,\Xi]\,\langle \Pi_z(\boldsymbol{X}\,\Xi), \varphi_z^\lambda\rangle \\
&\quad + \langle U_\gamma(z,\cdot) - u'(z)\,(\cdot - z)_1, \partial_{e_1}\varphi_z^\lambda\rangle\,\Upsilon_z'[\boldsymbol{X}\,\Xi]\,\langle \Pi_z(\boldsymbol{X}\,\Xi), \varphi_z^\lambda\rangle \\
&\quad - \sum_{\rho\in\mathcal{T}_{\gamma-2}} \frac{\Upsilon_z[\rho]}{\rho!}\,\langle \Pi_z\,\mathcal{I}'(\rho), \varphi_z^\lambda\rangle\,\Upsilon_z'[\boldsymbol{X}\,\Xi]\,\langle \Pi_z(\boldsymbol{X}\,\Xi), \varphi_z^\lambda\rangle,
\end{aligned}
$$

and set

$$b_\lambda^\Pi(z) := \Upsilon_z'[\boldsymbol{X}\,\Xi]\,\langle \Pi_z(\boldsymbol{X}\,\Xi), \varphi_z^\lambda\rangle. \tag{3.43}$$

Then, we obtain the identity

$$
\begin{aligned}
&\langle \Pi_z(\hat{\sigma}(U_\gamma)(z)\,\Xi), \varphi_z^\lambda\rangle - (\partial_{e_1} u_\lambda)(z)\,b_\lambda^\Pi(z) \\
&= \sum_{\tau\in\mathcal{T}_\gamma^+\setminus\{\boldsymbol{X}\}} \frac{\Upsilon_z[\tau\,\Xi]}{\tau!}\,\langle \Pi_z(\tau\,\Xi), \varphi_z^\lambda\rangle - \sum_{\rho\in\mathcal{T}_{\gamma-2}} \frac{\Upsilon_z[\rho]}{\rho!}\,\langle \Pi_z\,\mathcal{I}'(\rho), \varphi_z^\lambda\rangle\,\Upsilon_z'[\boldsymbol{X}\,\Xi]\,\langle \Pi_z(\boldsymbol{X}\,\Xi), \varphi_z^\lambda\rangle. \\
&\quad + \langle U_\gamma(z,\cdot) - u'(z)\,(\cdot - z)_1, \partial_{e_1}\varphi_z^\lambda\rangle\,\Upsilon_z'[\boldsymbol{X}\,\Xi]\,\langle \Pi_z(\boldsymbol{X}\,\Xi), \varphi_z^\lambda\rangle.
\end{aligned} \tag{3.44}
$$

We bound each of the terms in (3.44). For $\tau = \mathbf{1}$ we have the required estimate since $\Upsilon_z[\Xi] = \sigma(u(z))$ is bounded:

$$|\Upsilon_z[\Xi]\langle \Pi_z \Xi, \varphi_z^\lambda\rangle| \lesssim \lambda^{\alpha-2} \|(\Pi,\Gamma)\|.$$

Given $\tau \in \mathcal{T}_\gamma^+ \setminus \{\mathbf{1}, \boldsymbol{X}\}$ we have by Lemma 2.7 and Corollary 3.3

$$\begin{aligned}
\left|\frac{\Upsilon_z[\tau\,\Xi]}{\tau!}\langle \Pi_z(\tau\,\Xi), \varphi_z^\lambda\rangle\right| &\lesssim |u'(z)|^{|\mathfrak{n}(\tau)|}\,\lambda^{|\tau\Xi|}\,[\Pi;\tau\,\Xi]\\
&\lesssim (\lambda_t^{-1}\,\|u\|_{P_t})^{|\mathfrak{n}(\tau)|}\,\lambda^{|\tau\Xi|}\,|||(\Pi,\Gamma)|||^{\mathfrak{C}(\tau\Xi)}\\
&\sim \lambda^{|\tau\Xi|}\,\|u\|_{P_t}^{\frac{1}{\alpha}(\gamma-1+\alpha)|\mathfrak{n}(\tau)|}\,(1+|||(\Pi,\Gamma)|||)^{\frac{|\mathfrak{n}(\tau)|}{\alpha}}\,|||(\Pi,\Gamma)|||^{\mathfrak{C}(\tau\Xi)}\\
&\lesssim \lambda^{\alpha-2}\,|||(\Pi,\Gamma)|||\left(\lambda^{|\tau|}\,\|u\|_{P_t}^{\frac{1}{\alpha}(\gamma-1+\alpha)|\mathfrak{n}(\tau)|}\,(1+|||(\Pi,\Gamma)|||)^{\frac{|\tau|}{\alpha}}\right)\\
&= \lambda^{\alpha-2}\,|||(\Pi,\Gamma)|||\left(\lambda\,(1+|||(\Pi,\Gamma)|||)^{\frac{1}{\alpha}}\,\|u\|_{P_t}^{\frac{1}{\alpha}(\gamma-1+\alpha)\frac{|\mathfrak{n}(\tau)|}{|\tau|}}\right)^{|\tau|}\\
&= \lambda^{\alpha-2}\,|||(\Pi,\Gamma)|||\,\big(\lambda\,\lambda_t^{-1}\,\|u\|_{P_t}^{\mathcal{E}_\tau(\gamma)}\big)^{|\tau|},
\end{aligned} \tag{3.45}$$

where for $\tau \in \mathcal{T}_\gamma^+ \setminus \{\mathbf{1}, \boldsymbol{X}\}$ we define:

$$\mathcal{E}_\tau(\gamma) := \frac{1}{\alpha}\,(\gamma-1+\alpha)\,\frac{|\mathfrak{n}(\tau)|}{|\tau|} - \frac{(\gamma-1)}{\alpha} = \frac{\alpha\,|\mathfrak{n}(\tau)| - (\gamma-1)\,(|\tau| - |\mathfrak{n}(\tau)|)}{\alpha\,|\tau|}, \tag{3.46}$$

and where we used the identity

$$\frac{|\mathfrak{n}(\tau)|}{\alpha} + \mathfrak{C}(\tau\,\Xi) = 1 + \frac{1}{\alpha}(\alpha\,\mathfrak{C}(\tau) - |\mathfrak{n}(\tau)|) = 1 + \frac{|\tau|}{\alpha}. \tag{3.47}$$

Given $\rho \in \mathcal{T}_{\gamma-2}$ we have (recall that $\Upsilon'$ is a bounded function by Lemma 2.7)

$$\begin{aligned}
&\left|\frac{\Upsilon_z[\rho]}{\rho!}\langle \Pi_z \mathcal{I}'(\rho), \varphi_z^\lambda\rangle\,\Upsilon_z'[\boldsymbol{X}\,\Xi]\,\langle \Pi_z(\boldsymbol{X}\,\Xi), \varphi_z^\lambda\rangle\right|\\
\lesssim\;& |u'(z)|^{|\mathfrak{n}(\rho)|}\,[\Pi;\mathcal{I}'(\rho)]\,\lambda^{|\mathcal{I}'(\rho)|}\,[\Pi;\boldsymbol{X}\,\Xi]\,\lambda^{|\boldsymbol{X}\Xi|}\\
\lesssim\;& \lambda^{\alpha-2+|\mathcal{I}(\rho)|}\,|||(\Pi,\Gamma)|||^{\mathfrak{C}(\rho)+\mathfrak{C}(\boldsymbol{X}\Xi)}\,(\lambda_t^{-1}\,\|u\|_{P_t})^{|\mathfrak{n}(\rho)|}\\
\lesssim\;& \lambda^{\alpha-2}\,|||(\Pi,\Gamma)|||\left(\lambda^{|\mathcal{I}(\rho)|}\,(1+|||(\Pi,\Gamma)|||)^{\mathfrak{C}(\rho)}\left((1+|||(\Pi,\Gamma)|||)^{\frac{1}{\alpha}|\mathfrak{n}(\rho)|}\,\|u\|_{P_t}^{\frac{1}{\alpha}(\gamma-1+\alpha)|\mathfrak{n}(\rho)|}\right)\right)\\
=\;& \lambda^{\alpha-2}\,|||(\Pi,\Gamma)|||\left(\lambda\,(1+|||(\Pi,\Gamma)|||)^{\frac{\alpha\,\mathfrak{C}(\rho)+|\mathfrak{n}(\rho)|}{\alpha\,|\mathcal{I}(\rho)|}}\,\|u\|_{P_t}^{\frac{1}{\alpha}(\gamma-1+\alpha)\frac{|\mathfrak{n}(\rho)|}{|\mathcal{I}(\rho)|}}\right)^{|\mathcal{I}(\rho)|}\\
=\;& \lambda^{\alpha-2}\,|||(\Pi,\Gamma)|||\left(\lambda\,\lambda_t^{-1}\,\|u\|_{P_t}^{\mathcal{E}_{\mathcal{I}(\rho)}(\gamma)}\right)^{|\mathcal{I}(\rho)|},
\end{aligned} \tag{3.48}$$

where $\mathcal{E}_{\mathcal{I}(\rho)}(\gamma)$ is defined as in (3.46) since $\mathcal{I}(\rho) \in \mathcal{T}_\gamma^+$, and where we used that

$$\mathfrak{C}(\rho) = \mathfrak{C}(\mathcal{I}(\rho)) = \frac{1}{\alpha}\,(|\mathcal{I}(\rho)| - |\mathfrak{n}(\mathcal{I}(\rho))|).$$

For the last term in (3.44) we have by [MW20b, eq. (2.10)] and Theorem 3.4

$$\begin{aligned}
&|\langle U_\gamma(z,\cdot) - u'(z)(\cdot - z)_1, \partial_{e_1}\varphi_z^\lambda\rangle\,\Upsilon_z'[\boldsymbol{X}\,\Xi]\,\langle \Pi_z(\boldsymbol{X}\,\Xi), \varphi_z^\lambda\rangle|\\
\lesssim\;& \lambda^{\gamma-1}\,[U_\gamma;\mathbf{1}]_{B_\lambda(z)}\,\lambda^{|\boldsymbol{X}\Xi|}\,[\Pi;\boldsymbol{X}\,\Xi]\\
\lesssim\;& \lambda^{|\Xi|+\gamma}\,\lambda_t^{-\gamma}\,\|u\|_{P_t}\,|||(\Pi,\Gamma)|||\\
\sim\;& \lambda^{\alpha-2}\,|||(\Pi,\Gamma)|||\left(\lambda\,\lambda_t^{-1}\,\|u\|_{P_t}^{\frac{1}{\gamma}}\right)^{\gamma}.
\end{aligned} \tag{3.49}$$

We now choose $\lambda \in (0, \lambda_t]$ so that the factors in (3.45), (3.49), and (3.48) are bounded. Since $\|u\|_{P_t} \geqslant 1$, it is sufficient to require

$$\lambda \leqslant \lambda_t \|u\|_{P_t}^{-1/\gamma} \wedge \min_{\tau \in \mathcal{T}_\gamma^+ \setminus \{1, \mathbf{X}\}} \lambda_t \|u\|_{P_t}^{-\mathcal{E}_\tau(\gamma)}.$$

The contribution corresponding to $\tau = 1$ was treated separately above. Moreover, the exponents $\mathcal{E}_{\mathcal{I}(\rho)}(\gamma)$ appearing in (3.48) are included in the same maximum. A direct inspection of the finite set $\mathcal{T}_\gamma^+ \setminus \{1, \mathbf{X}\}$, using (3.46), gives

$$\max_{\tau \in \mathcal{T}_\gamma^+ \setminus \{1, \mathbf{X}\}} \mathcal{E}_\tau(\gamma) \leqslant \frac{2 - \gamma}{\alpha + 1}.$$

and therefore, since $\|u\|_{P_t} \geqslant 1$, we need to impose the following conditions on $\lambda$:

$$\lambda \lesssim \lambda_t \, \|u\|_{P_t}^{-\frac{1}{\gamma}} \wedge \min_{\tau \in \mathcal{T}_\gamma^+ \setminus \{\mathbf{1}, \boldsymbol{X}\}} \lambda_t \, \|u\|_{P_t}^{-\mathcal{E}_\tau(\gamma)} = \lambda_t \, \|u\|_{P_t}^{-\left(\frac{1}{\gamma} \vee \max_{\tau \in \mathcal{T}_\gamma^+ \setminus \{\mathbf{1}, \boldsymbol{X}\}} \mathcal{E}_\tau(\gamma)\right)},$$

where we used that for $\tau = \mathbf{1}$ in (3.45) the bound automatically holds since $|\mathbf{1}| = 0$. We claim that the minimum is attached at the first term, to see this observe that $\mathcal{E}_\tau(\gamma)$ as in (3.46) is decreasing in $|\tau|$ and therefore the maximum for $\tau \in \mathcal{T}_\gamma^+ \setminus \{\mathbf{1}, \boldsymbol{X}\}$ is attached at $\rho = \mathcal{I}(\Xi)$ or $\rho = \mathcal{I}(\boldsymbol{X}\,\Xi)$, i.e.,

$$\max_{\tau \in \mathcal{T}_\gamma^+ \setminus \{\mathbf{1}, \boldsymbol{X}\}} \mathcal{E}_\tau(\gamma) = \mathcal{E}_{\mathcal{I}(\Xi)}(\gamma) \vee \mathcal{E}_{\mathcal{I}(\Xi \boldsymbol{X})} = \frac{-(\gamma - 1)\,\alpha}{\alpha^2} \vee \frac{\alpha - (\gamma - 1)\,\alpha}{\alpha\,(\alpha + 1)} = \frac{2 - \gamma}{\alpha + 1},$$

since the first term is negative. Since $\gamma > 2 - \alpha$, $\gamma < 2$, and $\alpha < 1$, we have

$$\gamma\,(2 - \gamma) < 2\,\alpha < \alpha + 1,$$

and therefore

$$\frac{2 - \gamma}{\alpha + 1} < \frac{1}{\gamma}.$$

Consequently,

$$\frac{1}{\gamma} \vee \max_{\tau \in \mathcal{T}_\gamma^+ \setminus \{1, \mathbf{X}\}} \mathcal{E}_\tau(\gamma) = \frac{1}{\gamma}.$$

Thus, whenever $0 < \lambda \leqslant \lambda_t \|u\|_{P_t}^{-1/\gamma}$, all three terms are bounded by $\lambda^{\alpha - 2} \|(\Pi, \Gamma)\|$, which proves the lemma. □

For the proof of Lemma 3.15 we will need the following multi-scale decomposition of the reconstruction operator. Recall that given a modelled distribution $F \in \mathcal{D}^\gamma$ for $\gamma > 0$, we have the identity

$$\langle \mathcal{R}\,F - \Pi_x\,F(x), \varphi_x^\lambda \rangle = \sum_{n=0}^{\infty} \iint \left\langle \Pi_z\,F(z) - \Pi_y\,F(y), \varphi_z^{\lambda 2^{-(n+1)}} \right\rangle \psi_y^{\lambda 2^{-(n+1)}}(z)\,\varphi_x^{\lambda, n}(y)\,\mathrm{d}y\,\mathrm{d}z, \tag{3.50}$$

and since

$$\Pi_z\,F(z) - \Pi_y\,F(y) = \Pi_z(F(z) - \Gamma_{zy}\,F(y)) = \sum_{\tau \in \mathcal{T}} \langle \tau, F(z) - \Gamma_{zy}\,F(y) \rangle\,\Pi_z\,\tau,$$

we can conclude that

$$\begin{aligned}
&\langle \mathcal{R}\,F - \Pi_x\,F(x), \varphi_x^\lambda \rangle \\
= &\sum_{\tau \in \mathcal{T}} \sum_{n=0}^{\infty} \iint \langle \tau, F(z) - \Gamma_{zy}\,F(y) \rangle \left\langle \Pi_z\,\tau, \varphi_z^{\lambda 2^{-(n+1)}} \right\rangle \psi_y^{\lambda 2^{-(n+1)}}(z)\,\varphi_x^{\lambda, n}(y)\,\mathrm{d}y\,\mathrm{d}z \\
= &\sum_{\tau \in \mathcal{T}} \sum_{n=0}^{\infty} \Lambda_{n, \lambda}[H^\tau; \Pi\,\tau](x),
\end{aligned} \tag{3.51}$$

where $\Lambda_{n,\lambda}$ takes as an argument a pair $(H,G)$ where $H$ is function valued germ, $G$ is a distribution valued germ and

$$\Lambda_{n,\lambda}[H,G](x) := \iint H(y,z)\,\langle G(z), \varphi_z^{\lambda 2^{-(n+1)}}\rangle \psi_y^{\lambda 2^{-(n+1)}}(z)\,\varphi_x^{\lambda,n}(y)\,\mathrm{d}y\mathrm{d}z. \tag{3.52}$$

The particular choice of the germ $H^\tau$ in (3.51) is chosen as

$$H^\tau(y,z) = \langle \tau, F(z) - \Gamma_{zy} F(y)\rangle.$$

By the reconstruction theorem [MW20b, Theorem 2.8] we have that

$$\Lambda_\lambda[H,G] := \sum_{n=0}^{\infty} \Lambda_{n,\lambda}[H,G] \tag{3.53}$$

is well defined and satisfies the following estimate

$$|\Lambda_\lambda[H,G](x)| \lesssim \lambda^{\beta_1+\beta_2}\,\|G\|_{\beta_1;B_\lambda(x)}\,[H]_{\beta_2;B_\lambda(x)} \tag{3.54}$$

as long as the seminorms on the right are finite for some $\beta_1+\beta_2>0$.

**Proof of Lemma 3.15.** The reconstruction theorem in the form of (3.51) allows us to write

$$\langle \mathcal{R}(\hat{\sigma}(U_\gamma)\,\Xi) - \Pi_z(\hat{\sigma}(U_\gamma)(z)\,\Xi), \varphi_z^\lambda\rangle = \sum_{\tau\in\mathcal{T}_\gamma^+} \Lambda_\lambda[H^\tau, \Pi(\tau\,\Xi)](z)$$

where the reconstruction operator $\Lambda_\lambda$ operator is defined as in (3.53) and

$$H^\tau(x,y) := \langle \tau\,\Xi, \hat{\sigma}(U_\gamma(y))\,\Xi - \Gamma_{yx}\,\hat{\sigma}(U_\gamma(x))\,\Xi\rangle = \langle \tau, \hat{\sigma}(U_\gamma(y)) - \Gamma_{yx}\,\hat{\sigma}(U_\gamma(x))\rangle.$$

By the results of Lemma 2.5 we have that the operator $\Lambda_\lambda[\,\cdot\,,\Pi(\tau\,\Xi)]$ is well defined for any of the germs coming from the decomposition of $H^\tau$ in Lemma 2.9. We start by identifying the function $b_\lambda^{\mathcal{R}}$. The germ $H_\emptyset^{\mathbf{1}}$ as defined in (2.20) for $\tau=\mathbf{1}$ and $n=0$ can be expressed as a Taylor remainder:

$$\begin{aligned} H_\emptyset^{\mathbf{1}}(x,y) &= \Upsilon[\Xi](u(x)+u'(x)\,(y-x)_1) - \Upsilon[\Xi](u(x)) - (\Upsilon[\Xi])^{(1)}(u(x))\,u'(x)\,(y-x)_1\\ &= \sigma(u(x)+u'(x)\,(y-x)_1) - \sigma(u(x)) - \sigma^{(1)}(u(x))\,u'(x)\,(y-x)_1\\ &= (u'(x)\,(y-x)_1)^2 \int_0^1 (1-\lambda)\;\sigma^{(2)}\,(u(x)+\lambda\,u'(x)\,(y-x)_1)\,\mathrm{d}\lambda\\ &=: u'(x)\,\tilde{H}^{\mathbf{1}}(x,y), \end{aligned} \tag{3.55}$$

which leads to the decomposition of $H^{\mathbf{1}}$ as

$$\begin{aligned} H^{\mathbf{1}}(x,y) - (\partial_{e_1} u_\lambda)(z)\,\tilde{H}^{\mathbf{1}}(x,y) &= H_1^{\mathbf{1}}(x,y) + \sum_{\beta\in\partial A_{\mathbf{1}}} H_\beta^{\mathbf{1}}(x,y) + \sum_{n\neq 0} \frac{1}{n!} \sum_{\substack{\tau_1,\ldots,\tau_n\in\mathcal{T}\\ \sum_{i=1}^n |\mathcal{I}(\tau_i)|<\gamma}} H^{\mathbf{1}}_{(\tau_i)_{i=1}^n}(x,y)\\ &\quad + H_\emptyset^{\mathbf{1}}(x,y) - (\partial_{e_1} u_\lambda)(z)\,\tilde{H}^{\mathbf{1}}(x,y). \end{aligned} \tag{3.56}$$

On $B_\lambda(z)\subset P_t$ we have the estimate

$$\begin{aligned} |\tilde{H}^{\mathbf{1}}(x,y)| &\lesssim \left| u'(x)\,((y-x)_1)^2 \int_0^1 (1-\lambda)\;\sigma^{(2)}\,(u(x)+\lambda\,u'(x)\,(y-x)_1)\,\mathrm{d}\lambda\right|\\ &\lesssim \|u'\|_{B_\lambda(z)}\,d(x,y)^2 \lesssim \lambda_t^{-1}\,\|u\|_{P_t}\,d(x,y)^2, \end{aligned} \tag{3.57}$$

where we used Corollary 3.3 to bound the generalised derivative. Since $2+|\Xi|=\alpha>0$ we have that

$$b_\lambda^{\mathcal{R},\mathbf{1}}(z):=\Lambda_\lambda[\tilde{H}^{\mathbf{1}},\Pi\Xi](z) \tag{3.58}$$

is well-defined and the identity

$$\Lambda_\lambda[(\partial_{e_1}u_\lambda)(z)\,\tilde{H}^{\mathbf{1}},\Pi\Xi](z)=(\partial_{e_1}u_\lambda)(z)\,\Lambda_\lambda[\tilde{H}^{\mathbf{1}},\Pi\Xi](z)$$

holds. The other term contributing to the transport term will come from $\tau=\boldsymbol{X}$. Consider the germ $H_2^{\boldsymbol{X}}$ as defined in Lemma 2.9 and write it as a Taylor remainder

$$\begin{aligned}
& H_2^{\boldsymbol{X}}(x,y)\\
= \;& \Upsilon'[\boldsymbol{X}\Xi](u(y)-\langle\mathbf{1},U_\gamma(y)-\Gamma_{yx}U_\gamma(y)\rangle)-\Upsilon'[\boldsymbol{X}\Xi](u(y)-\langle\mathbf{1},U_{\gamma-1}(y)-\Gamma_{yx}U_{\gamma-1}(y)\rangle)\\
= \;& \sigma^{(1)}(\langle\mathbf{1},\Gamma_{yx}U_\gamma(y)\rangle)-\sigma^{(1)}(\langle\mathbf{1},\Gamma_{yx}U_{\gamma-1}(y)\rangle)\\
= \;& \langle\mathbf{1},\Gamma_{yx}(U_\gamma(x)-U_{\gamma-1}(x))\rangle\int_0^1\sigma^{(2)}(\langle\mathbf{1},\Gamma_{yx}(\lambda\,U_\gamma(x)+(1-\lambda)\,U_{\gamma-1}(x))\rangle)\,\mathrm{d}\lambda\\
= \;& \left(u'(x)\,(y-x)_1+\sum_{\rho\in\mathcal{T}_{\gamma-3,\gamma-2}}\frac{\Upsilon_x[\rho]}{\rho!}\,\gamma_{yx}(\mathcal{I}(\rho))\right)\tilde{H}_2^{\boldsymbol{X}}(x,y)
\end{aligned}$$

where $\tilde{H}_2^{\boldsymbol{X}}(x,y)=\int_0^1\sigma^{(2)}(\langle\mathbf{1},\Gamma_{yx}(\lambda\,U_\gamma(x)+(1-\lambda)\,U_{\gamma-1}(x))\rangle)\,\mathrm{d}\lambda$. We define

$$\hat{H}^{\boldsymbol{X}}(x,y):=u'(x)\,(y-x)_1\,\tilde{H}_2^{\boldsymbol{X}}(x,y),$$

which leads to the decomposition of $H^{\boldsymbol{X}}$ as

$$\begin{aligned}
& H^{\boldsymbol{X}}(x,y)-(\partial_{e_1}u_\lambda)(z)\,\hat{H}^{\boldsymbol{X}}(x,y)\\
= \;& u'(y)\left(H_1^{\boldsymbol{X}}(x,y)+\sum_{\beta\in\partial A_{\boldsymbol{X}}}H_\beta^{\boldsymbol{X}}(x,y)\right)+\sum_{n\in\mathbb{N}}\frac{1}{n!}\sum_{\substack{\tau_1,\ldots,\tau_n\in\mathcal{T}\\ \sum_{i=1}^n|\mathcal{I}(\tau_i)|<\gamma-1}}H^{\tau}_{(\tau_i)_{i=1}^n}(x,y)\\
& +u'(y)\left(u'(x)\,(y-x)_1+\sum_{\rho\in\mathcal{T}_{\gamma-3,\gamma-2}}\frac{\Upsilon_x[\rho]}{\rho!}\,\gamma_{yx}(\mathcal{I}(\rho))\right)\tilde{H}_2^{\boldsymbol{X}}(x,y)-(\partial_{e_1}u_\lambda)(z)\,\hat{H}^{\boldsymbol{X}}(x,y)\\
= \;& u'(y)\left(H_1^{\boldsymbol{X}}(x,y)+\sum_{\beta\in\partial A_{\boldsymbol{X}}}H_\beta^{\boldsymbol{X}}(x,y)\right)+\sum_{n\in\mathbb{N}}\frac{1}{n!}\sum_{\substack{\tau_1,\ldots,\tau_n\in\mathcal{T}\\ \sum_{i=1}^n|\mathcal{I}(\tau_i)|<\gamma-1}}H^{\tau}_{(\tau_i)_{i=1}^n}(x,y)\\
& +u'(y)\sum_{\rho\in\mathcal{T}_{\gamma-3,\gamma-2}}\frac{\Upsilon_x[\rho]}{\rho!}\,\gamma_{yx}(\mathcal{I}(\rho))\,\tilde{H}_2^{\boldsymbol{X}}(x,y)+(u'(z)-(\partial_{e_1}u_\lambda)(z)+u'(y)-u'(z))\,\hat{H}^{\boldsymbol{X}}(x,y).
\end{aligned} \tag{3.59}$$

On $B_\lambda(z)\subset P_t$ we have the estimate

$$|\hat{H}^{\boldsymbol{X}}(x,y)|\lesssim\|u'\|_{B_\lambda(z)}\,d(x,y)|\tilde{H}_2^{\boldsymbol{X}}(x,y)|\lesssim\|u'\|_{B_\lambda(z)}\,d(x,y)\lesssim\lambda_t^{-1}\,\|u\|_{P_t}\,d(x,y), \tag{3.60}$$

and since $1+|\boldsymbol{X}\,\Xi|=\alpha>0$ we have that

$$b_\lambda^{\mathcal{R},\boldsymbol{X}}(z):=\Lambda_\lambda[\hat{H}^{\boldsymbol{X}},\Pi(\boldsymbol{X}\,\Xi)] \tag{3.61}$$

is well-defined and the following identity holds:

$$\Lambda_\lambda[(\partial_{e_1}u_\lambda)(z)\,\hat{H}^{\boldsymbol{X}},\Pi(\boldsymbol{X}\,\Xi)](z)=(\partial_{e_1}u_\lambda)(z)\,\Lambda_\lambda[\hat{H}^{\boldsymbol{X}},\Pi\,(\boldsymbol{X}\,\Xi)](z).$$

We set

$$b_\lambda^{\mathcal{R}}(z):=b_\lambda^{\mathcal{R},\mathbf{1}}(z)+b_\lambda^{\mathcal{R},\boldsymbol{X}}(z) \tag{3.62}$$

and obtain the identity

$$\begin{aligned} & \langle \mathcal{R}(\hat{\sigma}(U_\gamma)\,\Xi) - \Pi_z(\hat{\sigma}(U_\gamma)(z)\,\Xi), \varphi_z^\lambda \rangle - (\partial_{e_1} u_\lambda)(z)\, b_\lambda^{\mathcal{R}}(z) \\ = & \sum_{\tau \in \mathcal{T}_\gamma^+ \setminus \{\mathbf{1}, \boldsymbol{X}\}} \Lambda_\lambda[H^\tau, \Pi(\tau\,\Xi)](z) \\ & + \Lambda_\lambda[H^{\mathbf{1}} - (\partial_{e_1} u_\lambda)(z)\, \tilde{H}^{\mathbf{1}}, \Pi(\Xi)](z) + \Lambda_\lambda[H^{\boldsymbol{X}} - (\partial_{e_1} u_\lambda)(z)\, \hat{H}^{\boldsymbol{X}}, \Pi(\boldsymbol{X}\,\Xi)](z). \end{aligned}$$

We proceed to estimate each term. For any $\tau \in \mathcal{T}_\gamma^+ \setminus \{\mathbf{1}, \boldsymbol{X}\}$ we have by (3.54)

$$\begin{aligned} |\Lambda_\lambda[H^\tau, \Pi(\tau\,\Xi)](z)| & \lesssim \lambda^{\gamma - |\tau| + |\tau\Xi|}\, [\hat{\sigma}(U_\gamma); \tau]_{B_\lambda(z)}\, [\Pi; \tau\,\Xi] \\ & \lesssim \lambda^{\alpha - 2}\, \lambda^\gamma\, [\hat{\sigma}(U_\gamma); \tau]_{B_\lambda(z)}\, |||(\Pi, \Gamma)|||^{\mathfrak{C}(\tau\Xi)} \\ & \lesssim \lambda^{\alpha - 2}\, |||(\Pi, \Gamma)|||\, \lambda^\gamma\, [\hat{\sigma}(U_\gamma); \tau]_{B_\lambda(z)}\, (1 + |||(\Pi, \Gamma)|||)^{\mathfrak{C}(\tau)} \end{aligned} \tag{3.63}$$

If $\tau$ has no decorations, then by Corollary 3.3

$$\begin{aligned} & \lambda^\gamma\, [\hat{\sigma}(U_\gamma); \tau]_{B_\lambda(z)}\, (1 + |||(\Pi, \Gamma)|||)^{\mathfrak{C}(\tau)} \\ \lesssim & \lambda^\gamma \Big( \lambda_t^{-(\gamma - |\tau|)}\, \|u\|_{P_t}^{1 + \mathbb{1}_{\gamma - |\tau| > 1} (\gamma - |\tau|)(1 - \gamma^{-1})} \Big)\, (1 + |||(\Pi, \Gamma)|||)^{\mathfrak{C}(\tau)} \\ \sim & \lambda^\gamma \Big( (1 + |||(\Pi, \Gamma)|||)^{\frac{\gamma - |\tau|}{\alpha}}\, \|u\|_{P_t}^{\frac{(\gamma - 1)}{\alpha}(\gamma - |\tau|)}\, \|u\|_{P_t}^{1 + \mathbb{1}_{\gamma - |\tau| > 1} (\gamma - |\tau|)(1 - \gamma^{-1})} \Big)\, (1 + |||(\Pi, \Gamma)|||)^{\mathfrak{C}(\tau)} \\ \lesssim & \Big( \lambda\, (1 + |||(\Pi, \Gamma)|||)^{\frac{1}{\alpha}}\, \|u\|_{P_t}^{\frac{1}{\alpha}(\gamma - 1) + \mathcal{E}^\tau(\gamma)} \Big)^\gamma \\ = & \big( \lambda\, \lambda_t^{-1}\, \|u\|_{P_t}^{\mathcal{E}^\tau(\gamma)} \big)^\gamma \end{aligned} \tag{3.64}$$

where

$$\mathcal{E}^\tau(\gamma) := \left( \frac{(\gamma - 1)}{\alpha}\, (\gamma - |\tau|) + 1 + \mathbb{1}_{\gamma - |\tau| > 1}\, (\gamma - |\tau|)\, (1 - \gamma^{-1}) \right) \gamma^{-1} - \frac{(\gamma - 1)}{\alpha}.$$

If $\tau$ has decorations, then by Corollary 3.3

$$\begin{aligned} & \lambda^\gamma\, [\hat{\sigma}(U_\gamma); \tau]_{B_\lambda(z)}\, |||(\Pi, \Gamma)|||^{\mathfrak{C}(\tau)} \\ \lesssim & \lambda^\gamma \big( \lambda_t^{-(\gamma - |\tau| + 1)}\, \|u\|_{P_t}^{1 + (\gamma - |\tau| + 1)(1 - \gamma^{-1})} \big)\, |||(\Pi, \Gamma)|||^{\mathfrak{C}(\tau)} \\ \sim & \lambda^\gamma \Big( (1 + |||(\Pi, \Gamma)|||)^{\frac{\gamma - |\tau| + 1}{\alpha}}\, \|u\|_{P_t}^{\frac{(\gamma - 1)}{\alpha}(\gamma - |\tau| + 1)}\, \|u\|_{P_t}^{1 + (\gamma - |\tau| + 1)(1 - \gamma^{-1})} \Big)\, |||(\Pi, \Gamma)|||^{\frac{|\tau| - 1}{\alpha}} \\ = & \Big( \lambda\, (1 + |||(\Pi, \Gamma)|||)^{\frac{1}{\alpha}}\, \|u\|_{P_t}^{\frac{1}{\alpha}(\gamma - 1) + \mathcal{E}^\tau(\gamma)} \Big)^\gamma \\ = & \big( \lambda\, \lambda_t^{-1}\, \|u\|_{P_t}^{\mathcal{E}^\tau(\gamma)} \big)^\gamma \end{aligned} \tag{3.65}$$

where

$$\mathcal{E}^\tau(\gamma) := \left( \frac{(\gamma - 1)}{\alpha}\, (\gamma - |\tau| + 1) + 1 + (\gamma - |\tau| + 1)\, (1 - \gamma^{-1}) \right) \gamma^{-1} - \frac{(\gamma - 1)}{\alpha}.$$

From (3.64) and (3.65) we can conclude on (3.63)

$$\sum_{\tau \in \mathcal{T}_\gamma^+ \setminus \{\mathbf{1}, \boldsymbol{X}\}} |\Lambda_\lambda[H^\tau, \Pi(\tau\,\Xi)](z)| \lesssim \lambda^{\alpha - 2}\, |||(\Pi, \Gamma)||| \sum_{\tau \in \mathcal{T}_\gamma^+ \setminus \{\mathbf{1}, \boldsymbol{X}\}} \big( \lambda\, \lambda_t^{-1}\, \|u\|_{P_t}^{\mathcal{E}^\tau(\gamma)} \big)^\gamma. \tag{3.66}$$

For the term coming from $\tau = \mathbf{1}$ we need to estimate (3.56). For the first three terms we have by (2.30), (2.32), and (2.41) and Corollary 3.3

$$\begin{aligned} |H_1^{\mathbf{1}}(x, y)| & \lesssim \lambda_t^{-\gamma}\, \|u\|_{P_t}\, d(x, y)^\gamma, \\ |H_\beta^{\mathbf{1}}(z, y)| & \lesssim \lambda^{|\mathfrak{a}(\beta)| - \gamma}\, (\lambda_t^{-1}\, \|u\|_{P_t})^{|\mathfrak{n}(\mathfrak{a}(\beta))|}\, |||(\Pi, \Gamma)|||^{\mathfrak{C}(\mathfrak{a}(\beta))}\, d(x, y)^\gamma, \\ |H_{(\tau_i)_{i=1}^n}^{\mathbf{1}}(x, y)| & \lesssim (\lambda_t^{-1}\, \|u\|_{P_t})^{|\mathfrak{n}(\prod_{i=1}^n \mathcal{I}(\tau_i))| + \gamma - |\prod_{i=1}^n \mathcal{I}(\tau_i)|}\, |||(\Pi, \Gamma)|||^{\mathfrak{C}(\prod_{i=1}^n \mathcal{I}(\tau_i))}\, d(x, y)^\gamma, \end{aligned}$$

and since by assumption $\gamma + |\Xi| = \gamma + \alpha - 2 > 0$ then

$$\begin{aligned}
&\left| \Lambda_\lambda[H_1^{\mathbf{1}}, \Pi\Xi] + \sum_{\beta \in \partial A_{\mathbf{1}}} \Lambda_\lambda[H_\beta^{\mathbf{1}}, \Pi\Xi] + \sum_{n \neq 0} \frac{1}{n!} \sum_{\substack{\tau_1, \ldots, \tau_n \in \mathcal{T} \\ \sum_{i=1}^n |\mathcal{I}(\tau_i)| < \gamma}} \Lambda_\lambda[H^{\mathbf{1}}_{(\tau_i)_{i=1}^n}, \Pi\Xi] \right| \\
\lesssim\ & \lambda^{\gamma+\alpha-2} [\Pi; \Xi]\, \lambda_t^{-\gamma} \|u\|_{P_t} + \sum_{\beta \in \partial A_{\mathbf{1}}} \lambda^{\gamma+\alpha-2} [\Pi; \Xi]\, \lambda^{|a(\beta)|-\gamma} (\lambda_t^{-1} \|u\|_{P_t})^{|\mathfrak{n}(a(\beta))|} \, |||(\Pi, \Gamma)|||^{\mathfrak{E}(a(\beta))} \\
&+ \sum_{\beta \in A_{\mathbf{1}} \setminus \{0\}} \lambda^{\gamma+\alpha-2} [\Pi; \Xi]\, (\lambda_t^{-1} \|u\|_{P_t})^{|\mathfrak{n}(a(\beta))| + \gamma - |a(\beta)|} \, |||(\Pi, \Gamma)|||^{\mathfrak{E}(a(\beta))} \\
\lesssim\ & \lambda^{\alpha-2} |||(\Pi, \Gamma)|||\, \lambda^\gamma \lambda_t^{-\gamma} \|u\|_{P_t} + \lambda^{\alpha-2} |||(\Pi, \Gamma)||| \sum_{\tau \in a(\partial A_{\mathbf{1}})} \lambda^{|\tau|} (\lambda_t^{-1} \|u\|_{P_t})^{|\mathfrak{n}(\tau)|} \, |||(\Pi, \Gamma)|||^{\mathfrak{E}(\tau)} \\
&+ \lambda^{\alpha-2} |||(\Pi, \Gamma)||| \sum_{\tau \in a(A_{\mathbf{1}} \setminus \{0\})} \lambda^\gamma (\lambda_t^{-1} \|u\|_{P_t})^{|\mathfrak{n}(\tau)| + \gamma - |\tau|} \, |||(\Pi, \Gamma)|||^{\mathfrak{E}(\tau)} \\
\lesssim\ & \lambda^{\alpha-2} |||(\Pi,\Gamma)||| \left( \lambda \lambda_t^{-1} \|u\|_{P_t}^{\frac{1}{\gamma}} \right)^\gamma + \lambda^{\alpha-2} |||(\Pi,\Gamma)||| \sum_{\tau \in a(\partial A_{\mathbf{1}})} \left( \lambda (\lambda_t^{-1} \|u\|_{P_t})^{\frac{|\mathfrak{n}(\tau)|}{|\tau|}} (1 + |||(\Pi,\Gamma)|||)^{\frac{\mathfrak{E}(\tau)}{|\tau|}} \right)^{|\tau|} \\
&+ \lambda^{\alpha-2} |||(\Pi, \Gamma)||| \sum_{\tau \in a(A_{\mathbf{1}}) \setminus \{\mathbf{1}\}} \left( \lambda (\lambda_t^{-1} \|u\|_{P_t})^{\frac{|\mathfrak{n}(\tau)| + \gamma - |\tau|}{\gamma}} (1 + |||(\Pi, \Gamma)|||)^{\frac{\mathfrak{E}(\tau)}{\gamma}} \right)^\gamma. \qquad (3.67)
\end{aligned}$$

For $\tau \in a(\partial A_{\mathbf{1}})$ we have

$$\begin{aligned}
\lambda (\lambda_t^{-1} \|u\|_{P_t})^{\frac{|\mathfrak{n}(\tau)|}{|\tau|}} (1 + |||(\Pi, \Gamma)|||)^{\frac{\mathfrak{E}(\tau)}{|\tau|}} &\sim \lambda \|u\|_{P_t}^{\frac{1}{\alpha}(\gamma - 1 + \alpha) \frac{|\mathfrak{n}(\tau)|}{|\tau|}} (1 + |||(\Pi, \Gamma)|||)^{\frac{1}{\alpha} \frac{|\mathfrak{n}(\tau)|}{|\tau|} + \frac{\mathfrak{E}(\tau)}{|\tau|}} \\
&= \lambda \|u\|_{P_t}^{\frac{1}{\alpha}(\gamma - 1 + \alpha) \frac{|\mathfrak{n}(\tau)|}{|\tau|}} (1 + |||(\Pi, \Gamma)|||)^{\frac{|\mathfrak{n}(\tau)|}{\alpha |\tau|} + \frac{|\tau| - |\mathfrak{n}(\tau)|}{\alpha |\tau|}} \\
&= \lambda \|u\|_{P_t}^{\mathcal{E}_\tau^{\mathbf{1}}(\gamma) + \frac{(\gamma - 1)}{\alpha}} (1 + |||(\Pi, \Gamma)|||)^{\frac{1}{\alpha}} \\
&= \lambda \lambda_t^{-1} \|u\|_{P_t}^{\mathcal{E}_\tau^{\mathbf{1}}(\gamma)},
\end{aligned}$$

where

$$\mathcal{E}_\tau^{\mathbf{1}}(\gamma) := \frac{1}{\alpha}(\gamma - 1 + \alpha) \frac{|\mathfrak{n}(\tau)|}{|\tau|} - \frac{(\gamma - 1)}{\alpha} = \frac{\alpha\, |\mathfrak{n}(\tau)| - (\gamma - 1)\, (|\tau| - |\mathfrak{n}(\tau)|)}{\alpha\, |\tau|}.$$

Similarly, for $\tau \in a(A_{\mathbf{1}}) \setminus \{\mathbf{1}\}$ we have

$$\begin{aligned}
\lambda (\lambda_t^{-1} \|u\|_{P_t})^{\frac{|\mathfrak{n}(\tau)| + \gamma - |\tau|}{\gamma}} (1 + |||(\Pi, \Gamma)|||)^{\frac{\mathfrak{E}(\tau)}{\gamma}} &\sim \lambda \|u\|_{P_t}^{\frac{1}{\alpha}(\gamma - 1 + \alpha) \frac{|\mathfrak{n}(\tau)| + \gamma - |\tau|}{\gamma}} (1 + |||(\Pi, \Gamma)|||)^{\frac{1}{\alpha} \frac{|\mathfrak{n}(\tau)|}{\gamma} + \frac{\mathfrak{E}(\tau)}{\gamma}} \\
&= \lambda \lambda_t^{-1} \|u\|_{P_t}^{\tilde{\mathcal{E}}_\tau^{\mathbf{1}}(\gamma)}
\end{aligned}$$

where

$$\tilde{\mathcal{E}}_\tau^{\mathbf{1}}(\gamma) := \frac{1}{\alpha}(\gamma - 1 + \alpha) \frac{|\mathfrak{n}(\tau)| + \gamma - |\tau|}{\gamma} - \frac{(\gamma - 1)}{\alpha} = 1 - \frac{(\gamma - 1 + \alpha)\, (|\tau| - |\mathfrak{n}(\tau)|)}{\alpha\, \gamma}.$$

We conclude on (3.67)

$$\begin{aligned}
&\left| \Lambda_\lambda[H_1^{\mathbf{1}}, \Pi\Xi] + \sum_{\beta \in \partial A_{\mathbf{1}}} \Lambda_\lambda[H_\beta^{\mathbf{1}}, \Pi\Xi] + \sum_{n \neq 0} \frac{1}{n!} \sum_{\substack{\tau_1, \ldots, \tau_n \in \mathcal{T} \\ \sum_{i=1}^n |\mathcal{I}(\tau_i)| < \gamma}} \Lambda_\lambda[H^{\mathbf{1}}_{(\tau_i)_{i=1}^n}, \Pi\Xi] \right| \qquad (3.68) \\
\lesssim\ & \lambda^{\alpha-2} |||(\Pi,\Gamma)||| \left( \left( \lambda \lambda_t^{-1} \|u\|_{P_t}^{\frac{1}{\gamma}} \right)^\gamma + \sum_{\tau \in a(\partial A_{\mathbf{1}})} \left( \lambda \lambda_t^{-1} \|u\|_{P_t}^{\mathcal{E}_\tau(\gamma)} \right)^{|\tau|} + \sum_{\tau \in a(A_{\mathbf{1}}) \setminus \{\mathbf{1}\}} \left( \lambda \lambda_t^{-1} \|u\|_{P_t}^{\tilde{\mathcal{E}}_\tau(\gamma)} \right)^{|\tau|} \right).
\end{aligned}$$

For the last term in (3.56) we use (3.55), then do the replacement $u'(x)=u'(z)+u'(x)-u'(z)$ to fix the basepoint, and replace $u'(z)$ with the real derivative of the regularisation using (3.41) with which we obtain

$$
\begin{aligned}
& |H^{\mathbf{1}}_{\emptyset}(x,y)-(\partial_{e_1}u_\lambda)(z)\,\tilde{H}^{\mathbf{1}}(x,y)|\\
=\; & |(u'(x)-(\partial_{e_1}u_\lambda)(z))\,\tilde{H}^{\mathbf{1}}(x,y)|\\
=\; & \left|\langle U_\gamma(z,\cdot)-u'(z)\,(\cdot-z)_1,\partial_{e_1}\varphi_z^\lambda\rangle-\sum_{\rho\in\mathcal{T}_{\gamma-2}}\frac{\Upsilon_z[\rho]}{\rho!}\,\langle\Pi_z\mathcal{I}'(\rho),\varphi_z^\lambda\rangle+(u'(x)-u'(z))\right|\\
\lesssim\; & |\langle U_\gamma(z,\cdot)-u'(z)\,(\cdot-z)_1,\partial_{e_1}\varphi_z^\lambda\rangle|+\sum_{\rho\in\mathcal{T}_{\gamma-2}}|\Upsilon_z[\rho]\,\langle\Pi_z\mathcal{I}'(\rho),\varphi_z^\lambda\rangle|+|u'(x)-u'(z)|\\
\lesssim\; & \lambda^{\gamma-1}\,[U_\gamma;\mathbf{1}]_{B_\lambda(z)}+\sum_{\rho\in\mathcal{T}_{\gamma-2}}\|u'\|_{B_\lambda(z)}^{|\mathfrak{n}(\rho)|}\,\lambda^{|\mathcal{I}'(\rho)|}\,[\Pi;\mathcal{I}'(\rho)]+\lambda^{\gamma-1}\,[U_\gamma;\boldsymbol{X}]_{B_\lambda(z)}\\
& +\sum_{\rho\in\mathcal{T}_{-1,\gamma-2}}\|u'\|_{B_\lambda(z)}^{|\mathfrak{n}(\rho)|}\,\lambda^{|\mathcal{I}'(\rho)|}\,[\Gamma\,\mathcal{I}(\rho);\boldsymbol{X}]\\
\lesssim\; & \lambda^{\gamma-1}\,\lambda_t^{-\gamma}\,\|u\|_{P_t}+\sum_{\rho\in\mathcal{T}_{\gamma-2}}(\lambda_t^{-1}\,\|u\|_{P_t})^{|\mathfrak{n}(\rho)|}\,\lambda^{|\mathcal{I}'(\rho)|}\,|||(\Pi,\Gamma)|||^{\mathfrak{C}(\rho)}+\lambda^{\gamma-1}\,\lambda_t^{-\gamma}\,\|u\|_{P_t}\\
\lesssim\; & \lambda^{\gamma-1}\,\lambda_t^{-\gamma}\,\|u\|_{P_t}+\sum_{\rho\in\mathcal{T}_{\gamma-2}}(\lambda_t^{-1}\,\|u\|_{P_t})^{|\mathfrak{n}(\rho)|}\,\lambda^{|\mathcal{I}'(\rho)|}\,|||(\Pi,\Gamma)|||^{\mathfrak{C}(\rho)},
\end{aligned}
\tag{3.69}
$$

where we used Corollary 3.3 and that by [EW24, (2.22) from Lemma 2.19] the generalised gradient $u'$ can be modelled to order $\gamma-1$ as

$$
\begin{aligned}
|u'(x)-u'(z)|\;\lesssim\; & \left|u'(x)-u'(z)-\sum_{\rho\in\mathcal{T}_{-1,\gamma-2}}\frac{\Upsilon_z[\rho]}{\rho!}\,\gamma_{xz}(\mathcal{I}'(\rho))\right|+\sum_{\rho\in\mathcal{T}_{-1,\gamma-2}}\left|\frac{\Upsilon_z[\rho]}{\rho!}\,\gamma_{xz}(\mathcal{I}'(\rho))\right|\\
\lesssim\; & \lambda^{\gamma-1}\,[U_\gamma;\boldsymbol{X}]_{B_\lambda(z)}+\sum_{\rho\in\mathcal{T}_{-1,\gamma-2}}\|u'\|_{B_\lambda(z)}^{|\mathfrak{n}(\rho)|}\,\lambda^{|\mathcal{I}'(\rho)|}\,[\Gamma\,\mathcal{I}(\rho);\boldsymbol{X}].
\end{aligned}
\tag{3.70}
$$

Combining (3.69) with (3.57) then implies by (3.54) that

$$
\begin{aligned}
& |\Lambda_\lambda[H^{\mathbf{1}}_{\emptyset}-(\partial_{e_1}u_\lambda)(z)\,\tilde{H}^{\mathbf{1}},\Pi\Xi](z)|\\
\lesssim\; & \lambda^{2+|\Xi|}\,\lambda_t^{-1}\,\|u\|_{P_t}\,[\Pi;\Xi]\,\lambda^{\gamma-1}\,\lambda_t^{-\gamma}\,\|u\|_{P_t}\\
& +\sum_{\rho\in\mathcal{T}_{\gamma-2}}\lambda^{2+|\Xi|}\,\lambda_t^{-1}\,\|u\|_{P_t}\,[\Pi;\Xi]\,(\lambda_t^{-1}\,\|u\|_{P_t})^{|\mathfrak{n}(\rho)|}\,\lambda^{|\mathcal{I}'(\rho)|}\,|||(\Pi,\Gamma)|||^{\mathfrak{C}(\rho)}\\
\lesssim\; & \lambda^{\alpha-2}\,|||(\Pi,\Gamma)|||\,\lambda^{\gamma+1}\,\lambda_t^{-(\gamma+1)}\,\|u\|_{P_t}^2\\
& +\lambda^{\alpha-2}\,|||(\Pi,\Gamma)|||\sum_{\rho\in\mathcal{T}_{\gamma-2}}\lambda^{1+|\mathcal{I}(\rho)|}\,|||(\Pi,\Gamma)|||^{\mathfrak{C}(\rho)}\,\lambda_t^{-(1+|\mathfrak{n}(\rho)|)}\,\|u\|_{P_t}^{1+|\mathfrak{n}(\rho)|}\\
\sim\; & \lambda^{\alpha-2}\,|||(\Pi,\Gamma)|||\left(\left(\lambda\,\lambda_t^{-1}\,\|u\|_{P_t}^{\frac{2}{\gamma+1}}\right)^{\gamma+1}+\sum_{\rho\in\mathcal{T}_{\gamma-2}}\left(\lambda\,\lambda_t^{-1}\,\|u\|_{P_t}^{\hat{\mathcal{E}}^{\mathbf{1}}_\rho(\gamma)}\right)^{1+|\mathcal{I}(\rho)|}\right),
\end{aligned}
\tag{3.71}
$$

where

$$
\hat{\mathcal{E}}^{\mathbf{1}}_\rho(\gamma):=\frac{\alpha\,(1+|\mathfrak{n}(\rho)|)-(\gamma-1)\,(|\mathcal{I}(\rho)|-|\mathfrak{n}(\rho)|)}{\alpha\,(1+|\mathcal{I}(\rho)|)}.
$$

We can now conclude that for the term coming from $\tau=\mathbf{1}$ without the transport term we have

$$
\begin{aligned}
& |\Lambda_\lambda[H^{\mathbf{1}}-(\partial_{e_1}u_\lambda)(z)\,\tilde{H}^{\mathbf{1}},\Pi\Xi](z)|\\
\lesssim\; & \lambda^{\alpha-2}|||(\Pi,\Gamma)|||\left(\left(\lambda\lambda_t^{-1}\|u\|_{P_t}^{\frac{1}{\gamma}}\right)^{\gamma}+\sum_{\tau\in\mathfrak{a}(\partial A_{\mathbf{1}})}\left(\lambda\lambda_t^{-1}\|u\|_{P_t}^{\mathcal{E}^{\mathbf{1}}_\tau(\gamma)}\right)^{|\tau|}+\sum_{\tau\in\mathfrak{a}(A_{\mathbf{1}})\setminus\{\mathbf{1}\}}\left(\lambda\lambda_t^{-1}\|u\|_{P_t}^{\tilde{\mathcal{E}}^{\mathbf{1}}_\tau(\gamma)}\right)^{|\tau|}\right)\\
& +\lambda^{\alpha-2}\,|||(\Pi,\Gamma)|||\left(\left(\lambda\,\lambda_t^{-1}\,\|u\|_{P_t}^{\frac{2}{\gamma+1}}\right)^{\gamma+1}+\sum_{\rho\in\mathcal{T}_{\gamma-2}}\left(\lambda\,\lambda_t^{-1}\,\|u\|_{P_t}^{\hat{\mathcal{E}}^{\mathbf{1}}_\rho(\gamma)}\right)^{1+|\mathcal{I}(\rho)|}\right).
\end{aligned}
\tag{3.72}
$$

For the term coming from $\tau = \boldsymbol{X}$ we need to estimate (3.59). We have by (2.34), (2.41), and (2.42) and Corollary 3.3 for $\beta \in \partial A_{\boldsymbol{X}}$ and $\{\tau_i\}_{i=1}^n \subset \mathcal{T}$ such that $\sum_{i=1}^n |\mathcal{I}(\tau_i)| < \gamma - 1$ for $n \neq 0$ that

$$\begin{aligned} |u'(y) H_1^{\boldsymbol{X}}(x,y)| &\lesssim d(x,y)^{\gamma} \|u'\|_{B_\lambda(z)} [U_\gamma; \mathbf{1}]_{B_\lambda(z)} \\ &\lesssim d(x,y)^{\gamma} \lambda_t^{-(1+\gamma)} \|u\|_{P_t}^2, \end{aligned}$$

$$|u'(y) H_\beta^{\boldsymbol{X}}(x,y)| \lesssim d(x,y)^{|a(\beta)|} \lambda_t^{-1} \|u\|_{P_t} |||(\Pi,\Gamma)|||^{\mathfrak{C}(a(\beta))},$$

$$\begin{aligned} |H_{(\tau_i)_{i=1}^n}^{\boldsymbol{X}}(x,y)| &\lesssim d(x,y)^{\gamma-1} [U_{\gamma - \sum_{i=1}^n |\mathcal{I}(\tau_i)|}; \boldsymbol{X}]_{B_\lambda(z)} |||(\Pi,\Gamma)|||^{\mathfrak{C}(\prod_{i=1}^n \mathcal{I}(\tau_i))} \\ &\lesssim d(x,y)^{\gamma-1} \lambda_t^{-(\gamma - \sum_{i=1}^n |\mathcal{I}(\tau_i)|)} \|u\|_{P_t} |||(\Pi,\Gamma)|||^{\mathfrak{C}(\prod_{i=1}^n \mathcal{I}(\tau_i))}. \end{aligned}$$

Similarly we have for $\rho \in \mathcal{T}_{\gamma-3,\gamma-2}$

$$\begin{aligned} \left| u'(y) \frac{\Upsilon_x[\rho]}{\rho!} \gamma_{yx}(\mathcal{I}(\rho)) \tilde{H}_2^{\boldsymbol{X}}(x,y) \right| &\lesssim \|u'\|_{B_\lambda(z)}^{1+|\mathfrak{n}(\rho)|} d(x,y)^{|\mathcal{I}(\rho)|} [\Gamma \mathcal{I}(\rho); \mathbf{1}] \\ &\lesssim d(x,y)^{|\mathcal{I}(\rho)|} (\lambda_t^{-1} \|u\|_{P_t})^{1+|\mathfrak{n}(\rho)|} |||(\Pi,\Gamma)|||^{\mathfrak{C}(\rho)}, \end{aligned}$$

and by Lemma 3.24 and (3.60) and (3.69) we have

$$\begin{aligned} &|(u'(z) - (\partial_{e_1} u_\lambda)(z) + u'(y) - u'(z)) \hat{H}^{\boldsymbol{X}}(x,y)| \\ \lesssim &\left( \left| \langle U_\gamma(z,\cdot) - u'(z)(\cdot - z)_1, \partial_{e_1} \varphi_z^\lambda \rangle - \sum_{\rho \in \mathcal{T}_{\gamma-2}} \frac{\Upsilon_z[\rho]}{\rho!} \langle \Pi_z \mathcal{I}'(\rho), \varphi_z^\lambda \rangle \right| + |u'(y) - u'(z)| \right) \lambda_t^{-1} \|u\|_{P_t} d(x,y) \\ \lesssim &\left( \lambda^{\gamma-1} \lambda_t^{-\gamma} \|u\|_{P_t} + \sum_{\rho \in \mathcal{T}_{\gamma-2}} (\lambda_t^{-1} \|u\|_{P_t})^{|\mathfrak{n}(\rho)|} \lambda^{|\mathcal{I}'(\rho)|} |||(\Pi,\Gamma)|||^{\mathfrak{C}(\rho)} \right) \lambda_t^{-1} \|u\|_{P_t} d(x,y) \\ \lesssim &\left( \lambda^{\gamma-1} \lambda_t^{-(\gamma+1)} \|u\|_{P_t}^2 + \sum_{\rho \in \mathcal{T}_{\gamma-2}} (\lambda_t^{-1} \|u\|_{P_t})^{1+|\mathfrak{n}(\rho)|} \lambda^{|\mathcal{I}'(\rho)|} |||(\Pi,\Gamma)|||^{\mathfrak{C}(\rho)} \right) d(x,y) \end{aligned}$$

All these bounds on the germs imply by (3.54) (all the sums of regularities against $|\boldsymbol{X}\,\Xi|$ are positive)

$$\begin{aligned} &|\Lambda_\lambda[H^{\boldsymbol{X}} - (\partial_{e_1} u_\lambda)(z) \hat{H}^{\boldsymbol{X}}, \Pi(\boldsymbol{X}\,\Xi)](z)| \\ \lesssim\; &\lambda^{\gamma + |\boldsymbol{X}\Xi|} |||(\Pi,\Gamma)||| \lambda_t^{-(1+\gamma)} \|u\|_{P_t}^2 + \sum_{\beta \in \partial A_{\boldsymbol{X}}} \lambda^{|a(\beta)| + |\boldsymbol{X}\Xi|} \lambda_t^{-1} \|u\|_{P_t} |||(\Pi,\Gamma)|||^{\mathfrak{C}(a(\beta))+1} \\ &+ \sum_{\beta \in A_{\boldsymbol{X}}} \lambda^{\gamma - 1 + |\boldsymbol{X}\Xi|} \lambda_t^{-(\gamma - |a(\beta)|)} \|u\|_{P_t} |||(\Pi,\Gamma)|||^{\mathfrak{C}(a(\beta))+1} \\ &+ \sum_{\rho \in \mathcal{T}_{\gamma-3,\gamma-2}} \lambda^{|\mathcal{I}(\rho)| + |\boldsymbol{X}\Xi|} (\lambda_t^{-1} \|u\|_{P_t})^{1+|\mathfrak{n}(\rho)|} |||(\Pi,\Gamma)|||^{\mathfrak{C}(\rho)+1} \\ &+ \lambda^{1+|\boldsymbol{X}\Xi|} |||(\Pi,\Gamma)||| \left( \lambda^{\gamma-1} \lambda_t^{-(\gamma+1)} \|u\|_{P_t}^2 + \sum_{\rho \in \mathcal{T}_{\gamma-2}} (\lambda_t^{-1} \|u\|_{P_t})^{1+|\mathfrak{n}(\rho)|} \lambda^{|\mathcal{I}'(\rho)|} |||(\Pi,\Gamma)|||^{\mathfrak{C}(\rho)} \right) \\ \lesssim\; &\lambda^{\alpha-2} |||(\Pi,\Gamma)||| \left( \left( \lambda \lambda_t^{-1} \|u\|_{P_t}^{\frac{2}{\gamma+1}} \right)^{\gamma+1} + \sum_{\tau \in a(\partial A_{\boldsymbol{X}})} \left( \lambda (\lambda_t^{-1} \|u\|_{P_t})^{\frac{1}{|\tau|+1}} (1 + |||(\Pi,\Gamma)|||)^{\frac{\mathfrak{C}(\tau)}{|\tau|+1}} \right)^{|\tau|+1} \right) \\ &+ \lambda^{\alpha-2} |||(\Pi,\Gamma)||| \sum_{\tau \in a(A_{\boldsymbol{X}})} \left( \lambda \lambda_t^{-\frac{(\gamma - |\tau|)}{\gamma}} \|u\|_{P_t}^{\frac{1}{\gamma}} (1 + |||(\Pi,\Gamma)|||)^{\frac{\mathfrak{C}(\tau)}{\gamma}} \right)^{\gamma} \\ &+ \lambda^{\alpha-2} |||(\Pi,\Gamma)||| \sum_{\rho \in \mathcal{T}_{\gamma-3,\gamma-2}} \left( \lambda (\lambda_t^{-1} \|u\|_{P_t})^{\frac{1+|\mathfrak{n}(\rho)|}{|\mathcal{I}(\rho)|+1}} (1 + |||(\Pi,\Gamma)|||)^{\frac{\mathfrak{C}(\rho)}{|\mathcal{I}(\rho)|+1}} \right)^{|\mathcal{I}(\rho)|+1}. \end{aligned}$$

For $\tau \in \mathfrak{a}(\partial A_{\boldsymbol{X}})$ we have that $|\mathfrak{n}(\tau)| = 0$ (see Proof of Lemma 2.5) and therefore

$$\begin{aligned}\lambda\,(\lambda_t^{-1}\,\|u\|_{P_t})^{\frac{1}{|\tau|+1}}\,(1+|\!|\!|(\Pi,\Gamma)|\!|\!|)^{\frac{\mathfrak{e}(\tau)}{|\tau|+1}} &\sim \lambda\,\|u\|_{P_t}^{\frac{1}{\alpha}(\gamma-1+\alpha)\frac{1}{|\tau|+1}}\,(1+|\!|\!|(\Pi,\Gamma)|\!|\!|)^{\frac{1}{\alpha}\frac{1}{|\tau|+1}+\frac{\mathfrak{e}(\tau)}{|\tau|+1}}\\ &= \lambda\,\|u\|_{P_t}^{\frac{1}{\alpha}(\gamma-1+\alpha)\frac{1}{|\tau|+1}}\,(1+|\!|\!|(\Pi,\Gamma)|\!|\!|)^{\frac{1}{\alpha(|\tau|+1)}+\frac{|\tau|}{\alpha(|\tau|+1)}}\\ &= \lambda\,\|u\|_{P_t}^{\mathcal{E}_\tau^{\boldsymbol{X}}(\gamma)+\frac{(\gamma-1)}{\alpha}}\,(1+|\!|\!|(\Pi,\Gamma)|\!|\!|)^{\frac{1}{\alpha}}\\ &= \lambda\,\lambda_t^{-1}\,\|u\|_{P_t}^{\mathcal{E}_\tau^{\boldsymbol{X}}(\gamma)},\end{aligned}$$

where

$$\mathcal{E}_\tau^{\boldsymbol{X}}(\gamma) := \frac{1}{\alpha}(\gamma-1+\alpha)\,\frac{1}{|\tau|+1} - \frac{(\gamma-1)}{\alpha} = \frac{\alpha-(\gamma-1)\,|\tau|}{\alpha\,(|\tau|+1)}.$$

For $\tau \in \mathfrak{a}(A_{\boldsymbol{X}})$ we have

$$\begin{aligned}\lambda\,\lambda_t^{-\frac{(\gamma-|\tau|)}{\gamma}}\,\|u\|_{P_t}^{\frac{1}{\gamma}}\,(1+|\!|\!|(\Pi,\Gamma)|\!|\!|)^{\frac{\mathfrak{e}(\tau)}{\gamma}} &\sim \lambda\,\|u\|_{P_t}^{\frac{(\gamma-1)}{\alpha}\frac{(\gamma-|\tau|)}{\gamma}+\frac{1}{\gamma}}\,(1+|\!|\!|(\Pi,\Gamma)|\!|\!|)^{\frac{1}{\alpha}\frac{\gamma-|\tau|}{\gamma}+\frac{\mathfrak{e}(\tau)}{\gamma}}\\ &= \lambda\,\|u\|_{P_t}^{\frac{(\gamma-1)(\gamma-|\tau|)+\alpha}{\alpha\gamma}}\,(1+|\!|\!|(\Pi,\Gamma)|\!|\!|)^{\frac{\gamma-|\tau|}{\alpha\gamma}+\frac{|\tau|}{\alpha\gamma}}\\ &= \lambda\,\|u\|_{P_t}^{\tilde{\mathcal{E}}_\tau^{\boldsymbol{X}}(\gamma)+\frac{(\gamma-1)}{\alpha}}\,(1+|\!|\!|(\Pi,\Gamma)|\!|\!|)^{\frac{1}{\alpha}}\\ &= \lambda\,\lambda_t^{-1}\,\|u\|_{P_t}^{\tilde{\mathcal{E}}_\tau^{\boldsymbol{X}}(\gamma)},\end{aligned}$$

where

$$\tilde{\mathcal{E}}_\tau^{\boldsymbol{X}}(\gamma) := \frac{(\gamma-1)\,(\gamma-|\tau|)+\alpha}{\alpha\,\gamma} - \frac{(\gamma-1)}{\alpha} = \frac{\alpha-(\gamma-1)\,|\tau|}{\alpha\,\gamma}.$$

For $\rho \in \mathcal{T}_{\gamma-3,\gamma-2}$ we have

$$\begin{aligned}&\lambda(\lambda_t^{-1}\|u\|_{P_t})^{\frac{1+|\mathfrak{n}(\rho)|}{|\mathcal{I}(\rho)|+1}}(1+|\!|\!|(\Pi,\Gamma)|\!|\!|)^{\frac{\mathfrak{e}(\rho)}{|\mathcal{I}(\rho)|+1}}\\ \sim\;& \lambda\|u\|_{P_t}^{\frac{1}{\alpha}(\gamma-1+\alpha)\frac{1+|\mathfrak{n}(\rho)|}{|\mathcal{I}(\rho)|+1}}(1+|\!|\!|(\Pi,\Gamma)|\!|\!|)^{\frac{1}{\alpha}\frac{1+|\mathfrak{n}(\rho)|}{|\mathcal{I}(\rho)|+1}+\frac{\mathfrak{e}(\rho)}{|\mathcal{I}(\rho)|+1}}\\ =\;& \lambda\|u\|_{P_t}^{\frac{1}{\alpha}(\gamma-1+\alpha)\frac{1+|\mathfrak{n}(\rho)|}{|\mathcal{I}(\rho)|+1}}(1+|\!|\!|(\Pi,\Gamma)|\!|\!|)^{\frac{1+|\mathfrak{n}(\rho)|}{\alpha(|\mathcal{I}(\rho)|+1)}+\frac{|\mathcal{I}(\rho)|-|\mathfrak{n}(\rho)|}{\alpha(|\mathcal{I}(\rho)|+1)}}\\ =\;& \lambda\|u\|_{P_t}^{\mathcal{E}_\rho^{\boldsymbol{X}}(\gamma)+\frac{(\gamma-1)}{\alpha}}(1+|\!|\!|(\Pi,\Gamma)|\!|\!|)^{\frac{1}{\alpha}}\\ =\;& \lambda\lambda_t^{-1}\|u\|_{P_t}^{\mathcal{E}_\rho^{\boldsymbol{X}}(\gamma)},\end{aligned}$$

where

$$\mathcal{E}_\rho^{\boldsymbol{X}}(\gamma) := \frac{1}{\alpha}(\gamma-1+\alpha)\,\frac{1+|\mathfrak{n}(\rho)|}{|\mathcal{I}(\rho)|+1} - \frac{(\gamma-1)}{\alpha} = \frac{\alpha\,(1+|\mathfrak{n}(\rho)|)-(\gamma-1)\,(|\mathcal{I}(\rho)|-|\mathfrak{n}(\rho)|)}{\alpha\,(|\mathcal{I}(\rho)|+1)}.$$

We conclude that

$$\begin{aligned}&|\Lambda_\lambda[H^{\boldsymbol{X}}-(\partial_{e_1}u_\lambda)(z)\hat{H}^{\boldsymbol{X}},\Pi(\boldsymbol{X}\Xi)](z)| \\ \lesssim\;& \lambda^{\alpha-2}|\!|\!|(\Pi,\Gamma)|\!|\!|\Big(\lambda\lambda_t^{-1}\|u\|_{P_t}^{\frac{2}{\gamma+1}}\Big)^{\gamma+1}+\lambda^{\alpha-2}|\!|\!|(\Pi,\Gamma)|\!|\!|\sum_{\tau\in\mathfrak{a}(\partial A_{\boldsymbol{X}})}\Big(\lambda\lambda_t^{-1}\|u\|_{P_t}^{\mathcal{E}_\tau^{\boldsymbol{X}}(\gamma)}\Big)^{|\tau|+1}\\ &+\lambda^{\alpha-2}|\!|\!|(\Pi,\Gamma)|\!|\!|\sum_{\tau\in\mathfrak{a}(A_{\boldsymbol{X}})}\Big(\lambda\lambda_t^{-1}\|u\|_{P_t}^{\tilde{\mathcal{E}}_\tau^{\boldsymbol{X}}(\gamma)}\Big)^{\gamma}+\lambda^{\alpha-2}|\!|\!|(\Pi,\Gamma)|\!|\!|\sum_{\rho\in\mathcal{T}_{\gamma-3,\gamma-2}}\Big(\lambda\lambda_t^{-1}\|u\|_{P_t}^{\mathcal{E}_\rho^{\boldsymbol{X}}(\gamma)}\Big)^{|\mathcal{I}(\rho)|+1}.\end{aligned} \tag{3.73}$$

We want to choose $\lambda \in (0, \lambda_t)$ such that (3.66), (3.72), and (3.73) are bounded by $\lambda^{\alpha-2} |||(\Pi, \Gamma)|||$, and therefore, since $\|u\|_{P_t} \geqslant 1$ we need to impose the condition $\lambda \lesssim \lambda_t \|u\|_{P_t}^{-\mathcal{E}}$, where

$$\begin{aligned}\mathcal{E} &= \max_{\tau \in \mathcal{T}_\gamma^+ \setminus \{\mathbf{1}, \boldsymbol{X}\}} \mathcal{E}^\tau(\gamma) \vee \frac{1}{\gamma} \vee \max_{\tau \in a(\partial A_{\mathbf{1}})} \mathcal{E}_\tau^{\mathbf{1}}(\gamma) \vee \max_{\tau \in a(A_{\mathbf{1}}) \setminus \{\mathbf{1}\}} \tilde{\mathcal{E}}_\tau^{\mathbf{1}}(\gamma) \vee \frac{2}{\gamma+1} \vee \max_{\rho \in \mathcal{T}_{\gamma-2}} \hat{\mathcal{E}}_\rho^{\mathbf{1}}(\gamma) \\ &\quad \vee \max_{\tau \in a(\partial A_{\boldsymbol{X}})} \mathcal{E}_\tau^{\boldsymbol{X}}(\gamma) \vee \max_{\tau \in a(A_{\boldsymbol{X}})} \tilde{\mathcal{E}}_\tau^{\boldsymbol{X}}(\gamma) \vee \max_{\rho \in \mathcal{T}_{\gamma-3,\gamma-2}} \mathcal{E}_\rho^{\boldsymbol{X}}(\gamma).\end{aligned}$$

We claim that $\mathcal{E} = \frac{2}{\gamma+1}$. Clearly, since $\gamma > 1$ we have that $\frac{1}{\gamma} \leqslant \frac{2}{\gamma+1}$. Moreover, one can easily see from the explicit form of each exponent that they are decreasing in the homogeneity of the corresponding tree and therefore

$$\begin{aligned}
\max_{\tau \in \mathcal{T}_\gamma^+ \setminus \{\mathbf{1}, \boldsymbol{X}\}} \mathcal{E}^\tau(\gamma) &\leqslant \mathcal{E}^{\mathcal{I}(\Xi)}(\gamma) \vee \mathcal{E}^{\mathcal{I}(\boldsymbol{X}\Xi)}(\gamma) = \frac{2-\gamma}{\gamma} \vee \frac{\gamma - \alpha(\gamma-1)}{\gamma^2} \leqslant \frac{1}{\gamma}, \\
\max_{\tau \in a(\partial A_{\mathbf{1}})} \mathcal{E}_\tau^{\mathbf{1}}(\gamma) &\leqslant \mathcal{E}_{\mathcal{I}(\Xi)}^{\mathbf{1}}(\gamma) \vee \mathcal{E}_{\mathcal{I}(\boldsymbol{X}\Xi)}^{\mathbf{1}}(\gamma) = \frac{-(\gamma-1)\alpha}{\alpha^2} \vee \frac{\alpha - (\gamma-1)\alpha}{\alpha(\alpha+1)} = \frac{2-\gamma}{\alpha+1}, \\
\max_{\tau \in a(A_{\mathbf{1}}) \setminus \{\mathbf{1}\}} \tilde{\mathcal{E}}_\tau^{\mathbf{1}}(\gamma) &\leqslant \tilde{\mathcal{E}}_{\mathcal{I}(\Xi)}^{\mathbf{1}}(\gamma) \vee \tilde{\mathcal{E}}_{\mathcal{I}(\boldsymbol{X}\Xi)}^{\mathbf{1}}(\gamma) = 1 - \frac{(\gamma-1+\alpha)\alpha}{\alpha\gamma} = \frac{1-\alpha}{\gamma}, \\
\max_{\rho \in \mathcal{T}_{\gamma-2}} \hat{\mathcal{E}}_\rho^{\mathbf{1}}(\gamma) &\leqslant \hat{\mathcal{E}}_{\Xi}^{\mathbf{1}}(\gamma) \vee \hat{\mathcal{E}}_{\boldsymbol{X}\Xi}^{\mathbf{1}}(\gamma) = \frac{\alpha - (\gamma-1)\alpha}{\alpha(1+\alpha)} \vee \frac{2\alpha - (\gamma-1)\alpha}{\alpha(2+\alpha)} = \frac{3-\gamma}{\alpha+2}, \\
\max_{\tau \in a(\partial A_{\boldsymbol{X}})} \mathcal{E}_\tau^{\boldsymbol{X}}(\gamma) &\leqslant \mathcal{E}_{\mathcal{I}(\Xi)}^{\boldsymbol{X}}(\gamma) = \frac{\alpha - (\gamma-1)\alpha}{\alpha(\alpha+1)} = \frac{2-\gamma}{\alpha+1}, \\
\max_{\tau \in a(A_{\boldsymbol{X}})} \tilde{\mathcal{E}}_\tau^{\boldsymbol{X}}(\gamma) &\leqslant \tilde{\mathcal{E}}_{\mathbf{1}}^{\boldsymbol{X}}(\gamma) \vee \tilde{\mathcal{E}}_{\mathcal{I}(\Xi)}^{\boldsymbol{X}}(\gamma) = \frac{1}{\gamma} \vee \frac{\alpha - (\gamma-1)\alpha}{\alpha\gamma} = \frac{1}{\gamma} \vee \frac{2-\gamma}{\gamma} = \frac{1}{\gamma}, \\
\max_{\rho \in \mathcal{T}_{\gamma-3,\gamma-2}} \mathcal{E}_\rho^{\boldsymbol{X}}(\gamma) &\leqslant \mathcal{E}_{\Xi}^{\boldsymbol{X}}(\gamma) \vee \mathcal{E}_{\boldsymbol{X}\Xi}^{\boldsymbol{X}}(\gamma) = \frac{\alpha - (\gamma-1)\alpha}{\alpha(\alpha+1)} \vee \frac{2\alpha - (\gamma-1)\alpha}{\alpha(\alpha+2)} = \frac{3-\gamma}{\alpha+2},
\end{aligned}$$

where we used that trees in $a(A_{\boldsymbol{X}} \cup \partial A_{\boldsymbol{X}})$ have no polynomial decorations, and from here one can directly verify that all these terms are dominated by $\frac{2}{\gamma+1}$ if $\gamma \in (2-\alpha, 2)$ and $\alpha \in (0,1)$. □

## 4. Differentiability of the Flow

For the differentiability statements in this section, we assume $m \geqslant 2$. Set

$$b(v) := \beta v - v|v|^{m-1}, \qquad b'(v) := \beta - m|v|^{m-1}. \tag{4.1}$$

The derivative $b'$ is locally Lipschitz on bounded sets.

We introduce the definitions for *singular* modelled distributions (see [Hai14, Section 6]).

**Definition 4.1**. *For fixed $0 < s < T$ and any $\gamma > 0$, $\eta \in \mathbb{R}$ the space $\mathcal{D}_{(s,T]}^{\gamma,\eta}$ consists of all functions $f: [s,T] \times \mathbb{R} \to T$ such that*

$$|||f|||_{\gamma,\eta;(s,T]} := \|f\|_{\gamma,\eta;(s,T]} + \sup_{(x,y) \in D_{(s,T]}} \sup_{\tau \in \mathcal{T}} \frac{|\langle \tau, f(x) - \Gamma_{yx} f(y) \rangle|}{d(x,y)^{\gamma - |\tau|} (|x_0 - s| \wedge |y_0 - s|)^{\frac{1}{2}(\eta - \gamma)}} < +\infty, \tag{4.2}$$

*where* $D_{(s,T]} = \left\{ (x,y) \in ((s,T] \times \mathbb{R})^2 : x \neq y, d(x,y) \leqslant (|x_0 - s| \wedge |y_0 - s|)^{\frac{1}{2}} \right\}$

$$\|f\|_{\gamma,\eta;(s,T]} := \sup_{x \in (s,T] \times \mathbb{R}} \sup_{\tau \in \mathcal{T}} \frac{|\langle \tau, f(x) \rangle|}{|x_0 - s|^{\frac{1}{2}(\eta - |\tau|) \wedge 0}}. \tag{4.3}$$

*In particular, when $s = 0$ we use the notations $\mathcal{D}_T^{\gamma,\eta}$, $|||f|||_{\gamma,\eta;T}$ and $\|f\|_{\gamma,\eta;T}$.*

For smooth models, the local-in-time mild solutions from [Hai14, Proposition 7.11] satisfy the weak formulation used in Theorem 1.1 and Theorem 3.13, and hence are global in time. The same holds for the Itô model, since it was constructed in [HP15] as the limit of smooth models.

**DEFINITION 4.2.** *Fix $\alpha\in(0,1),\beta\in\mathbb{R},m>\frac{2-\alpha}{\alpha}\varepsilon_\alpha$, $\gamma\in(2-\alpha,2)$ as in Theorem 1.1. Assume that $\sigma\in C_b^{k+1}$ where $k=\lceil\gamma/\alpha\rceil$. Let $(\Pi,\Gamma)$ be a 1-periodic admissible model, as in [Hai14, Definition 5.9], on the regularity structure described in Section 2.1. We denote by $\Phi^{(\Pi,\Gamma)}=\Phi$ the flow induced by the modelled solutions to (1.3), i.e., for any $0\leqslant s<T$ we have that $\Phi(s,\,\cdot\,;\,\cdot\,)\colon[s,T]\times L^\infty(\mathbb{T})\to\mathcal{D}^{\gamma,0}_{(s,T]}$ is a solution to the fixed point problem*

$$\Phi(s,\,\cdot\,;h)=(\mathcal{K}_{\gamma-2}+R_\gamma\mathcal{R})\boldsymbol{R}^+(\tilde{b}(\Phi(s,\,\cdot\,;h))+\hat{\sigma}(\Phi(s,\,\cdot\,;h))\,\Xi)+\mathscr{T}_\gamma Ph$$

*where $\mathscr{T}_\gamma$ is the lift of a function $u$ to the modelled distribution defined by its Taylor expansion of order $\lfloor\gamma\rfloor$, i.e.*

$$(\mathscr{T}_\gamma u)(x)=\sum_{|k|<\gamma}(D^k u)(z)\,\frac{\boldsymbol{X}^k}{k!},$$

*$Ph$ is the harmonic extension of $h$, $\mathcal{P}=(\mathcal{K}_{\gamma-2}+R_\gamma\mathcal{R})\boldsymbol{R}^+$ is the lift of the integration operator to modelled distributions as in [Hai14, Section 7.3], and $b$ is the drift defined by (4.1).*

*Then $\varphi:=\mathcal{R}\Phi$ defines a flow*

$$\varphi(s,t;\,\cdot\,)\colon L^\infty(\mathbb{T})\to L^\infty(\mathbb{T})$$

*for $0\leqslant s\leqslant t\leqslant T$, such that $\varphi(s,s;\,\cdot\,)=\mathrm{Id}_{L^\infty(\mathbb{T})}$ and which satisfies the flow property*

$$\varphi(r,t;\varphi(s,r;u_0))=\varphi(s,t;u_0),\qquad\forall 0\leqslant s\leqslant r\leqslant t,\quad u_0\in L^\infty(\mathbb{T}).$$

For simplicity we assume that $s=0$ and remove it from the dependence on $\Phi$ and $\varphi$.

**DEFINITION 4.3.** *If $\sigma\in C_b^k(\mathbb{R})$ for $k=\lceil\gamma/\alpha\rceil+2$ and fix $u_0\in L^\infty(\mathbb{T})$, $T>0$. For any $h\in L^\infty(\mathbb{T})$ let $W_{u_0}[h]\in\mathcal{D}_T^{\gamma,0}$ be the solution to the linear fixed-point problem*

$$W_{u_0}[h]\;=\;(\mathcal{K}_{\gamma-2}+R_\gamma\mathcal{R})\boldsymbol{R}^+((\tilde{b}(\Phi(s,\,\cdot\,;h))+\hat{\sigma}^{(1)}(\Phi(\,\cdot\,;u_0))\,\Xi)\,W_{u_0}[h])+Gh. \tag{4.4}$$

**THEOREM 4.4.** *Under the assumptions of Definitions 4.2 and 4.3, there exists $C(u_0)>0$ such that the estimate*

$$|||\Phi(\,\cdot\,;u_0+h)-\Phi(\,\cdot\,;u_0)-W_{u_0}[h]|||_{\gamma,0;T}\lesssim\|h\|^2_{L^\infty(\mathbb{T})} \tag{4.5}$$

*holds uniformly for $\|h\|_{L^\infty(\mathbb{T})}\leqslant 1$. In particular, the flow $\varphi$ is differentiable, with respect to the initial condition, in any direction $h\in L^\infty(\mathbb{T})$. Moreover, this directional derivative $(D_{u_0}\Phi)[h]$ satisfies the identity $(D_{u_0}\varphi)[h]=\mathcal{R}(W_{u_0}[h])$ and solves the linearisation of the original equation, i.e.,*

$$\begin{aligned}(\partial_t-\partial_x^2)\,\mathcal{R}(W_{u_0}[h])\;&=\;D_{u_0}((\partial_t-\partial_x^2)\,\varphi(\,\cdot\,;u_0))[h]\\&=\;(\beta-m\,|u|^{m-1})\,w+\mathcal{R}(\hat{\sigma}^{(1)}(\Phi(\,\cdot\,;u_0))\,W_{u_0}[h]\,\Xi)\end{aligned} \tag{4.6}$$

*with initial condition $\mathcal{R}(W_{u_0}[h])(0)=h$.*

We postpone the proof of Theorem 4.4 to Section 4.1.

## 4.1. Proof of Theorem 4.4

Since $b'$ is locally Lipschitz on bounded sets, for every $R>0$ there exists $C_R$ such that

$$|b\,(v+r)-b(v)-b'(v)\,r|\leqslant C_R|r|^2,\qquad|v|,|v+r|\leqslant R.$$

Thus the difference of the two drift terms is its linearisation plus a quadratic remainder, which is handled by the deterministic fixed-point estimates. We therefore focus below on the singular part. The damping remains important for obtaining the uniform bound on $\|U_h\|_{\gamma,0;T}$ for $\|h\|_{L^\infty(\mathbb{T})} \leqslant 1$, using the a priori estimates from Theorem 1.1.

**LEMMA 4.5**. *We have the estimate*

$$|||\Phi(\,\cdot\,;u_0+h)-\Phi(\,\cdot\,;u_0)|||_{\gamma,0;T} \lesssim C^{\lceil T/T_1 \rceil}\,\|h\|_{L^\infty(\mathbb{T})},$$

*uniformly on $h \in L^\infty(\mathbb{T})$, where $C>1$ and $T_1^{-1}$ depends polynomially on $|||\Phi(\,\cdot\,;u_0+h)|||_{\gamma,0;T}$.*

**Proof.** We have that $\Phi(\,\cdot\,;u_0+h)-\Phi(\,\cdot\,;u_0)$ satisfies in $\mathcal{D}_T^{\gamma,0}$

$$\begin{aligned} & \Phi(\,\cdot\,;u_0+h)-\Phi(\,\cdot\,;u_0) \\ = \ & (\mathcal{K}_{\gamma-2}+R_\gamma\mathcal{R})\boldsymbol{R}^+(\hat{b}(\Phi(\,\cdot\,;u_0+h))-\hat{b}(\Phi(\,\cdot\,;u_0))+(\hat{\sigma}(\Phi(\,\cdot\,;u_0+h))-\hat{\sigma}(\Phi(\,\cdot\,;u_0)))\Xi)+Gh \end{aligned} \tag{4.7}$$

By [Hai14, Theorem 7.1] one has, for $\kappa=\alpha>0$ in the notation there, the bound

$$|||\mathcal{K}_{\gamma-2}\boldsymbol{R}^+((\hat{\sigma}(\Phi(\,\cdot\,;u_0+h))-\hat{\sigma}(\Phi(\,\cdot\,;u_0)))\Xi)|||_{\gamma,0;T} \lesssim T^{\alpha/2}|||(\hat{\sigma}(\Phi(\,\cdot\,;u_0+h))-\hat{\sigma}(\Phi(\,\cdot\,;u_0)))\Xi|||_{\gamma-2,|\Xi|;T},$$

since $\bar{\eta}=(\eta\wedge|\Xi|)+2-\kappa=(|\Xi|\wedge|\Xi|)+2-\alpha=|\Xi|+2-\alpha=0$. Analogously, by [Hai14, Lemma 7.3] one has the bound on the smooth remainder

$$|||R_\gamma\mathcal{R}\boldsymbol{R}^+((\hat{\sigma}(\Phi(\,\cdot\,;u_0+h))-\hat{\sigma}(\Phi(\,\cdot\,;u_0)))\Xi)|||_{\gamma,0;T} \lesssim T\,|||(\hat{\sigma}(\Phi(\,\cdot\,;u_0+h))-\hat{\sigma}(\Phi(\,\cdot\,;u_0)))\Xi|||_{\gamma-2,|\Xi|;T}.$$

Since $\Xi \in \mathcal{D}_T^{\tilde{\gamma},\tilde{\gamma}}$ for any $\tilde{\gamma}>0$ (see [Hai14, Proof of Lemma 9.1]) and $\Xi$ lies in a sector of regularity $|\Xi|$ (since $\Gamma\Xi=\Xi$), one can see that (a simplification of) the argument in [Hai14, Proposition 6.12] implies the estimate

$$|||(\hat{\sigma}(\Phi(\,\cdot\,;u_0+h))-\hat{\sigma}(\Phi(\,\cdot\,;u_0)))\,\Xi|||_{\gamma-2,|\Xi|;T} \lesssim |||\hat{\sigma}(\Phi(\,\cdot\,;u_0+h))-\hat{\sigma}(\Phi(\,\cdot\,;u_0))|||_{\gamma-\alpha,0;T},$$

where we used that $\hat{\sigma}(\Phi(\,\cdot\,;u_0+h))-\hat{\sigma}(\Phi(\,\cdot\,;u_0))$ lies in a function-like sector which implies that $\gamma-2=((\gamma-\alpha)+|\Xi|)\wedge(\tilde{\gamma}+0)$ and $|\Xi|=(0+|\Xi|)\wedge(\tilde{\gamma}+0)\wedge(0+\tilde{\gamma})$ when choosing $\tilde{\gamma}>0$ large enough. Moreover, one has the embedding $\mathcal{D}_T^{\gamma-\alpha,0} \hookrightarrow \mathcal{D}_T^{\gamma,0}$ with the bound

$$|||\hat{\sigma}(\Phi(\,\cdot\,;u_0+h))-\hat{\sigma}(\Phi(\,\cdot\,;u_0))|||_{\gamma-\alpha,0;T} \lesssim |||\hat{\sigma}(\Phi(\,\cdot\,;u_0+h))-\hat{\sigma}(\Phi(\,\cdot\,;u_0))|||_{\gamma,0;T}.$$

From [Hai14, Proposition 6.13] we know that the map $\hat{\sigma}\colon\mathcal{D}_T^{\gamma,0}(V)\to\mathcal{D}_T^{\gamma,0}(V)$ (where $V$ is a function-like sector) is locally Lipschitz, which implies the estimate

$$|||\hat{\sigma}(\Phi(\,\cdot\,;u_0+h))-\hat{\sigma}(\Phi(\,\cdot\,;u_0))|||_{\gamma,0;T} \lesssim C(|||\Phi(\,\cdot\,;u_0+h)|||_{\gamma,0;T})\,|||\Phi(\,\cdot\,;u_0+h)-\Phi(\,\cdot\,;u_0)|||_{\gamma,0;T}.$$

At last, a trivial modification of [Hai14, Lemma 7.5] using the contractivity of the heat semigroup in $L^\infty$ allows us to conclude the bound on the initial condition

$$|||Gh|||_{\gamma,0;T} \lesssim \|h\|_{L^\infty(\mathbb{T})}.$$

Combining this and using (4.7) one concludes that for any $T_1 \in (0,T]$

$$\begin{aligned} & |||\Phi(\,\cdot\,;u_0+h)-\Phi(\,\cdot\,;u_0)|||_{\gamma,0;T_1} \\ \leqslant \ & |||\mathcal{K}_{\gamma-2}\boldsymbol{R}^+((\hat{\sigma}(\Phi(\,\cdot\,;u_0+h))-\hat{\sigma}(\Phi(\,\cdot\,;u_0)))\,\Xi)|||_{\gamma,0;T_1}+|||R_\gamma\mathcal{R}\boldsymbol{R}^+((\hat{\sigma}(\Phi(\,\cdot\,;u_0+h))-\hat{\sigma}(\Phi(\,\cdot\,; \\ & u_0)))\,\Xi)|||_{\gamma,0;T_1}+|||Gh|||_{\gamma,0;T_1} \\ \lesssim \ & C(|||\Phi(\,\cdot\,;u_0)|||_{\gamma,0;T},|||\Phi(\,\cdot\,;u_0+h)|||_{\gamma,0;T})\,(T_1^{\alpha/2}+T_1)\,|||\Phi(\,\cdot\,;u_0+h)-\Phi(\,\cdot\,;u_0)|||_{\gamma,0;T_1}+\|h\|_{L^\infty(\mathbb{T})}, \end{aligned}$$

and therefore choosing $0<T_1\ll 1$ small enough one concludes that

$$|||\Phi(\,\cdot\,;u_0+h)-\Phi(\,\cdot\,;u_0)|||_{\gamma,0;T_1}\leqslant C\,\|h\|_{L^\infty(\mathbb{T})}.$$

Since $T_1^{-1}$ dependence on $h$ is only polynomially in $|||\Phi(\,\cdot\,;u_0+h)|||_{\gamma,0;T}$ (on the full interval $(0,T]$), which is finite by assumption, a bootstrapping argument allows us to conclude

$$|||\Phi(\,\cdot\,;u_0+h)-\Phi(\,\cdot\,;u_0)|||_{\gamma,0;T}\lesssim C^{\lceil T/T_1\rceil}\,\|h\|_{L^\infty(\mathbb{T})}.$$ □

**LEMMA 4.6**. *Given $U,\tilde{U}\in\mathcal{D}_T^{\gamma,0}$ and $\sigma\in C_b^k(\mathbb{R})$ for $k=\lceil\gamma/\alpha\rceil+2$, there exists $C>0$ depending polynomially on $\|U\|_{\gamma,0;T},\|U+\tilde{U}\|_{\gamma,0;T}$ such that*

$$|||\hat{\sigma}(U+\tilde{U})-\hat{\sigma}(U)-\hat{\sigma}^{(1)}(U)\,\tilde{U}|||_{\gamma,0;T}\leqslant C\,|||\tilde{U}|||_{\gamma,0;T}^2.$$

**Proof.** Let $U,\tilde{U}\in\mathcal{D}_T^{\gamma,0}$ and for $s\in[0,1]$ define $U_s:=U+s\,\tilde{U}\in\mathcal{D}_T^{\gamma,0}$. Denote by $u:=\langle\mathbf{1},U\rangle$ and let $\tilde{u}$ and $u_s$ be defined analogously. By definition we have that

$$\hat{\sigma}(U_s)=\sum_{k\geqslant 0}\frac{\sigma^{(k)}(u_s)}{k!}\,(U_s-u_s\,\mathbf{1})^k\in\mathcal{D}_T^{\gamma,0},$$

with the sum being finite. For every $z\in\mathbb{R}^{1+d}$, the map $s\mapsto\hat{\sigma}(U_s)(z)\in T$ is then differentiable and

$$\begin{aligned}
\frac{\mathrm{d}}{\mathrm{d}s}\,\hat{\sigma}(U_s) &= \sum_{k\geqslant 0}\frac{\sigma^{(k+1)}(u_s)\,\frac{\mathrm{d}}{\mathrm{d}s}u_s}{k!}\,(U_s-u_s\,\mathbf{1})^k+\sum_{k\geqslant 0}\frac{\sigma^{(k)}(u_s)}{k!}\,\frac{\mathrm{d}}{\mathrm{d}s}(U-u_s\,\mathbf{1})^k\\
&= \sum_{k\geqslant 0}\frac{\sigma^{(k+1)}(u_s)\,\tilde{u}}{k!}\,(U_s-u_s\,\mathbf{1})^k+\sum_{k\geqslant 1}\frac{\sigma^{(k)}(u_s)}{(k-1)!}\,(U-u_s\,\mathbf{1})^{k-1}\,(\tilde{U}-\tilde{u}\,\mathbf{1})\\
&= \left(\sum_{k\geqslant 0}\frac{\sigma^{(k+1)}(u_s)\,\tilde{u}}{k!}\,(U_s-u_s\,\mathbf{1})^k\right)(\tilde{u}\,\mathbf{1}+\tilde{U}-\tilde{u}\,\mathbf{1})\\
&= \hat{\sigma}^{(1)}(U_s)\,\tilde{U}.
\end{aligned}$$

Integrating over $s\in[0,1]$ and subtracting $\hat{\sigma}^{(1)}(U)\,\tilde{U}$ we obtain the pointwise identity

$$\hat{\sigma}(U+\tilde{U})-\hat{\sigma}(U)-\hat{\sigma}^{(1)}(U)\,\tilde{U}=\int_0^1(\hat{\sigma}^{(1)}(U_s)-\hat{\sigma}^{(1)}(U))\,\tilde{U}\,\mathrm{d}s.$$

Since $U_s=s\,(U+\tilde{U})+(1-s)\,U$, it is a convex combination of $U$ and $U+\tilde{U}$, then by the local Lipschitz continuity from [Hai14, Proposition 6.13] one obtains the estimate

$$|||\hat{\sigma}^{(1)}(U_s)-\hat{\sigma}(U)|||_{\gamma,0;T}\lesssim C\,|||U+\tilde{U}|||_{\gamma,0;T}=C\,s\,|||\tilde{U}|||_{\gamma,0;T}\leqslant C\,|||\tilde{U}|||_{\gamma,0;T},$$

for $C>0$ depending polynomially in $|||U|||_{\gamma,0;T},|||U+\tilde{U}|||_{\gamma,0;T}$, and therefore

$$\begin{aligned}
|||\hat{\sigma}(U+\tilde{U})-\hat{\sigma}(U)-\hat{\sigma}^{(1)}(U)\,\tilde{U}|||_{\gamma,0;T} &\lesssim \sup_{s\in[0,1]}|||(\hat{\sigma}^{(1)}(U_s)-\hat{\sigma}(U))\,\tilde{U}|||_{\gamma,0;T}\\
&\lesssim \sup_{s\in[0,1]}|||\hat{\sigma}^{(1)}(U_s)-\hat{\sigma}(U)|||_{\gamma,0;T}\,|||\tilde{U}|||_{\gamma,0;T}\\
&\lesssim C(|||U|||_{\gamma,0;T},|||U+\tilde{U}|||_{\gamma,0;T})\,|||\tilde{U}|||_{\gamma,0;T}^2,
\end{aligned}$$

which concludes the result. □

**Proof of Theorem 4.4.** We set the notation $\Delta\Phi_h := \Phi(\,\cdot\,;u_0+h) - \Phi(\,\cdot\,;u_0)$ and $R_h := \Delta\Phi_h - W_{u_0}[h]$. Then we have the identity

$$\begin{aligned} & \hat{\sigma}(\Phi(\,\cdot\,;u_0+h)) - \hat{\sigma}(\Phi(\,\cdot\,;u_0)) - \hat{\sigma}^{(1)}(\Phi(\,\cdot\,;u_0))\,W_{u_0}[h] \\ = \;& \hat{\sigma}(\Phi(\,\cdot\,;u_0)+\Delta\Phi_h) - \hat{\sigma}(\Phi(\,\cdot\,;u_0)) - \hat{\sigma}^{(1)}(\Phi(\,\cdot\,;u_0))\,\Delta\Phi_h + \sigma^{(1)}(\Phi(\,\cdot\,;u_0))\,R_h, \end{aligned}$$

and therefore the remainder $R_h \in \mathcal{D}_T^{\gamma,0}$ satisfies the equation

$$R_h = (\mathcal{K}_{\gamma-2} + R_\gamma\mathcal{R})\boldsymbol{R}^+((\hat{\sigma}(\Phi(\,\cdot\,;u_0)+\Delta\Phi_h) - \hat{\sigma}(\Phi(\,\cdot\,;u_0)) - \hat{\sigma}^{(1)}(\Phi(\,\cdot\,;u_0))\Delta\Phi_h)\Xi + \sigma^{(1)}(\Phi(\,\cdot\,;u_0))R_h\Xi).$$

Using analogous arguments to Lemma 4.5 we have the estimates

$$\begin{aligned} & |||(\mathcal{K}_{\gamma-2} + R_\gamma\mathcal{R})\boldsymbol{R}^+((\hat{\sigma}(\Phi(\,\cdot\,;u_0)+\Delta\Phi_h) - \hat{\sigma}(\Phi(\,\cdot\,;u_0)) - \hat{\sigma}^{(1)}(\Phi(\,\cdot\,;u_0))\,\Delta\Phi_h)\,\Xi)|||_{\gamma,0;T} \\ \lesssim\; & (T^{\alpha/2}+T)\,|||(\hat{\sigma}(\Phi(\,\cdot\,;u_0)+\Delta\Phi_h) - \hat{\sigma}(\Phi(\,\cdot\,;u_0)) - \hat{\sigma}^{(1)}(\Phi(\,\cdot\,;u_0))\,\Delta\Phi_h)\,\Xi|||_{\gamma-2,|\Xi|;T} \\ \lesssim\; & |||\hat{\sigma}(\Phi(\,\cdot\,;u_0)+\Delta\Phi_h) - \hat{\sigma}(\Phi(\,\cdot\,;u_0)) - \hat{\sigma}^{(1)}(\Phi(\,\cdot\,;u_0))\,\Delta\Phi_h|||_{\gamma-\alpha,0;T} \\ \lesssim\; & |||\hat{\sigma}(\Phi(\,\cdot\,;u_0)+\Delta\Phi_h) - \hat{\sigma}(\Phi(\,\cdot\,;u_0)) - \hat{\sigma}^{(1)}(\Phi(\,\cdot\,;u_0))\,\Delta\Phi_h|||_{\gamma,0;T} \\ \lesssim\; & C(|||\Phi(\,\cdot\,;u_0+h)|||_{\gamma,0;T})\,|||\Delta\Phi_h|||_{\gamma,0;T}^2, \end{aligned}$$

where the last estimate follows from Lemma 4.6, and

$$\begin{aligned} |||(\mathcal{K}_{\gamma-2} + R_\gamma\mathcal{R})\boldsymbol{R}^+(\sigma^{(1)}(\Phi(\,\cdot\,;u_0))\,R_h\,\Xi)|||_{\gamma,0;T} \;\lesssim\; & (T^{\alpha/2}+T)\,|||\sigma^{(1)}(\Phi(\,\cdot\,;u_0))\,R_h\,\Xi|||_{\gamma-2,|\Xi|;T} \\ \lesssim\; & T^{\alpha/2}\,|||\sigma^{(1)}(\Phi(\,\cdot\,;u_0))\,R_h|||_{\gamma-\alpha,0;T} \\ \lesssim\; & T^{\alpha/2}\,|||\sigma^{(1)}(\Phi(\,\cdot\,;u_0))\,R_h|||_{\gamma,0;T} \\ \lesssim\; & T^{\alpha/2}\,|||\sigma^{(1)}(\Phi(\,\cdot\,;u_0))|||_{\gamma,0;T}\,|||R_h|||_{\gamma,0;T} \\ \lesssim\; & T^{\alpha/2}\,|||R_h|||_{\gamma,0;T}. \end{aligned}$$

We conclude that for any $T_2 \in (0,T)$,

$$|||R_h|||_{\gamma,0;T_2} \;\lesssim\; |||\Delta\Phi_h|||_{\gamma-\alpha,0;T_2}^2 + C(|||\Phi(\,\cdot\,;u_0)|||_{\gamma,0;T})\,T_2^{\alpha/2}\,|||R_h|||_{\gamma,0;T_2},$$

and therefore choosing $T_2 \in (0,T_1]$ small enough, where $T_1$ is as in the proof of Lemma 4.5, one can conclude that

$$|||R_h|||_{\gamma,0;T_2} \lesssim |||\Delta\Phi_h|||_{\gamma,0;T_2}^2 \lesssim |||\Delta\Phi_h|||_{\gamma,0;T_1}^2 \lesssim \|h\|_{L^\infty(\mathbb{T})}^2.$$

Another bootstrapping argument allows us to conclude for the full interval $[0,T]$ the bound:

$$|||\Phi(\,\cdot\,;u_0+h) - \Phi(\,\cdot\,;u_0) - W_{u_0}[h]|||_{\gamma,0;T} \lesssim C^{\lceil T/T_2\rceil}\,\|h\|_{L^\infty(\mathbb{T})}^2,$$

where $T_2^{-1}$ depends polynomially on $|||\Phi(\,\cdot\,;u_0+h)|||_{\gamma,0;T}$, and which is uniformly bounded for $\|h\|_{L^\infty(\mathbb{T})} \leqslant 1$ by Lemma 4.7. The differentiability of $\varphi$ now follows, for example, by the reconstruction theorem (see also Lemma 4.8). The proof of (4.6) is done in Section 4.2. ☐

**LEMMA 4.7**. *We have* $|||\Phi(\,\cdot\,;u_0+h)|||_{\gamma,0;T} \lesssim 1$ *for all* $\|h\|_{L^\infty(\mathbb{T})} \leqslant 1$.

**Proof.** For simplicity we set $T=1$ and denote $u_h := \langle \mathbf{1}, \Phi(\,\cdot\,;u_0+h)\rangle$, and $u_h' = \langle \boldsymbol{X}, \Phi(\,\cdot\,;u_0+h)\rangle$. By Theorem 3.20 for small times the $L^\infty$-norm of the solution is bounded by $\|h\|_{L^\infty(\mathbb{T})} + 1 \lesssim 1$, and by Theorem 1.1 for large times we have the coming down from infinity bound which is independent of the initial condition, and in particular, we can conclude that $\|u_h\|_{[0,T]} \lesssim 1$ (with the proportionality constant depending on the fixed model), which controls the coefficient $\tau = \mathbf{1}$ in (4.3). For the coefficient $\tau = \boldsymbol{X}$ in (4.3) let $z \in (0,1] \times \mathbb{T}$ and let $\lambda_z := \lambda_0 \wedge \frac{1}{2}|z_0|^{\frac{1}{2}}$, where $\lambda_0$ is as in Definition 3.1. Then $\lambda_z^2 \leqslant \frac{1}{4} z_0 < z_0$ which implies that $\lambda_z \in (0,\lambda_0]$ and $z \in P_{\lambda_z}$ and therefore Corollary 3.3 applies as

$$\lambda_z\,|u_h'(z)| \leqslant \lambda_z\,\|u_h'\|_{B_{\frac{1}{2}\lambda_z}(z)} \lesssim \|u_h\|_{P_0}.$$

If $\lambda_z = \lambda_0$ then

$$|z_0|^{\frac{1}{2}} |u_h'(z)| \leqslant |z_0|^{\frac{1}{2}} \lambda_z^{-1} \|u_h\|_{P_0} = |z_0|^{\frac{1}{2}} \lambda_0^{-1} \|u_h\|_{P_0} \sim |z_0|^{\frac{1}{2}} \|u\|_{P_0}^E \lesssim \|u\|_{P_0}^E$$

for some $E > 0$. On the other hand, if $\lambda_z = \frac{1}{2} |z_0|^{\frac{1}{2}}$ then

$$|z_0|^{\frac{1}{2}} |u_h'(z)| \sim \lambda_z |u_h'(z)| \lesssim \|u_h\|_{P_0} \lesssim 1.$$

For non-polynomial coefficients we have that $\Upsilon[\tau] \sim (u_h')^{|\mathfrak{n}(\tau)|}$ and therefore, by the previous case,

$$\frac{|\langle \tau, \Phi( z ; u_0 + h) \rangle|}{|z_0|^{\frac{1}{2}(0 - |\tau|) \wedge 0}} \sim |u_h'(z)|^{|\mathfrak{n}(\tau)|} |z_0|^{\frac{1}{2}|\mathfrak{n}(\tau)|} |z_0|^{\frac{1}{2}(|\tau| - |\mathfrak{n}(\tau)|)} \lesssim |z_0|^{\frac{1}{2}|\mathfrak{n}(\tau)|} |u_h'(z)|^{|\mathfrak{n}(\tau)|} \lesssim 1$$

which shows that $\|\Phi( \cdot ; u_0 + h)\|_{\gamma, 0; T} \lesssim 1$. For the terms in (4.2) one argues analogously using Theorem 3.4 and Corollary 3.3 for the polynomial components, whilst for the non-polynomial coefficients one uses in addition the coherence in the form of Lemma 2.6 to write them in terms of $\sigma(\Phi( \cdot ; u_0 + h))$. □

## 4.2. Renormalised equation for the linearisation

In this section we prove the identity (4.6) from Theorem 4.4. First we show that the estimate (4.5) implies the differentiability of all the coefficients of the modelled distribution $\Phi( \cdot ; u_0)$, which in particular implies the differentiability of the flow $\varphi$ since $\varphi( \cdot ; u_0) = \mathcal{R}(\Phi( \cdot ; u_0)) = \langle \mathbf{1}, \Phi( \cdot ; u_0) \rangle$.

**Lemma 4.8**. *For any $x \in (0, T) \times \mathbb{T}$ and $\tau \in \mathcal{T}$, $(D_{u_0} \langle \tau, \Phi(x; u_0) \rangle)[h](x)$ is well defined and satisfies the identity $(D_{u_0} \langle \tau, \Phi( \cdot ; u_0) \rangle)[h] = \langle \tau, W_{u_0}[h] \rangle$.*

**Proof.** By definition of the seminorms Definition 4.1 and the estimate (4.5) from Theorem 4.4:

$$\begin{aligned} |\langle \tau, (\Phi( \cdot ; u_0 + h) - \Phi( \cdot ; u_0) - W_{u_0}[h])(x) \rangle| &\lesssim x_0^{\frac{1}{2}((0 - |\tau|) \wedge 0)} \|\Phi( \cdot ; u_0 + h) - \Phi( \cdot ; u_0) - W_{u_0}[h]\|_{\gamma, 0; T} \\ &\lesssim x_0^{-\frac{1}{2}|\tau|} |||\Phi( \cdot ; u_0 + h) - \Phi( \cdot ; u_0) - W_{u_0}[h]|||_{\gamma, 0; T} \\ &\lesssim x_0^{-\frac{1}{2}|\tau|} \|h\|_{L^\infty(\mathbb{T})}^2. \end{aligned}$$

In particular, for each $x$ we have that $\langle \tau, \Phi(x; u_0) \rangle$ is differentiable with respect to $u_0$ in the direction $h$ and we have the claimed identity. □

**Definition 4.9**. *Given $F: L^\infty(\mathbb{T}) \rightarrow \mathcal{D}_T^{\gamma, 0}$ such that $(D_{u_0} \langle \tau, F(x) \rangle)[h]$ is well-defined for any $\tau \in \mathcal{T}$ and $h \in L^\infty(\mathbb{T})$ we define*

$$(\hat{D}_{u_0} F)[h] := \sum_{\tau \in \mathcal{T}} \frac{(D_{u_0} \langle \tau, F \rangle)[h]}{\tau!} \tau, \qquad h \in L^\infty(\mathbb{T}).$$

**Remark 4.10.** In [CFG17, Sch23] they consider the lift of the Malliavin derivative to the space of modelled distributions and obtain results analogous to this. Since the derivative with respect to the initial condition that we consider vanishes on the model terms, this leads to significant simplifications with respect to those works.

The following result tells us that $\hat{D}_{u_0}$ indeed satisfies the Leibniz rule and the chain rule.

**LEMMA 4.11**. *We have that $\hat{D}^h_{u_0}$ satisfies Leibniz rule, i.e., for $F,G\colon L^\infty(\mathbb{T})\to\mathcal{D}^{\gamma,0}_T$ such that $\hat{D}^h_{u_0}$ is well-defined on them, and the product $FG\colon L^\infty(\mathbb{T})\to\mathcal{D}^{\gamma,0}_T$ is well-defined, we have the identity*

$$\hat{D}^h_{u_0}(FG)=\hat{D}^h_{u_0}(F)\,G+F\,\hat{D}^h_{u_0}(G).$$

*Moreover, if $F$ is function-like, and $\sigma\colon\mathbb{R}^{1+d}\to\mathbb{R}$ is sufficiently smooth to define $\hat{\sigma}(F)$, then*

$$\hat{D}^h_{u_0}(\hat{\sigma}(F))=\hat{\sigma}^{(1)}(F)\,\hat{D}^h_{u_0}(F).$$

**Proof.** Write $F(x)=\sum_\tau F^\tau(x)\,\tau$ and $G=\sum_\tau G^\tau(x)\,\tau$. Then $FG=\sum_{\tau_1,\tau_2}F^{\tau_1}G^{\tau_2}\tau_1\tau_2$ and therefore

$$\begin{aligned}
\hat{D}^h_{u_0}(FG) &= \sum_{\tau_1,\tau_2} D^h_{u_0}(F^{\tau_1}G^{\tau_2})\,\tau_1\,\tau_2\\
&= \sum_{\tau_1,\tau_2}((D^h_{u_0}F^{\tau_1})\,G^{\tau_2}+F^{\tau_1}\,(D^h_{u_0}G^{\tau_2}))\,\tau_1\,\tau_2\\
&= \sum_{\tau_1,\tau_2}(D^h_{u_0}F^{\tau_1})\tau_1\,G^{\tau_2}\,\tau_2+\sum_{\tau_1,\tau_2}F^{\tau_1}\,\tau_1\,(D^h_{u_0}G^{\tau_2})\,\tau_2\\
&= (\hat{D}^h_{u_0}F)\,G+F\,(\hat{D}^h_{u_0}G).
\end{aligned}$$

For the second part, we have that

$$\hat{\sigma}(F)=\sum_{k\geqslant 0}\frac{\sigma^{(k)}(f^{\mathbf{1}})}{k!}\,(F-f^{\mathbf{1}}\,\mathbf{1})^k$$

and therefore

$$\begin{aligned}
\hat{D}^h_{u_0}\,\hat{\sigma}(F) &= \sum_{k\geqslant 0}\left(\frac{D^h_{u_0}(\sigma^{(k)}(f^{\mathbf{1}}))}{k!}\,(F-f^{\mathbf{1}}\,\mathbf{1})^k+\frac{\sigma^{(k)}(f^{\mathbf{1}})}{k!}\,\hat{D}^h_{u_0}((F-f^{\mathbf{1}}\,\mathbf{1})^k)\right)\\
&= \sum_{k\geqslant 0}\left(\frac{\sigma^{(k+1)}(f^{\mathbf{1}})\,D^h_{u_0}f^{\mathbf{1}}}{k!}\,(F-f\,\mathbf{1})^k+\frac{\sigma^{(k)}(f^{\mathbf{1}})}{k!}\,k\,(F-f^{\mathbf{1}}\,\mathbf{1})^{k-1}\,\hat{D}^h_{u_0}(F-f^{\mathbf{1}}\,\mathbf{1})\right)\\
&= \left(\sum_{k\geqslant 0}\frac{\sigma^{(k+1)}(f^{\mathbf{1}})}{k!}(F-f^{\mathbf{1}}\mathbf{1})^k\right)D^h_{u_0}f^{\mathbf{1}}+\left(\sum_{k\geqslant 1}\frac{\sigma^{(k)}(f^{\mathbf{1}})}{(k-1)!}(F-f^{\mathbf{1}}\mathbf{1})^{k-1}\right)(\hat{D}^h_{u_0}F-(D^h_{u_0}f^{\mathbf{1}})\mathbf{1})\\
&= \hat{\sigma}^{(1)}(F)\,D^h_{u_0}\,f^{\mathbf{1}}+\left(\sum_{k\geqslant 0}\frac{\sigma^{(k+1)}(f^{\mathbf{1}})}{k!}\,(F-f^{\mathbf{1}}\,\mathbf{1})^k\right)(\hat{D}^h_{u_0}F-(D^h_{u_0}f^{\mathbf{1}})\,\mathbf{1})\\
&= \hat{\sigma}^{(1)}(F)\,D^h_{u_0}\,f^{\mathbf{1}}+\hat{\sigma}^{(1)}(F)\,(\hat{D}^h_{u_0}F-(D^h_{u_0}f^{\mathbf{1}})\,\mathbf{1})\\
&= \hat{\sigma}^{(1)}(F)\,\hat{D}^h_{u_0}F,
\end{aligned}$$

which concludes the proof. □

**LEMMA 4.12**. *Let $(\Pi,\Gamma)$ be a smooth model. Assume that $F\colon L^\infty(\mathbb{T})\to\mathcal{D}^{\gamma,0}_T$ is such that $\hat{D}_{u_0}F$ is well defined. Then, the reconstruction operator and the derivative with respect to the initial condition commute, i.e.,*

$$\mathcal{R}((\hat{D}_{u_0}F)[h])=(D_{u_0}(\mathcal{R}(F)))[h].$$

**Proof.** Since $(\Pi,\Gamma)$ is a smooth model, $\mathcal{R}F$ is a function and given by the evaluation at the diagonal, i.e.,

$$(\mathcal{R}F)(x)=\Pi_x(F(x))(x)=\sum_\tau F^\tau(x)\,(\Pi_x\tau)(x).$$

Therefore

$$\begin{aligned}(D_{u_0}(\mathcal{R}F))[h] &= \sum_{\tau\in\mathcal{T}} D_{u_0}(F^\tau(x)\,(\Pi_x\tau)(x))[h]\\ &= \sum_{\tau\in\mathcal{T}} ((D_{u_0}F^\tau(x))\,(\Pi_x\tau)(x) + F^\tau(x)\,D_{u_0}((\Pi_x\tau)(x)))[h]\\ &= \sum_{\tau\in\mathcal{T}} (D_{u_0}F^\tau(x))[h]\,(\Pi_x\tau)(x),\end{aligned}$$

since the model terms $\{\Pi_{\cdot}\tau\}_{\tau\in\mathcal{T}}$ are independent of the initial condition. On the other hand, since by definition, $(\hat{D}_{u_0}F)[h] = \sum_{\tau\in\mathcal{T}}(D_{u_0}F^\tau)[h]\,\tau$ then

$$\mathcal{R}((\hat{D}_{u_0}F)[h])(x) = \Pi_x((\hat{D}_{u_0}F(x))[h])(x) = \sum_\tau (D_{u_0}F^\tau)[h](x)\,(\Pi_x\tau)(x) = (D_{u_0}(\mathcal{R}F))[h],$$

which concludes the proof. □

**Proof of (4.6) from Theorem 4.4.** Using Lemma 4.11 we have that

$$\begin{aligned}(\hat{D}_{u_0}(\hat{\sigma}(\Phi(\,\cdot\,;u_0))\,\Xi))[h] &= ((\hat{D}_{u_0}\hat{\sigma}(\Phi(\,\cdot\,;u_0)))\,\Xi + \hat{\sigma}(\Phi(\,\cdot\,;u_0))\,\hat{D}_{u_0}\Xi)[h]\\ &= \hat{\sigma}^{(1)}(\Phi(\,\cdot\,;u_0))\,(\hat{D}_{u_0}\Phi(\,\cdot\,;u_0))\,\Xi\\ &= \hat{\sigma}^{(1)}(\Phi(\,\cdot\,;u_0))\,(W_{u_0}[h])\,\Xi,\end{aligned}$$

with the last identity following from Lemma 4.8 and Definition 4.9. By the choice of $\gamma > 2-\alpha$ the uniqueness of the reconstructions imply

$$\mathcal{R}(\hat{D}_{u_0}(\hat{\sigma}(U)\,\Xi)[h]) = \mathcal{R}(\hat{\sigma}^{(1)}(U)\,W_{u_0}[h]\,\Xi).$$

Using the commutation from Lemma 4.12 and applying the operator $(\partial_t - \partial_x^2)$ to the mild formulation (4.4) one obtains that

$$D_{u_0}(\mathcal{R}(\hat{b}(U) + \hat{\sigma}(U)\,\Xi))[h] = (\partial_t - \partial_x^2)\mathcal{R}(W_{u_0}[h]),$$

which concludes the proof. □

**Proof of Corollary 1.3.** Choose $\alpha = \frac{1}{2} - \kappa$, with $\kappa > 0$ sufficiently small, and then choose $\gamma > 2-\alpha$ sufficiently close to $2-\alpha$. The assumptions of Theorem 1.1 are then satisfied with $m = 3$, and $\lceil \frac{\gamma}{\alpha} \rceil = 4$. Let $Z = (\Pi, \Gamma)$ be the Itô model of [HP15]. On the event of probability one on which $Z$ is defined on every compact time interval, the fixed-point construction of Definition 4.2 can be restarted, and the a priori bound of Theorem 1.1 rules out finite-time explosion which yields the global flow. Its flow property follows from uniqueness of solutions. The argument of [HP15, Theorem 6.2] then identifies $\mathcal{R}(\hat{\sigma}(U)\,\Xi)$ with $\sigma(u)\,\xi$ in the Walsh–Itô sense. The same argument applied to $\hat{\sigma}^{(1)}(U)\,W_{u_0}[h]\,\Xi$ identifies its reconstruction with $\sigma'(u)\,w\;\xi$, where $w = \mathcal{R}(W_{u_0}[h])$ in the Walsh–Itô sense and shows the stated equation for the first variation. □

## 4.3. The stochastic Allen-Cahn equation

We now consider the particular case where $\xi \in \mathcal{C}^{\alpha-2}$ is an a.s. realisation of $(1+1)$ space-time white noise in $\mathbb{R}^+ \times \mathbb{T}$. This corresponds to $\alpha = 1/2 - \kappa \in (0,1)$ for some fixed $0 < \kappa \ll 1$. The lift of $\xi$ to a model $(\Pi, \Gamma)$, the Itô model, was done in [HP15] as a limit of smooth models. There are some small differences in the convention for the trees that generate the regularity structure used in [HP15], since noises there are seen as decorations on leaves instead of decorations on edges as we do. Following their notation, the noise is represented as $\Xi \leftrightarrow \circ$, (spatial) polynomials are denoted by $\boldsymbol{X} \leftrightarrow \boldsymbol{X}$ and when a noise is multiplied with this monomial it is denoted by $\boldsymbol{X}\,\Xi \longleftrightarrow \otimes$. Integration by $\mathcal{I}$ is denoted by a straight line, e.g., $\mathcal{I}(\Xi) \leftrightarrow$ ⚲. With this notation in mind, we have for $\gamma \in (2-\alpha, 2) = (3/2 + \kappa, 2)$ close enough to $2 - \alpha = 3/2 + \kappa$ that the elements in $\mathcal{T}_\gamma$ are

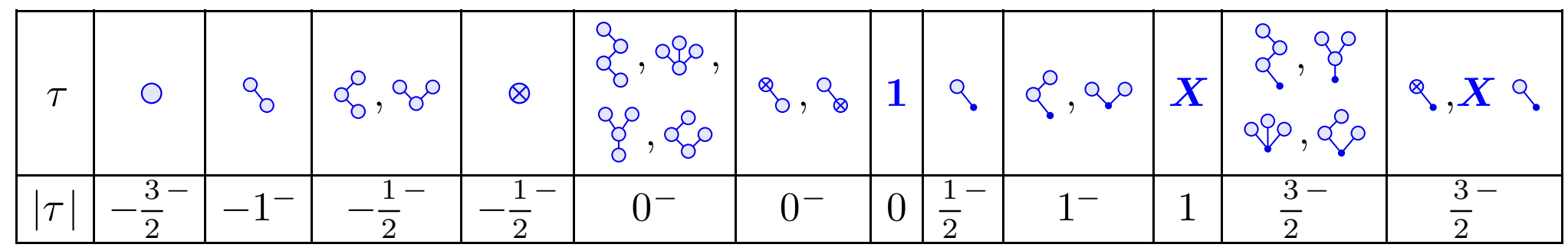

| $\tau$ | ○ | | | ⊗ | | | $\mathbf{1}$ | | | $\boldsymbol{X}$ | | , $\boldsymbol{X}$ |
|---|---|---|---|---|---|---|---|---|---|---|---|---|
| $\lvert\tau\rvert$ | $-\frac{3}{2}^-$ | $-1^-$ | $-\frac{1}{2}^-$ | $-\frac{1}{2}^-$ | $0^-$ | $0^-$ | $0$ | $\frac{1}{2}^-$ | $1^-$ | $1$ | $\frac{3}{2}^-$ | $\frac{3}{2}^-$ |

**Table 4.1.**

On the other hand, the elements in the structure group $\mathcal{T}_\gamma^+$ are where integration by $\mathcal{I}_{e_1} = \mathcal{I}'$ (corresponding to a spatial derivative of the heat kernel) is denoted by a zig-zag line, e.g., $\mathcal{I}'(\boldsymbol{X}\,\Xi) \leftrightarrow$ .

| $\mu$ | $\mathbf{1}$ | | | | $\boldsymbol{X}$ | , , , , $\boldsymbol{X}$ , $\boldsymbol{X}$ | , $\boldsymbol{X}$ , $\boldsymbol{X}$ |
|---|---|---|---|---|---|---|---|
| $\lvert\mu\rvert$ | $0$ | $\frac{1}{2}^-$ | $\frac{1}{2}^-$ | $1^-$ | $1$ | $\frac{3}{2}^-$ | $\frac{3}{2}^-$ |

**Table 4.2.**

The representation of a coherent modelled distribution $U := U_\gamma$ as in (2.6) is given explicitly by

$$U = u\,\mathbf{1} + \sigma(u) \; + (\sigma\,\sigma^{(1)})(u) \; + u'\,\boldsymbol{X} + (\sigma\,(\sigma^{(1)})^2)(u) \; + \frac{(\sigma^2\,\sigma^{(2)})(u)}{2} \; + \sigma^{(1)}(u)\,u' \; .$$

To reconstruct the singular product $\sigma(u)\,\xi$ one needs to model $\hat{\sigma}(U)$ ○ to degree $\gamma + |\Xi| = 0^+$ which, by coherence, has the form

$$\begin{aligned}
\hat{\sigma}(U) \circ \;\; = \;\; & \sigma(u) \circ + (\sigma\,\sigma^{(1)})(u) \; + (\sigma\,(\sigma^{(1)})^2)(u) \; + \frac{\sigma^2\,\sigma^{(2)}(u)}{2} \; + \sigma^{(1)}(u)\,u' \otimes \\
& + (\sigma\,(\sigma^{(1)})^3)(u) \; + \frac{(\sigma^3\,\sigma^{(3)})(u)}{6} \; + \frac{(\sigma^2\,\sigma^{(1)}\,\sigma^{(2)})(u)}{2} \\
& + (\sigma^2\,\sigma^{(1)}\,\sigma^{(2)})(u) \; + (\sigma^{(1)})^2(u)\,u' \; + (\sigma\,\sigma^{(2)})(u)\,u' \; .
\end{aligned} \tag{4.8}$$

For $0 < \varepsilon \ll 1$ we consider a smooth mollification of the noise $\xi_\varepsilon$ and let $\{(\Pi^\varepsilon, \Gamma^\varepsilon)\}_{\varepsilon>0}$ be the smooth renormalised models converging to the Itô model $(\Pi, \Gamma)$ as constructed in [HP15]. Denote by $\varphi_\varepsilon$ the flow (see Definition 4.2) associated to the renormalised equation

$$(\partial_t - \partial_x^2)\,u_\varepsilon = u_\varepsilon - u_\varepsilon^3 + \mathcal{R}_\varepsilon(\hat{\sigma}(U_\varepsilon) \circ). \tag{4.9}$$

Since $(\Pi_\varepsilon, \Gamma_\varepsilon)$ is a function-valued model, we can recover the reconstruction operator associated to it by evaluating the model at the diagonal. Following [HP15, Proposition 4.4] the only non-zero diagonal terms of the model appearing in (4.8), and which therefore contribute to the reconstruction operator and the renormalised equation, come from the trees $\{$ ○, , , $\}$ for which one has that

$$\Pi_x^\varepsilon(\circ)(x) = \xi_\varepsilon(x), \;\; \Pi_x^\varepsilon(\;)(x) = -c_\varepsilon, \;\; \Pi_x^\varepsilon\Big(\;\Big)(x) = -c^{(1)}, \;\; \Pi_x^\varepsilon\Big(\;\Big)(x) = -c^{(2)}. \tag{4.10}$$

It follows that

$$\mathcal{R}_\varepsilon(\hat{\sigma}(U_\varepsilon) \circ) = \sigma(u_\varepsilon(x))\xi_\varepsilon(x) - c_\varepsilon\,(\sigma\,\sigma^{(1)})(u_\varepsilon(x)) - c^{(1)}\,(\sigma(\sigma^{(1)})^3)(u_\varepsilon(x)) - c^{(2)}\,(\sigma^2\,\sigma^{(1)}\,\sigma^{(2)})(u_\varepsilon(x))$$

and therefore (4.9) can be written explicitly as

$$(\partial_t - \partial_x^2)\,u_\varepsilon = u_\varepsilon - u_\varepsilon^3 + \sigma(u_\varepsilon)\xi_\varepsilon - c_\varepsilon\,(\sigma\,\sigma^{(1)})(u_\varepsilon) - c^{(1)}\,(\sigma(\sigma^{(1)})^3)(u_\varepsilon) - c^{(2)}\,(\sigma^2\,\sigma^{(1)}\,\sigma^{(2)})(u_\varepsilon). \tag{4.11}$$

By formally differentiating this we see that the linearisation of the flow in a direction $h \in L^\infty(\mathbb{T})$ is the solution to the renormalised linear equation

$$
\begin{aligned}
(\partial_t - \partial_x^2)(w_\varepsilon[h]) &= w_\varepsilon[h] - 3\,u_\varepsilon^2\, w_\varepsilon[h] + \sigma^{(1)}(u_\varepsilon)\,\xi_\varepsilon\, w_\varepsilon[h] - c_\varepsilon\,(\sigma\,\sigma^{(1)})^{(1)}(u_\varepsilon)\, w_\varepsilon[h] \\
&\quad - c^{(1)}(\sigma(\sigma^{(1)})^3)^{(1)}(u_\varepsilon) w_\varepsilon[h] - c^{(2)}(\sigma^2\sigma^{(1)}\sigma^{(2)})^{(1)}(u_\varepsilon) w_\varepsilon[h],
\end{aligned} \tag{4.12}
$$

with initial condition $w_\varepsilon[h](0) = h$. Observe that in the limit $\varepsilon \to 0$ the product $\sigma^{(1)}(u)\,\xi$ appearing in (4.12) is singular and needs to be renormalised. Moreover, this renormalised product will inherit the same low regularity as the space-time white noise $\xi$, which will make the product with $w_\varepsilon[h]$ in (4.12) singular and in need of another renormalisation. However, by Theorem 4.4 one can write the right hand side of (4.12) as the reconstruction of a modelled distribution, and this representation is stable in the limit.

To understand the renormalisation in (4.12) better we recall that by [FH20, Section 14.2] for any $F \in C^4(\mathbb{R}, \mathbb{R})$, $F(u)$ can be modelled as

$$
\begin{aligned}
\hat{F}(U) &= F(u)\mathbf{1} + F^{(1)}(u)\sigma(u)\,\text{(tree)} + F^{(1)}(u)(\sigma\sigma^{(1)})(u)\,\text{(tree)} + \frac{1}{2}F^{(2)}(u)\sigma^2(u)\,\text{(tree)} + F^{(1)}(u)u'\boldsymbol{X} \\
&\quad + F^{(1)}(u)(\sigma(\sigma^{(1)})^2)(u)\,\text{(tree)} + \frac{1}{3!}F^{(3)}(u)\sigma^3(u)\,\text{(tree)} + F^{(1)}(u)\frac{(\sigma^2\sigma^{(2)})(u)}{2}\,\text{(tree)} \\
&\quad + F^{(2)}(u)\frac{(\sigma^2\sigma^{(1)})(u)}{2}\,\text{(tree)} + F^{(2)}(u)\sigma(u)u'\boldsymbol{X}\,\text{(tree)} + F^{(1)}(u)\sigma^{(1)}(u)u'\,\text{(tree)}.
\end{aligned}
$$

In particular, for $F = \sigma$, after multiplication with the noise, we recover (4.8), whilst for $F = \sigma^{(1)}$ we obtain that

$$
\begin{aligned}
\hat{\sigma}^{(1)}(U)\,\circ &= \sigma^{(1)}(u)\,\circ + (\sigma\,\sigma^{(2)})(u)\,\text{(tree)} + (\sigma\,\sigma^{(1)}\,\sigma^{(2)})(u)\,\text{(tree)} + \frac{(\sigma^2\,\sigma^{(3)})(u)}{2}\,\text{(tree)} + \sigma^{(2)}(u)\,u'\,\otimes \\
&\quad + (\sigma\,(\sigma^{(1)})^2\,\sigma^{(2)})(u)\,\text{(tree)} + \frac{(\sigma^3\,\sigma^{(4)})(u)}{6}\,\text{(tree)} + \frac{(\sigma^2\,(\sigma^{(2)})^2)(u)}{2}\,\text{(tree)} \\
&\quad + (\sigma^2\,\sigma^{(1)}\,\sigma^{(3)})(u)\,\text{(tree)} + (\sigma\,\sigma^{(3)})(u)\,u'\,\text{(tree)} + (\sigma^{(1)}\,\sigma^{(2)})(u)\,u'\,\text{(tree)},
\end{aligned}
$$

and when combined with (4.10) we can conclude that

$$
\begin{aligned}
\mathcal{R}_\varepsilon(\hat{\sigma}^{(1)}(U_\varepsilon)\,\circ)(x) &= \sigma^{(1)}(u_\varepsilon(x))\xi_\varepsilon(x) - c_\varepsilon\,(\sigma\,\sigma^{(2)})(u_\varepsilon(x)) \\
&\quad - c^{(1)}\,(\sigma\,(\sigma^{(1)})^2\,\sigma^{(2)})(u_\varepsilon(x)) - c^{(2)}\,(\sigma^2\,\sigma^{(1)}\,\sigma^{(3)})(u_\varepsilon(x)).
\end{aligned} \tag{4.13}
$$

Since, by Leibniz rule, we have the identities

$$
\begin{aligned}
(\sigma\,\sigma^{(1)})^{(1)} &= \sigma\,\sigma^{(2)} + (\sigma^{(1)})^2 \\
(\sigma\,(\sigma^{(1)})^3)^{(1)} &= \sigma\,(\sigma^{(1)})^2\,\sigma^{(2)} + (\sigma\,(\sigma^{(1)})^2)^{(1)}\,\sigma^{(1)} \\
(\sigma^2\,\sigma^{(1)}\,\sigma^{(2)})^{(1)} &= \sigma^2\,\sigma^{(1)}\,\sigma^{(3)} + (\sigma^2\,\sigma^{(1)})^{(1)}\,\sigma^{(2)},
\end{aligned}
$$

the renormalisation needed to define the product $\sigma^{(1)}(u)\,\xi$ is already included in (4.12), and the other terms are precisely there to renormalise the product with $w_\varepsilon[h]$.

**Remark 4.13.** One could view the linearised equation (4.12) as a singular linear PAM driven by the already renormalised noise $\mathcal{R}_\varepsilon(\hat{\sigma}^{(1)}(U)\,\Xi)$ as defined in (4.13). However, this approach would require a separate explicit renormalisation of the associated model and the derivation of its corresponding renormalisation constants.

**Acknowledgments.** Both authors are funded by the European Research Council (ERC) under the European Union's Horizon 2020 research and innovation programme (Grant agreement No. 101045082), and by the Deutsche Forschungsgemeinschaft (DFG, German Research Foundation) under Germany's Excellence Strategy EXC 2044-390685587, Mathematics Münster: Dynamics-Geometry-Structure.